\documentclass[12pt,oneside]{book}
\usepackage{setspace}
\usepackage{latexsym}
\usepackage{epsfig}
\usepackage{amsmath}
\usepackage{amsthm}
\usepackage{amssymb}
\usepackage{parskip}
\usepackage{hyperref}
\usepackage{color}
\usepackage{graphicx}

\usepackage{fancyhdr}
\fancyheadoffset{0in}

\def\cP{{\cal P}}

\newcommand{\beq}[1]{\begin{equation}\label{#1}}
\newcommand{\eeq}{\end{equation}}

\def\a{\alpha} \def\b{\beta} \def\d{\delta} \def\D{\Delta}
\def\e{\epsilon}    \def\g{\gamma}
\def\G{\Gamma}  
\def\z{\zeta} \def\th{\theta}  \def\Th{\Theta}  \def\l{\lambda}
  \def\n{\nu} \def\p{\pi}
\def\r{\rho}  \def\s{\sigma} 
\def\t{\tau} \def\om{\omega} \def\OM{\Omega} \def\Om{\Omega}

\newtheorem{theorem}{Theorem}

\newtheorem{conjecture}{Conjecture}
\newtheorem{lemma}[theorem]{Lemma}
\newtheorem{corollary}[theorem]{Corollary}

\newtheorem{definition}{Definition}
\newtheorem{proposition}{Proposition}
\newcommand{\proofstart}{{\bf Proof\hspace{2em}}}
\newcommand{\proofend}{\hspace*{\fill}\mbox{$\Box$}}
\def\cW{{\cal W}}

\newcommand{\ooi}{(1+o(1))}

\newcommand{\ul}[1]{\mbox{\boldmath$#1$}}

\newcommand{\wh}[1]{\widehat{#1}}

\newcommand{\rdown}[1]{{\lfloor #1 \rfloor}}

\newcommand{\brac}[1]{\left(#1\right)}

\newcommand{\bfrac}[2]{\left(\frac{#1}{#2}\right)}

\newcommand{\rai}{\rightarrow \infty}
\newcommand{\ra}{\rightarrow}
\newcommand{\rat}{{\textstyle \ra}}

\def\E{\mbox{{\bf E}}}
\def\Pr{\mbox{{\bf Pr}}}
\def\whp{{\bf whp}}

\newcommand{\ignore}[1]{}

\newcommand{\Gvc}{\G_v^{\circ}}

\newcommand{\Gv}{\G_v}
\def\uar{{\bf uar}} 

\def\cov{\mbox{{\bf COV}}}
\def\hit{\mbox{{\bf H}}}
\def\com{\mbox{{\bf COM}}}
\def\blank{\mbox{{\bf B}}}
\def\bcov{\mbox{{\bf BCOV}}}
\def\fretT{\mbox{{\bf R}}}

\begin{document}

\begin{titlepage}
\begin{center}
{\qquad}\\[2.0cm]
\rule{\linewidth}{0.5mm}\\[0.4cm] 
{ \Large\bfseries The Cover Time of Random Walks on Graphs}\\[0.4cm]
\rule{\linewidth}{0.5mm}\\[1.5cm] 

Mohammed Abdullah\\[0.5cm]
{\large King's College London}\\[1.5cm]

\vfill
{\small Submitted for the degree Doctor of Philosophy}\\[0.5cm]
{\large September 2011}

\end{center}

\end{titlepage}

\chapter*{\centering \begin{normalsize}Abstract\end{normalsize}}
\addcontentsline{toc}{chapter}{Abstract}
{\small A simple random walk on a graph is a sequence of movements from one vertex to another where at each step an edge is chosen uniformly at random from the set of edges incident on the current vertex, and then transitioned to next vertex. Central to this thesis is the cover time of the walk, that is, the expectation of the number of steps required to visit every vertex, maximised over all starting vertices. 
In our first contribution, we establish a relation between the cover times of a pair of graphs, and the cover time of their Cartesian product. This extends previous work on special cases of the Cartesian product, in particular, the square of a graph. We show that when one of the factors is in some sense larger than the other, its cover time dominates, and can become within a logarithmic factor of the cover time of the product as a whole. Our main theorem effectively gives conditions for when this holds. \ignore{The probabilistic technique that we introduce is more general and may be of independent interest, as might some of our lemmas.}The techniques and lemmas we introduce may be of independent interest.
In our second contribution, we determine the precise asymptotic value of the cover time of a random graph with given degree sequence. This is a graph picked uniformly at random from all simple graphs with that degree sequence. We also show that with high probability, a structural property of the graph called conductance, is bounded below by a constant. This is of independent interest.
Finally, we explore random walks with weighted random edge choices. We present a weighting scheme that has a smaller worst case cover time than a simple random walk. We give an upper bound for a random graph of given degree sequence weighted according to our scheme. We demonstrate that the speed-up (that is, the ratio of cover times) over a simple random walk can be unbounded.}
\clearpage
\chapter*{\centering \begin{normalsize}Acknowledgment\end{normalsize}}
\addcontentsline{toc}{chapter}{Acknowledgment }
{\small I firstly wish to express my deepest gratitude to my supervisor, Colin Cooper, whom I have been very fortunate to have known. It has been a pleasure to work with Colin, both as his student and as a research colleague. His guidance, patience and encouragement have been invaluable, and I am greatly indebted to him for the opportunities he has given me.

I also wish to thank my second supervisor, Tomasz Radzik. More than merely an excellent source of advice, Tomasz has been a colleague with whom I have greatly enjoyed working. Our research meetings have always been inspiring and productive, and a rich source of ideas. 

I wish to thank Alan Frieze, a co-author of one of my papers that forms a significant part of this thesis. I am grateful to Alan for his role in directly developing the field to which I have dedicated so much time, and for the opportunity to collaborate with him.

In the final year of my time as a PhD student, I have been fortunate to have met and worked with Moez Draief. Though our work together does not form part of this thesis, it has nevertheless been both highly compelling and enjoyable part of my time as a PhD student. I would like to thank Moez for our research collaboration.

Finally, I would like to thank my parents for innumerable reasons, but in particular for their encouragement and support in all its forms. It is to them that this thesis is dedicated.} 
\clearpage

\tableofcontents
\listoffigures

\chapter{Introduction}
Let $G=(V,E)$ be a finite, connected, undirected graph. Suppose we start at time step $t=0$ on some vertex $u \in V$ and choose an edge $e$ uniformly at random (\uar) from those incident on $u$. We then transition $e$ to the vertex that the other end of $e$ is incident on. We repeat this process at the next step, and so on. This is known as a \emph{simple random walk} (often abbreviated to \emph{random walk}) on $G$. We shall denote it by $\mathcal{W}_u$, where the subscript is the starting vertex. We write $\mathcal{W}_u(t)=x$ if the walk is at vertex $x$ at time step $t$. 

Immediately, a number of questions can be asked about this process. For example, 
\begin{description}
	\item[(1)] 	Does $\mathcal{W}_u$ visit every vertex in $G$?
	\item[(2)]	If so, how long does it take on average?
	\item[(3)] 	On average, how long does it take to visit a particular vertex $v$?
	\item[(4)] 	On average, how long does it take to come back to itself?
	\item[(5)] 	In the long run, do all vertices get visited roughly the same number of times, or are there differences?
	\item[(6)]	If there are differences, what is the proportion of the time spent at a particular vertex $v$ in the long run?
	\item[(7)] 	How do the answers to the above questions vary if we change the starting vertex?
	\item[(8)] 	How do the answers to the above questions vary for a different graph?
\end{description}

This thesis addresses all of these questions in one way or another for specific classes of graphs. However, the particular question that is the central motivation for this work is the following:

\emph{For a random walk $\mathcal{W}_u$ on a simple, connected, undirected graph $G=(V,E)$, what is the expected number of steps required
to visit all the vertices in $G$, maximised over starting vertices $u$?}

The following quantities, related to the questions above, are formally defined in chapter \ref{intro}. The expected time it takes $\mathcal{W}_u$ to visit every vertex of $G$ is the \emph{cover time from $u$}, $\cov_u[G]$, and the \emph{cover time} $\cov[G]=\max_{u}\cov_u[G]$. The expected time it takes $\mathcal{W}_u$ to visit some $v$ is the \emph{hitting time} $\hit[u,v]$, and when $v=u$, it is called the \emph{first return time}. 

These questions, much like the process itself, are easy to understand, yet they and many others have been been the focus of a great deal of study in the mathematics and computer science communities.  Some questions are easy to answer with basic probability theory, others are more involved and seem to require more sophisticated techniques. The difficulty usually varies according to what kind of answer we are looking for. For example, for the $2$-dimensional torus $\mathbb{Z}^2_n$ with $N=n^2$ vertices, $\cov[\mathbb{Z}^2_n]=O(N\log^2N)$ is not too difficult to show with some of the theory and techniques we present in chapters \ref{networktheory} and \ref{Techniques and Results for Hitting and Cover Times}. However, it was not until quite recently that a precise asymptotic result of $\cov[\mathbb{Z}^2_n]\sim \frac{N}{\pi}\log^2N$ was given by \cite{Dembo}. 

This thesis is concerned primarily with cover time.

\section{Applications}
Before we give an outline of the thesis, we mention the applications of random walks and in particular, cover times. Applications are not the focus of this thesis, but it is worth mentioning their role, particularly in algorithmic and networking areas. The classical application of random walks in an algorithmic context is a randomised $s-t$ connectivity algorithm. The problem, known as the \emph{s-t connectivity problem}, is as follows: Given a graph $G=(V,E)$, with $|V|=n$ and $|E|=m$, and two vertices $s,t\in V$, if there is a path in $G$ connecting $s$ and $t$, return ``true'' otherwise, return ``false''. This can be done in time $O(n+m)$ with, for example, breadth first search. However, the space requirement is $\Omega(n)$ for such an algorithm (or various others, such as depth first search). Take, for example, the case where $G$ is a path of length $n$ and $s$ and $t$ are ends of the path.  

With a random-walk based algorithm, we can present a randomized algorithm for the problem that requires $O(\log n)$ space. It relies upon the following proposition, proved in chapter \ref{Techniques and Results for Hitting and Cover Times}. See, e.g., \cite{mitz}.
\\

\begin{proposition}
For any connected, finite graph $G=(V,E)$, $\cov[G]<4|V||E|$. \label{stprop}
\end{proposition}

To avoid confusion with the name of the problem, we shall use the variable $\tau$ to stand for time in the walk process. 

The algorithm is as follows: Start a walk $\mathcal{W}_s$ on $G$ from vertex $s$. Assume time $\tau=0$ at the start of the walk. Stop the walk at time $\tau$ if (i)$\tau=8nm$ or (ii) $\mathcal{W}_s(\tau)=t$. If $\mathcal{W}_s(\tau)=t$, then output ``true''. Otherwise output  ``false''.

Observe, if the algorithm returns ``true'', then it is correct, since there must be a path. If it returns ``false'', then it may be wrong, since it may have simply not visited $t$ even though it could have. The question is, what is the probability that the algorithm returns an incorrect answer given that an $s-t$ path exists? Suppose the random variable $X$ counts the number of steps the random walk from $s$ takes before visiting $t$. By Markov's inequality, 
\[
\Pr(X>T)\leq \frac{\E[X]}{T}.
\]
Now, $\E[X] = \hit[s,t]$ and $\hit[s,t]\leq \cov[G]$. So using Proposition \ref{stprop}, $\Pr(X>T)\leq 4nm/T$, hence $\Pr(X>8nm)\leq 1/2$.

The size of the input, that is, the graph, is $n$ vertices and $m\leq n^2$ edges. This algorithm needs $O(\log n)$ space since it need only store enough bits to keep track of its position and maintain a counter $\tau$.  

In fact, a breakthrough 2004 paper \cite{Omer} showed that the problem can be solved by a deterministic algorithm using $O(\log n)$ space. Nevertheless, the randomised algorithm demonstrated here, as well as being very simple, remains a strong example of the role that the theory of random walks plays in applications.

There are many other applications of random walks, particularly in networks and distributed systems, where they have been applied to self-stabilization (\cite{Dolev}), sensor networks (\cite{Servetto}), peer-to-peer networks (\cite{Gkan}, \cite{cooperp2p}) and voting (\cite{coopermult}), amongst many others.

\section{Overview of the Thesis}
Roughly speaking, this thesis can be divided into two parts. Aside from this chapter, chapters \ref{intro}, \ref{Theory of Markov Chains and Random Walks}, \ref{networktheory} and \ref{Techniques and Results for Hitting and Cover Times} are drawn from the established literature. They provide a review of some results and vital theory required for the original contribution. The original contribution can be considered to be in chapters \ref{ch:The Cover Time of Cartesian Product Graphs},  \ref{CovDS} and \ref{Weighted Random Walks}.

\subsection{Background}
In chapter \ref{intro}, we give definitions and basic lemmas for graphs, random walks and Markov chains. We also give definitions of \emph{weighted} random walks. These differ from simple random walks in how the next edge to be transitioned is chosen. Rather than choosing edges \uar\, the probabilities are weighted by a weight assigned to the edge. This is more general than simple random walks, which are weighted random walks in which all edge weights are the same. Weighted random walks are the subject of chapter \ref{Weighted Random Walks}. 

Simple random walks on graphs are special cases of weighted random walks, which are in turn, special cases of Markov chains. Markov chains play a part in a fundamental lemma of chapter \ref{ch:The Cover Time of Cartesian Product Graphs}, and we use theorems from the literature based on Markov chains in chapters \ref{CovDS} and \ref{Weighted Random Walks}. In chapter \ref{Theory of Markov Chains and Random Walks}, we give an account of the theory of Markov chains and random walks relevant to our work. Much of the theory presented in the framework of Markov chains is vital to sections of the thesis where only simple random walks are considered, and would have had to be written in a similar form had the more general presentation not been given. This, in conjunction with our use of Markov chains in various parts is why we chose to give a presentation in terms of Markov chains supplemented with explanations of how the general theory specialises for random walks.  We also demonstrate a characterisation of Markov chains that are equivalent to weighted random walks.

In chapter \ref{networktheory} we present the electrical network metaphor of random walks on graphs. This is a framework in which a theory has been built up to describe properties and behaviours of random walks in a different way. It provides a means of developing an intuition about random walks, and provides a tool kit of useful lemmas and theorems. Much of chapter \ref{ch:The Cover Time of Cartesian Product Graphs} is built on the material in chapter \ref{networktheory}. In chapter \ref{CovDS}, the tools of electrical network theory are exploited in a number of proofs.

In chapter \ref{Techniques and Results for Hitting and Cover Times} we present detailed proofs of hitting and cover times for specific, simple graph structures. We then give general techniques for bounding these parameters, including some we use in our proofs. We then give some bounds from the literature.

\subsection{Original Contribution}
In chapter \ref{ch:The Cover Time of Cartesian Product Graphs} we present the first section of our original contribution. We study random walks on the Cartesian product $F$ of a pair of graphs $G$ and $H$. We refer the reader to chapter \ref{ch:The Cover Time of Cartesian Product Graphs} for a definition of the \emph{Cartesian product}. After giving definitions and context (including related work from the literature), we describe a probabilistic technique which we use to analyse the cover time. We then present a number of lemmas relating to the effective resistance of products of graphs. We apply the above to the problem of the cover time, and develop bounds on the cover time of the product $F$ in terms of properties of the factors $G$ and $H$. The resulting theorem can be used to demonstrate when the cover time of one of the factors dominates the other, and becomes of the same order or within a logarithmic factor of the cover time of the product as a whole. The probabilistic technique we introduce and the effective resistance lemmas may be of independent interest. This chapter is based on joint work with Colin Cooper and Tomasz Radzik, published in \cite{CartProd}.

In chapter \ref{CovDS}, we give a precise asymptotic result for the cover time of a random graph with given degree sequence $\mathbf{d}$. That is, if a graph $G$ on $n$ vertices is picked uniformly at random from the set $\mathcal{G}(\mathbf{d})$ of all simple graphs with vertices having pre-specified degrees $\mathbf{d}$, then with high probability (\whp), the cover time tends to a value $\tau$ that we present. The phrase \emph{with high probability} means with probability tending to $1$ as $n$ tends to $\infty$. After giving an account of the necessary theory, we give a proof that a certain structural property of graphs, known as the \emph{conductance}, is, \whp, bounded below by a constant for a $G$ chosen \uar\ from $\mathcal{G}(\mathbf{d})$. This allows us to use some powerful theory from the literature to analyse the problem. We continue our analysis with further study of the structural properties of the graphs, and the behavior of random walks on them. We use these results to show that \whp, no vertex is unvisited by  time $\tau+\epsilon$, where $\epsilon$ is some quantity that tends to $0$. For the lower bound, we show that at time $\tau-\epsilon$, there is at least one unvisited vertex, \whp. This chapter is based on joint work with Colin Cooper and Alan Frieze, published in \cite{CovDS}.

Finally, in chapter \ref{Weighted Random Walks}, we investigate weighted random walks. Graph edges are given non-negative weights and the probability that an edge $e$ is transitioned from a vertex $u$ is $w(e)/w(u)$, where $w(e)$ is the weight of $e$ and $w(u)$ is the total weight of edges incident on $u$. We present from the literature a weighting scheme that has a worst case cover time better than a simple random walk. We then present our own weighting scheme, and show that it also has this property. We give an upper bound for the cover time of a weighted random walk on a random graph of given degree sequence weighted according to our scheme. We demonstrate that the speed-up (that is, the ratio of cover times) over a simple random walk can be unbounded. This chapter is based on joint work with Colin Cooper.

\chapter{Definitions and Notation}\label{intro}

\section{Graphs}\label{graphdefns}
A \emph{graph} $G$ is a tuple $(V,E)$, where the \emph{vertex set} $V$ is a set of objects called \emph{vertices} and the \emph{edge set} $E$ is a set of two-element tuples $(u,v)$ or two-element sets $\{u,v\}$ on members of $V$. The members of $E$ are called \emph{edges}. Graphs can be \emph{directed} or \emph{undirected}. In a directed graph the tuple $(u,v)$ is considered ordered, so $(u,v)$ and $(v,u)$ are two different edges. In an undirected graph, in line with the conventions of set notation, the edge $\{u,v\}$ can be written as $\{v,u\}$. However, in the standard conventions of the literature, edges of undirected graphs are usually written in tuple form, and tuples are considered unordered. Thus, the edge $\{u,v\}$ is written $(u,v)$, and is included only once in $E$. We shall use this convention throughout most of the thesis, and be explicit when departing from it. Furthermore, we may use the notation $u \in G$ and $e \in G$ to stand for $u \in V$ and $e\in E$ respectively.

In this thesis, we deal only with finite graphs, that is, both $V$ and $E$ are finite, and we shall sometimes use the notation $V(G)$ and $E(G)$ to denote the vertex and edge set respectively of a graph $G$.

A \emph{loop} $(u,u)$ is an edge from a vertex $u$ to itself. In a \emph{multigraph}, a pair of vertices $u,v$ can have more than one edge between them, and each edge is included once in $E$. In this case, $E$ is a multiset. A graph is \emph{simple} if it does not contain loops and is not a multigraph.

For a vertex $u \in V$, denote by $N(u)\subseteq V$ the \emph{neighbour set} of $u$, 
\begin{equation}
N(u) = \{v \in V : (u,v)\in E\}.
\end{equation} 
Denote by $d(u)$ or $d_u$ the \emph{degree} of $u$. This is the number of ends of edges incident on $u$, i.e., 
\begin{equation}
d(u) = |\{e \in E : e=(u,x), x \neq u\}| + 2|\{e \in E : e = (u,u)\}|. \label{nonSimpGraphDegree}
\end{equation} 

The second term in the sum shows that we count loops twice, since a loop has two ends incident on $u$. When $G$ is simple it is seen that $d(u) = |N(u)|$.

When a graph is directed, there is an \emph{in-degree} and an \emph{out-degree}, taking on the obvious definitions. 
  
If $|V|=n$ then $G$ can be represented as a $n\times n$ matrix $\mathbf{A}=[a_{i,j}]$, called the \emph{adjacency matrix}. Without loss of generality, assume that the vertices are labelled $1$ to $n$, then in $A$, $a_{i,i}=2l$ where $l$ is the number of loops from $i$ to itself, and for $i\neq j$, $a_{i,j}$ is the number of edges between $i$ and $j$.

A \emph{walk} in a graph $G$ is a sequence of (not necessarily distinct) vertices in $G$,  $(v_0, v_1, v_2, \ldots)$ or $(v_0, v_1, v_2, \ldots, v_t)$ if the sequence is finite. A vertex $v_i$ is followed by $v_{i+1}$ only if $(v_i, v_{i+1}) \in E$. A walk is a \emph{path} if and only if no vertex appears more than once in the sequence. If $v_0=v_t$ and this is the only vertex that repeats then then the walk is a \emph{cycle}.

If $w = (v_0, v_1, v_2, \ldots, v_t)$ is a walk, the \emph{length} $\ell(w) = t$, of the walk is one less than the number of elements in the sequence. The \emph{distance} between $u,v$ is $D(u,v) = \min\{\ell(\rho) : \rho \text{ is a path from $u$ to $v$}\}$. The \emph{diameter} of $G$ is $D(G) = \max \{D(u,v) : u,v \in V\}$. We may write $D_G$ for $D(G)$.

A \emph{subgraph} $G'$ of a graph $G$, is a graph such that $V(G')\subseteq V(G)$ and $E(G')\subseteq E(G)$. We write $G'\subseteq G$.

The following simple lemma is very useful and common in the study of graphs. See, e.g., \cite{deistel}.
\\
\begin{lemma}[\textbf{Handshaking Lemma}]\label{Handshaking Lemma}
For an undirected graph $G=(V,E)$, $\sum_{u\in V}d(u) = 2|E|$.
\end{lemma}
Note, there is no requirement that $G$ be connected or simple. 

\proofstart
Using equation (\ref{nonSimpGraphDegree}),
\begin{equation}
\sum_{u\in V}d(u)= \sum_{u\in V}|\{e \in E : e=(u,x), x \neq u\}| + 2 \sum_{u\in V}|\{e \in E : e = (u,u)\}|\label{eq:kjsjhbs}
\end{equation}
If $u,x\in V$ and $u \neq x$ then $e=(u,x)\in E$, if and only if $e=(x,u)\in E$, because the graph is undirected (though a particular edge is only included once in $E$). Hence, in the sum $\sum_{u\in V}|\{e \in E : e=(u,x), x \neq u\}|$
each edge $(u,x)$ is included twice. Thus, (\ref{eq:kjsjhbs}) becomes $\sum_{\substack{e=(u,v)\in E\\u\neq v}}2 + 2\sum_{e=(u,u)\in E}1 = 2|E|$.
\proofend

\subsection{Weighted Graphs}\label{weighted graphs}
Graphs may be \emph{weighted}, where in the context of this thesis, weights are non-negative real numbers assigned to edges in the graph. They can be represented as $G=(V,E,c)$ where $c : E \rightarrow \mathbb{R^+}$ is the weight function. We further define 
\begin{equation}
c(u) = \sum_{\substack {e = (u,x)\\ x\neq u,\ e \in E}}c(e) +\sum_{\substack{e=(u,u)\\ e \in E }}2c(e) \label{c(u)}
\end{equation}
and
\begin{equation}
c(G) = \sum_{u\in V}c(u). \label{c(G)}
\end{equation}
The \emph{weight} of the graph, $w(G)$, is
\begin{equation}
w(G) = \sum_{e\in E}c(e). \label{w(G)}
\end{equation}    
By the same arguments as the Handshaking Lemma \ref{Handshaking Lemma}, $w(G) = c(G)/2$.

For convenience, we may also define 
\begin{equation}
c(u,v) = \sum_{\substack{e = (u,v)\\ e \in E, }}c(e)\label{condSum}
\end{equation}
that is, the total weight of edges between $u,v$.

Analogous definitions can be given for directed graphs.
   
\subsection{Examples}\label{graphdefnsExamples}
We define some specific classes of graphs which will feature in subsequent chapters. All are simple, connected and undirected. Without loss of generality, we may assume that vertices are labelled $[0,n-1]$, where $n=|V|$.

\begin{description}
  \item[\emph{Complete graph}]\hfill \\
  										 The complete graph on $n$ vertices, denoted by $K_n$ is the graph such that $E=\{(u,v) : u, v\in V, u\neq v\}$, and so $|E| = \binom{n}{2} = \frac{n(n-1)}{2}$.
  \item[\emph{Path graph}]\hfill \\ 
  													The path graph on $n$ vertices, or \emph{$n$-path}, denoted by $P_n=(0,1,2,\ldots,n-1)$. $E=\{(0,1), (1,2), \ldots, (n-2, n-1)\}$. It has $|E| = |V|-1$.
  \item[\emph{Cycle graph}]\hfill\\
  												 The cycle graph (or simply, \emph{cycle}) on $n$ vertices, denoted by $\mathbb{Z}_n=(0,1,2,\ldots,n-1, 0)$. $E=\{(0,1), (1,2), \ldots, (n-2, n-1), (n-1,0)\}$. It is same as $P_n$, with an additional edge $(n-1,0)$ connecting the two ends. It has  $|E| = |V|$. 
\end{description}

We say a graph $G$ \emph{has a cycle} if $\mathbb{Z}_r \subseteq G$ for some $r$.

\textbf{Trees}\\
A tree $T$ is a graph satisfying any one of the following equivalent set of conditions (see, e.g., \cite{deistel}): $T$ is connected and has no cycles; $T$ has no cycles, and a cycle is formed if any edge is added to $T$; $T$ is connected, and it is not connected anymore if any edge is removed from $T$; Any two vertices in $T$ can be connected by a unique path; $T$ is connected and has $n - 1$ edges.

If, for a graph $G$, there is some tree $T \subseteq G$ such that $V(T)=V(G)$, then $T$ is called a \emph{spanning tree} of $G$.

\section{Markov Chains}\label{MarkovChainFirstDef}
Let $\Omega$ be some finite set. A \emph{Markov chain} is a sequence $\mathbf{X}=(X_0, X_1, \ldots)$ of random variables with $X_i \in \Omega$ having the Markov property, that is, for all $t \geq 0$
\begin{equation}
\Pr(X_{t+1}=x  \mid X_1=x_1, \ldots, X_t=x_t) = \Pr(X_{t+1}=x \mid X_t=x_t).\label{MemProp}
\end{equation}
If, in addition, we have 
\begin{equation}
\Pr(X_{t+1}=a \mid X_t=b) = \Pr(X_{t}=a \mid X_{t-1}=b)\label{TimeHom}
\end{equation}
for all $t \geq 1$, then the Markov chain is \emph{time-homogeneous}. Such Markov chains can be defined in terms of the tuple $\mathcal{M}=(\Omega, \mathbf{P}, X_0)$ where $\mathbf{P} = [P_{i,j}]$ is the $|\Omega| \times |\Omega|$ \emph{transition matrix} having entries $P_{i,j}=\mathbf{P}[i,j] = \Pr(X_{t+1}=j \mid X_t=i)$. The first element of the sequence $X_0$ is drawn from some distribution on $\Omega$, and in many applications this distribution is concentrated entirely on some known starting state.

Equations \eqref{MemProp} and \eqref{TimeHom} together express the fact that, given knowledge of $X_{t-1}$, we have a probability distribution on $X_t$, and this distribution is independent of the history of the chain before $X_{t-1}$. That is, if $X_{t-1}$ is known, any knowledge of $X_s$ for $s < t-1$ (should it exist) does not change the distribution on $X_t$. This is called the \emph{Markov property} or \emph{memoryless property}.

Without loss of generality, label the $n$ states of the Markov chain $[1,n]$. Let $\mathbf{p}(t) = [p_1(t),p_2(t),\ldots,p_n(t)]$ be the vector representing the distribution on states at time $t$. The first state $X_0$ will be drawn from some distribution $\mathbf{p}(0)$, possibly concentrated entirely on one state.  It is immediate that we have the relation 
\begin{equation*}
p_i(t) = \sum_{j = 1}^{n}p_j(t-1)P_{j,i}
\end{equation*}
for any state of the Markov chain $i$. Alternatively,
\begin{equation*}
\textbf{p}(t) = \mathbf{p}(t-1)\mathbf{P}.
\end{equation*}
For any $s \geq 1$, define the $s$-step transition probability 
\[
P_{i,j}^{(s)} = \Pr(X_{t+s} =j \mid X_t=i)
\]
and let $\mathbf{P^{(s)}} = [P_{i,j}^{(s)}]$ be the corresponding transition matrix. Observe that $\mathbf{P^{(1)}}=\mathbf{P}$ and that
\[
P_{i,j} ^{(s)}= \sum_{k = 1}^{n}P_{i,k}P_{k,j}^{(s-1)}.
\]
Thus, 
\[
\mathbf{P^{(s)}} = \mathbf{P}\mathbf{P^{(s-1)}}
\]
and by induction on $s$,
\[
\mathbf{P^{(s)}} = \mathbf{P}^s.
\]
The above is consistent with the idea that $\mathbf{P^{(0)}}= \mathbf{I}$, since this merely says that $P(X_0=i \mid X_0=i)  =1$.

For a Markov chain $\mathcal{M}$ with $X_0=i$, define $h_i(j) = \min\{t\geq 0: X_t=j\}$ and $h^+_i(j) = \min\{t\geq 1: X_t=j\}$. We define the following 
\begin{description}
  \item[\emph{Hitting time from $i$ to $j$}]\hfill \\
  										 $\hit[i,j] =\E[h_i(j)]$.
  \item[\emph{First return time to $i$}]\hfill \\
  										 $\fretT[i]=\E[h^+_i(i)]$.										 
  \item[\emph{Commute time  between $u$ and $v$}] \hfill \\
  										$\com[i,j]=\E[h_i(j) + h_j(i)] = \hit[i,j] + \hit[j,i]$ by linearity of expectation.
\end{description}

Observe that $\hit[i,i]=0$ and $\hit[i,j]\geq 1$ for $i\neq j$. Note, furthermore, that it is not generally the case that $\hit[i,j]=\hit[j,i]$, although in some classes of Markov chains it is (examples would be random walks on the complete graph or the cycle, as we shall see in chapter \ref{Techniques and Results for Hitting and Cover Times}).

The definition of hitting time can be generalised to walks starting according to some distribution $\mathbf{p}$ over the states: 
\[
\hit[\mathbf{p}, a] = \sum_{i}\mathbf{p}_i\hit[i,a].
\]

\section{Random Walks on Graphs}\label{Random Walks on Graphs}
A walk on a graph, as defined in section \ref{graphdefns}, is a sequence of vertices connected by edges $(v_0, v_1, v_2, \ldots)$. A \emph{random walk} is a walk which is the outcome of some random process, and a \emph{simple random walk} is a random walk in which the next edge transitioned is chosen uniformly at random from the edges incident on a vertex.

Random walks on graphs are a specialisation of Markov chains. For a graph $G=(V,E)$, the state space of the Markov chain is the set of vertices of the graph $V$, and a transition from a vertex $u$ is made by choosing uniformly at random (\uar) from the set of all incident edges and transitioning that edge. For undirected graphs, an edge can be traversed in either direction, and a loop counts twice. In a directed graph, the convention is that an edge is traversed in the direction of the arc, and so, where there is a directed loop, only one end can be transition.

We give formal definitions. Let $G=(V,E)$ be an unweighted, undirected simple graph. The Markov chain $\mathcal{M}_G=(V, \textbf{P}, .)$ has transition matrix $P_{u,v} = 1/d(u)$ if $(u,v) \in E$ otherwise $P_{u,v} = 0$.

More generally, when $G=(V,E,c)$ is weighted, undirected (and not necessarily simple)
\begin{equation}
P_{u,u} = \sum_{\substack{e=(u,u) \\ e \in E}}\frac{2c(e)}{c(u)} \label{P_uu}
\end{equation}  
and 

\begin{equation}
P_{u,v} = \sum_{\substack{e=(u,v)\\  e \in E}}\frac{c(e)}{c(u)} \label{P_uv}
\end{equation} 
if $v \neq u$. 

Observe, that by the above definitions, and the definition of $c(u)$, a walk on an unweighted graph is the same (meaning, has the same distribution) as a walk on a uniformly weighted graph (that is, all edges have the same weight). Conventionally, when an unweighted graph is treated as a weighted graph, edges are given unit weight.

Analogous definitions can be given for directed graphs.

Some more notation concerning random walks: Let $\mathcal{W}_{u}$ denote a random walk started from a vertex $u$ on a graph $G=(V,E)$.
Let $P_{u}^{(t)}(v)=\Pr(\mathcal{W}_{u}(t)=v)$. 

For a random walk $\mathcal{W}_{u}$, let $c_u = \max_{v \in V}h_u(v)$ where $h_u(v)$ was defined in section \ref{MarkovChainFirstDef}. In addition to the quantities defined in that section (which are defined also for random walks on graphs, since they are a type of Markov chain), we define the following

\begin{description}
  \item[\emph{Cover time of $G$ from $u$}]\hfill \\ 
  													$\cov_u[G]=\E[c_u]$.
  \item[\emph{Cover time of $G$}]\hfill \\
  												 $\cov[G]=\max_{u \in V}\cov_u[G]$.
\end{description}

\chapter{Theory of Markov Chains and Random Walks}\label{Theory of Markov Chains and Random Walks}

\section{Classification of States}\label{MCStateClassification}
The states of a Markov chain exhibit different behaviours in general. It is often the case that some of these properties can be ascertained by visual inspection of the graph of the chain, particularly when the graph is small. 

Much of what follows is standard material in an introduction to the topic. Aside from minor modifications, we quote heavily from \cite{mitz} for many of the following definitions and lemmas.

As before, we shall assume without loss of generality that the states of a chain with $n$ states are labelled $[1,n]$. 
\\

\begin{definition}[\cite{mitz}]
A state $j$ is \emph{accessible} from a state $i$ if $P_{i,j}^{(t)}>0$ for some integer $t\geq 0$. If two states $i$ and $j$ are accessible from each other we say they \emph{communicate} and we write $i \leftrightarrow j$.
\end{definition}
\begin{quote}
In the graphical representation of a chain $i \leftrightarrow j$ if and only if there is a directed path from $i$ to $j$ and there is a directed path from $j$ to $i$. 
\end{quote}
\hfill-Extracted from \cite{mitz} p.164, with minor modifications.

For random walks on undirected graphs, this is equivalent to a path existing between $i$ and $j$.
The following lemma is easy to confirm, and we omit the proof. 
\\

\begin{proposition}[\cite{mitz}]
The communicating relation defines an equivalence relation, that is, it is
\begin{enumerate}
	\item \emph{reflexive} - for any state $i$, $i \leftrightarrow i$;
	\item \emph{symmetric} - $i \leftrightarrow j \Rightarrow j \leftrightarrow i$;
	\item \emph{transitive} - $i \leftrightarrow j$ and $j \leftrightarrow k$ $\Rightarrow i \leftrightarrow k$. 
\end{enumerate}
\end{proposition}
Note that a self-loop is not required for a state to be reflexive, since $P_{i,i}^{(0)} = 1$ by definition. Thus, the communication relation partitions the states into disjoint equivalence classes called \emph{communicating classes}. The following corollary is a simple consequence
\\

\begin{corollary}
A chain cannot return to any communicating class it leaves. 
\end{corollary}
For random walks on undirected graphs, the communicating classes are the connected components of the graph.
\\

\begin{definition}[\cite{mitz}]
A Markov chain is \emph{irreducible} if all states belong to one communicating class.
\end{definition}
Random walks on undirected graphs are therefore irreducible if and only if the graph is connected (i.e., a single component). More generally, a Markov chain is irreducible if and only if the graphical representation is strongly connected (\cite{mitz}).

Denote by $f_{i,j}^{(t)}$ the probability that, starting at state $i$, the first time the chain visits state $j$ is $t$; that is
\[
f_{i,j}^{(t)} = \Pr(X_t = j \text{ and, for }1 \leq s \leq t-1, X_s \neq j \mid X_0 = i).
\]

\begin{definition}[\cite{mitz}]
A state $i$ is \emph{recurrent} if $\sum_{t \geq 1} f_{i,i}^{(t)} = 1$ and it is \emph{transient} if $\sum_{t \geq 1} f_{i,i}^{(t)} < 1$. A Markov chain is recurrent if every state in the chain is recurrent.   
\end{definition}

\begin{quote} 
A recurrent state is one which, if visited by the chain, will, with probability $1$, be visited again. Thus, if a recurrent state is ever visited, it is visited an infinite number of times. If a state is transient, there is some probability that the chain will never return to it, having visited it. For a chain at a transient state $i$, the number of future visits is given by a geometrically distributed random variable with parameter $p=\sum_{t \geq 1} f_{i,i}^{(t)}$. If one state in a communicating class is transient (respectively, recurrent) then all states in that class are transient (respectively, recurrent).
\end{quote}
\hfill-Extracted from \cite{mitz}, p.164, with minor modifications.

Recalling the definition of $\hit[i,j]$ from section \ref{MarkovChainFirstDef}, we have $\sum_{t \geq 1} t\cdot f_{i,j}^{(t)}=\hit[i,j]$ for $i\neq j$ and $\sum_{t \geq 1} t\cdot f_{i,i}^{(t)}=\fretT[i]$. It is not necessarily the case that a recurrent state has $\fretT[i] < \infty$;
\\

\begin{definition}[\cite{mitz}]
A recurrent state $i$ is \emph{positive recurrent} if $\fretT[i] < \infty$, otherwise it is \emph{null recurrent}. 
\end{definition}
An example of a Markov chain with a null recurrent state is given in \cite{mitz} chapter 7; it has an infinite number of states, in fact,
\\

\begin{lemma}[\cite{mitz}]\label{recurrentLemma}
In a finite Markov chain:
\begin{enumerate}
	\item At least one state is recurrent.
	\item All recurrent states are positive recurrent.
\end{enumerate}
\end{lemma}
The proof is left as an exercise, thus we include our own.

\proofstart
1. Since there are a finite number of communicating classes, and since once the chain leaves a communicating class it cannot return, it must eventually settle into one communicating class. Thus, at least one state in this class is visited an unbounded number of times after the chain enters it.

2. The communicating class of a recurrent state has no transition outside that class, since otherwise there would be a positive probability of no return. Consider a recurrent state $i$ of an $n$-state Markov chain. Let $\mathcal{C}(i)$ denote the communicating class of $i$. Let $p$ be the largest transition probability less than $1$ from any of the states in $\mathcal{C}(i)$. Any walk of the chain of length $\ell \geq n$ in $\mathcal{C}(i)$ that avoids $i$ must include at least one transition with probability at most $p$. Then for any $j\in \mathcal{C}(i), j\neq i$, 
\ignore{\begin{eqnarray*}
\hit[j,i]&=& \sum_{t \geq 1}\Pr(h_{j,i} \geq t)\\
 &\leq& n\sum_{m \geq 1}{p}^m\\
 &=& \frac{np}{1-p}\\
  &<& \infty.
\end{eqnarray*}
}
\begin{equation*}
\hit[j,i]= \sum_{t \geq 1}\Pr(h_{j,i} \geq t) \leq n\sum_{m \geq 0}{p}^m = \frac{n}{1-p}
  < \infty.
\end{equation*}
\proofend

The above gives the following
\\

\begin{corollary}
For a finite Markov chain, if any state $i$ is recurrent, then all of states of the communicating class of $i$ are positive recurrent.
\end{corollary}

We next discuss periodicity of Markov chains. As suggested by the name, periodicity is a notion of regular behaviour of Markov chains. As a simple example of periodic behaviour, consider a $2$-state Markov chain with states $\{1,2\}$, with each state having a transition to the other with probability $1$. If the chain starts at state $i\in\{1,2\}$, then it will be at state $i$ for all even time steps (including time $0$), and it will be at the other state at all odd times. This oscillatory behaviour means that the distribution of the chain on states can never converge, and this hints at the importance of periodicity - or lack of it.
\\

\begin{definition}[\cite{mitz}]
A state $j$ in a Markov chain is \emph{periodic} if there exists an integer $\Delta >1$ such that $\Pr(X_{t+s}=j\mid X_t=j)=0$ unless $s$ is divisible by $\Delta$. A discrete time Markov chain is periodic if any state in the chain is periodic. A state or chain that is not periodic is \emph{aperiodic}.  
\end{definition} 
There is an equivalent definition of an aperiodic state; \cite{Norris}  p.40, gives the following definition.
\\
\begin{definition}[\cite{Norris}]
A state $i$ is \emph{aperiodic} if $P_{i,i}^{(t)}>0$ for all sufficiently large $t$. 
\end{definition}
This is followed by the following theorem, which establishes an equivalence between the two definitions.
\\

\begin{theorem}[\cite{Norris}]
A state $i$ is aperiodic if and only if the set $S=\{t : P_{i,i}^{(t)}>0\}$ has no common divisor other than $1$. 
\end{theorem} 
In \cite{Norris}, the proof is left as an exercise to the reader, and so we present a proof below.

\proofstart\\
$\Rightarrow$ If $S$ has a common divisor $\Delta>1$, then any $t\geq 1$ that is not a multiple of $\Delta$ is not in $S$, and so  $P_{i,i}^{(t)}=0$ in this case.

$\Leftarrow$ Let $t_1=\min S$. Since $t_1$ has a finite number of factors, and since for any $S',S'' \subseteq S$ we have $S'\subseteq S''\Rightarrow \text{gcd}(S')\geq \text{gcd}(S'')$, we deduce that there must be some finite $S'\subseteq S$ with $t_1 \in S'$ such that $\text{gcd}(S') = 1$. By the extended Euclidean algorithm (see, e.g.\ \cite{numbertheorybook}) -or, in fact, B\'{e}zout's lemma - there must be some $a_j \in \mathbb{Z}$ such that 
\begin{equation}
a_1t_1 + a_2t_2 + \ldots a_rt_r = 1, \label{Euclid1}
\end{equation}
where the $t_j$ are members of $S'$.

Now $a_j = a'_j \pmod{t_1}$ for some $0\leq a'_j< t_1$, so substituting $a'_j$ for $a_j$ in equation (\ref{Euclid1}) and taking it modulo $t_1$ we have
\begin{equation*}
\ell = a'_2t_2 + \ldots a'_rt_r \equiv 1 \pmod{t_1} \label{Euclid2}.
\end{equation*}
If $t \equiv s\pmod{t_1}$ then $t-s\ell \equiv 0\pmod{t_1}$. Furthermore, if $t \geq (t_1-1)\ell$, then $t=mt_1+s\ell$ for some non-negative integer $m$. Since any positive sum of elements in $S$ is an element in $S$,  it follows that $P_{i,i}^{(t)}>0$ .
\proofend

Finally for this section, we define and discuss \emph{ergodicity}:
\\

\begin{definition}[\cite{mitz}]
An aperiodic, positive recurrent state is an \emph{ergodic} state. A Markov chain is ergodic if all of it's states are ergodic.
\end{definition}

As a corollary to the above theorems, we have
\\

\begin{corollary}[\cite{mitz}]\label{ergoCor}
Any finite, irreducible, and aperiodic Markov chain is an ergodic chain.
\end{corollary}
\proofstart
A finite chain has at least one recurrent state by Lemma \ref{recurrentLemma} and if the chain is irreducible, then all of its states are recurrent. In a finite chain, all recurrent states are positive recurrent by Lemma \ref{recurrentLemma} and thus all states of the chain are positive recurrent and aperiodic. The chain is therefore ergodic.
\proofend

The significance of ergodicity is made clear in the following section.

\section{The Stationary Distribution}\label{MCStationaryDist}
Recall that the distribution on the states evolves with this relation
\begin{equation*}
\textbf{p}(t+s) = \mathbf{p}(t)\mathbf{P}^s.
\end{equation*}
A fundamental question is when, if ever, there exists a distribution that remains invariant under the operation of post-multiplication by the transition matrix. 
\\

\begin{definition}[\cite{mitz}]
A \emph{stationary distribution} (also called an \emph{equilibrium distribution}) of a Markov chain is a probability distribution $\mathbf{\pi}$ such that
\[
\mathbf{\pi} = \mathbf{\pi P}.
\]
\end{definition}

A chain in the stationary distribution will continue to be so after subsequent transitions.

We now quote an important theorem from \cite{mitz}, but omit the proof, which although not difficult, is fairly lengthy.
\\

\begin{theorem}[\cite{mitz}]\label{ergodicTheorem}
Any finite, irreducible, and ergodic Markov chain has the following properties:
\begin{enumerate}
	\item the chain has a unique stationary distribution $\mathbf{\pi}=(\pi_1, \pi_2, \ldots \pi_n)$;
	\item for all $i$ and $j$, the limit $\lim_{t \rightarrow \infty}P_{i,j}^{(t)}$ exists and is independent of $i$;
	\item $\pi_j = \lim_{t \rightarrow \infty}P_{i,j}^{(t)} = 1/\fretT[j]$.
\end{enumerate}
\end{theorem}
There are a number of other proofs of Theorem \ref{ergodicTheorem}; \cite{Norris} gives a proof based on coupling; \cite{LPW} gives a different different proof also based on coupling; \cite{LovaszSurvey} gives a treatment specialised for connected undirected graphs using the framework of linear algebra, in particular, the eigenvalues of the transition matrix and the Perron-Frobenius theorem to show the existence and convergence to the stationary distribution.

For a random walk on an undirected graph that is finite, connected and not bipartite, the conditions of ergodicity and thus the conditions for Theorem \ref{ergodicTheorem} are satisfied.

\section{Random Walks on Undirected Graphs}\label{randwalksundirectedgraphs}
For a random walk on an undirected graph a pair of vertices $u,v$ communicate if there is a path between them. Furthermore, as stated above, the communicating classes are the connected components of the graph and so the random walk is irreducible if and only if the graph is a single connected component. For aperiodicity, we quote the following lemma from \cite{mitz}, along with the accompanying proof.
\\
\begin{lemma}[\cite{mitz}]\label{nonBipartiteErgo}
A random walk on an indirected graph $G$ is aperiodic if and only if $G$ is not bipartite.
\end{lemma} 
\proofstart
A graph is bipartite if and only if it does not have cycles with an odd number of edges. In an undirected graph, there is always a path of length $2$ from a vertex to itself. If the graph is bipartite, then the random walk is periodic with period $2$. If the graph is not bipartite, then it has an odd cycle and by traversing the cycle, an odd-length path can be constructed from any vertex to itself. It follows that the Markov chain is aperiodic.
\proofend

Thus, by Corollary \ref{ergoCor} and Lemma \ref{nonBipartiteErgo} it is seen that
\\

\begin{corollary}
A random walk on an undirected graph that is finite, connected and non-bipartite is ergodic.
\end{corollary}

The existence of the stationary distribution and the convergence of the walk to it is thus established. 
\\

\begin{theorem}[\cite{mitz}]\label{Th:RandomWalkStationDist}
A random walk on an undirected graph $G=(V,E)$ that is finite, connected, and not bipartite converges to the stationary distribution $\mathbf{\pi}$ where, for any vertex $v \in V$
\[
\pi_v = \frac{d(v)}{2|E|}
\] 
\end{theorem}

\proofstart
By the handshaking lemma, $\sum_{v\in V}d(v) = 2|E|$. Thus, it follows that 
\[
\sum_{v\in V}\pi_v = \sum_{v\in V}\frac{d(v)}{2|E|} = 1
\]
and so $\mathbf{\pi}$ is a proper probability distribution over $V$.
Let $\mathbf{P}$ be the transition matrix of the walk on $G$ and let $N(v)$ represent the neighbour set of $v$. The relation $\mathbf{\pi} = \mathbf{\pi P}$ is equivalent to 
\[
\pi_v = \sum_{u \in N(v)}\frac{d(u)}{2|E|}\frac{1}{d(u)} = \frac{d(v)}{2|E|} 
\]
and the theorem follows.
\proofend

The above proof is for simple graphs. It can be generalised for graphs with loops and/or parallel edges by with stationary probability $\frac{d(v)}{2|E|}$ where now $d(v)$ is given by (\ref{nonSimpGraphDegree}). The theorem can be further generalised to (not necessarily simple) weighted graphs with stationary probability $\pi_v = \frac{c(v)}{c(G)}$. See section \ref{Random Walks on Graphs} for relevant definitions.

\section{Time Reversal and a Characterisation of Random Walks on Undirected Graphs}
This section discusses a characterisation of Markov chains that precisely captures random walks on undirected graphs, including non-simple graphs (those with loops and/or parallel edges) as well as weighted undirected graphs. To do so we introduce the following definition, taken from \cite{Norris}
\\

\begin{definition}
Let $(X_t)_{0\leq t \leq T}$ be a (sub)sequence of states of a Markov chain $\mathcal{M}=(\Omega, \mathbf{P}, \mathbf{\pi})$ where $\mathbf{\pi}$ is the stationary distribution of $\mathcal{M}$. The \emph{time reversal} of $(X_t)_{0\leq t \leq T}$ is the sequence $(Y_t)_{0\leq t \leq T}$ where $Y_t = X_{T-t}$.
\end{definition}

The following theorem is given in \cite{Norris}, and we omit the proof.
\\

\begin{theorem}\label{revTheorem}
Let $(X_t)_{0\leq t \leq T}$ be a (sub)sequence of states of a Markov chain $\mathcal{M}=(\Omega, \mathbf{P}, \mathbf{\pi})$ where $\mathbf{\pi}$ is the stationary distribution of $\mathcal{M}$. Then $(Y_t)_{0\leq t \leq T}$ is a (sub)sequence of states of a Markov chain $\mathcal{\widehat{M}}=(\Omega, \mathbf{\widehat{P}}, \mathbf{\pi})$ where $\mathbf{\widehat{P}} = [\widehat{P}_{i,j}]$ is given by
\begin{equation}
\pi_j\widehat{P}_{j,i} = \pi_iP_{i,j} \text{ \, \, \, for all } i,j \label{revEq}
\end{equation}
and $\mathbf{\widehat{P}}$ is also irreducible with stationary distribution $\mathbf{\pi}$.
\end{theorem}
This leads us to the following definition
\\

\begin{definition}
If $\mathbf{\widehat{P}} = \mathbf{P}$ then the Markov chain $\mathcal{M}$ is said to be \emph{time reversible}. 
\end{definition}

If a Markov chain is reversible, then (\ref{revEq}) becomes 
\begin{equation}
\pi_jP_{j,i} = \pi_iP_{i,j} \text{ \, \, \, for all } i,j \label{detailedBalance}
\end{equation}
known as the \emph{detailed balance} condition.

Conversely, if the detailed balance condition is satisfied for some distribution $\mathbf{p}$, that is
\[
p_jP_{j,i} = p_iP_{i,j} \text{ \, \, \, for all } i,j 
\]
then $\mathbf{p}$ is the stationary distribution, which, along with irreducibility and Theorem \ref{revTheorem}, implies the following
\\

\begin{corollary}\label{detailedBalanceCor}
An irreducible Markov chain with a stationary distribution $\mathbf{\pi}$ is reversible if and only if it satisfies the detailed balance condition (\ref{detailedBalance}).
\end{corollary}

The next theorem characterises Markov chains as random walks. 
\\

\begin{theorem}\label{TimeRevMCRandomWalk}
Random walks on undirected weighted graphs are equivalent to time reversible Markov chains. That is, every random walk on a weighted graph is a time reversible Markov chain, and every time reversible Markov chain is a random walk on some weighted graph.
\end{theorem}
\proofstart
The transition matrix for a random walk on a weighted undirected graph $G=(V,E)$ defines 
\[
P_{u,v} = \frac{c(u,v)}{c(u)}
\]
where $c(u,v)$ is defined by (\ref{condSum}) (and therefore valid for non-simple graphs).
Thus 
\[
\pi_uP_{u,v} = \frac{c(u)}{c(G)}\frac{c(u,v)}{c(u)} = \frac{c(u,v)}{c(G)} = \frac{c(v,u)}{c(G)}=\frac{c(v)}{c(G)}\frac{c(v,u)}{c(v)} = \pi_vP_{v,u}.
\]
Hence, the detailed balance condition is satisfied and by Corollary \ref{detailedBalanceCor}, the random walk is a reversible Markov chain.

Consider some reversible Markov chain $\mathbf{M} = (\Omega,\mathbf{P},.)$ with transition matrix $\mathbf{P}$ and where, as before, we assume the states are labelled $[1,n]$. We define a weighted graph $G=(V,E)$ based on $M$ as follows: Let $V = \Omega$, let $(i,j) \in E$ if and only if $P_{i,j}>0$ and weight the edge $(i,j)$ as $c(i,j) = \pi_iP_{i,j}$. By reversibility, the detailed balance equations imply $c(i,j) = c(j,i)$, hence weights are consistent and the weighted graph is proper. The random walk on the graph, by construction, has transition matrix $\mathbf{P}$.
\proofend

\chapter{The Electrical Network Metaphor}\label{networktheory}
In this chapter we give an introduction to the electrical network metaphor of random walks on graphs and present some of the concepts  and results from the literature that are used in subsequent parts of this thesis. Although a purely mathematical construction, the metaphor of electrical networks facilitates the expression of certain properties and behaviours of random walks on networks, and provides a language for which to describe these properties and behaviours. The classical treatment of the topic is \cite{DoyleSnell}. The recent book \cite{LPW} presents material within the more general context of Markov chains. Other treatments of the topic are found in \cite{AlFi} and \cite{LovaszSurvey}. 

We first present some definitions. 

\section{Electrical Networks: Definitions}\label{Electrical Networks: Definitions}
An \emph{electrical network} is a connected, undirected, finite, graph $G=(V,E)$ where each edge $e \in E$ is has a non-negative weight $c(e)$. The weight is called the \emph{conductance}. The \emph{resistance} of an edge $e$, $r(e)$ is defined as the reciprocal of the conductance, $1/c(e)$ if $c(e)$ is finite, and is defined as $\infty$ if $c(e)=0$. It is quite often the case in the literature that in the context of electrical networks, the vertices of the network are referred to as \emph{nodes}. We shall use the terms `graph' and `network', and `vertex' and `node', interchangeably in the context of electrical network metaphor.  

A random walk on an electrical network is a standard random walk on a weighted graph, as per the definition of section \ref{Random Walks on Graphs}. A random walk on an electrical network are therefore equivalent to time reversible Markov chains by Theorem \ref{TimeRevMCRandomWalk}.

Since edges are always weighted in the context of electrical networks, we shall use the notion $G=(V,E,c)$ for the network, where the third element of the tuple is the weighting function on the edges.

\section{Harmonic Functions}
Given a network $G=(V,E,c)$, a function $f: V \rightarrow \mathbb{R}$ is \emph{harmonic} at $u \in V$ if it satisfies
\begin{equation}
f(u) = \sum_{v \in V}P_{u,v}f(v), \label{harmonicequation}
\end{equation}
where $P_{u,v}$ is defined by equations \eqref{P_uu} and \eqref{P_uv}.
\\

For some set $V_B \subset V$, called \emph{boundary nodes}, call the complement $V_I = V \setminus V_B$ \emph{internal nodes}.
\\

\begin{lemma}\label{maxminharmlemma}
For a function $f_{V_B}: V_B \rightarrow \mathbb{R}$, any extension of $f_{V_B}$ to $V$, $f: V \rightarrow \mathbb{R}$ that is harmonic on the internal nodes $V_I$ attains it's minimum and maximum values on the boundary. That is, there are some $b,B \in V_B$ such that for any $v \in V$, $f(b) \leq f(v) \leq f(B)$. 
\end{lemma}
\proofstart
We start with the upper bound. Let $M= \max_{x \in V} f(x)$ and let $V_M = \{x \in V : f(x)=M\}$.  If $V_B \cap V_M \neq \emptyset$ then we are done. If not, then $V_M \subseteq V_I$ and we choose some $x \in V_M$. Since $f$ is harmonic on $V_I$, $f(x)$ is a weighted average of its neighbours. It follows that $f(y)=M$ for each neighbour $y$ of $x$, i.e., $y \in V_M$. Iterating this repeatedly over neighbours, we see that any path in the network $x=x_0, x_1, \ldots, x_r=z$ must have the property that each $x_i \in V_M$. Since a network is connected by definition (see \ref{Electrical Networks: Definitions}), there must be a path from $x$ to some $z \in V_B$, in which case we get a contradiction, and therefore deduce that $V_B \cap V_M \neq \emptyset$. A similar argument holds for the  minimum. 
\proofend
\\

\begin{theorem}\label{uniqueharmonictheorem}
For a function $f_{V_B}: V_B \rightarrow \mathbb{R}$ the extension of $f_{V_B}$ to $V$, $f: V \rightarrow \mathbb{R}$ is unique if $f$ is harmonic on all the internal nodes $V_I$.
\end{theorem}
\proofstart
Suppose there are functions $f,g$ which extend $f_{V_B}$ and are harmonic on each node in $V_I$. Consider $h=f-g$. This function has $h(v) = 0$ for any $v \in V_B$, and is harmonic on $V_I$. It follows by Lemma \ref{maxminharmlemma} that $h(v)=0$ for any $v \in V_I$ as well. Therefore $g=h$. 
\proofend

The problem of extending a function $f_{V_B}$ to a function $f$ harmonic on $V_I$ is known as the \emph{Dirichlet problem}, in particular, the \emph{discrete} Dirichlet problem, (in contrast to the continuous analogue).  For electrical networks (and in fact more generally, for irreducible Markov chains), a solution to the Dirichlet problem always exists, as provided by the following function.
\\

\begin{theorem}
Let $G=(V,E,c)$ be a network, $V_B \subset V$, be a set of boundary nodes and $V_I = V \setminus V_B$ be the internal nodes. Let $f_{V_B}: V_B \rightarrow \mathbb{R}$ be a function on the boundary nodes. Let $\mathbf{X}(\mathcal{W})$ be a random variable that represents the first $v \in V_B$ visited by a weighted random walk $\mathcal{W}$ on $G=(V,E,c)$ started at some time. The function $f(v)=\E[f_{V_B}(X(\mathcal{W}_v))]$, where $v \in V$, extends $f_{V_B}$ and is harmonic on $V_I$.  
\end{theorem} 
\proofstart
For a node $v \in V_B $, 
\[
f(v) = \E[f_{V_B}(X(\mathcal{W}_v))] = \E[f_{V_B}(v)] = f_{V_B}(v).
\] 
thus, $f$ is consistent with $f_{V_B}$ on the boundary nodes.

For $u \in V_I$,
\begin{eqnarray*}
f(u) &=& \E[f_{V_B}(X(\mathcal{W}_u))]\\
 &=& \sum_{v \in V}\E[f_{V_B}(X(\mathcal{W}_u)) \mid \mathcal{W}_u(1)=v]P_{u,v}\\
 &=& \sum_{v \in V}\E[f_{V_B}(X(\mathcal{W}_v))]P_{u,v}\\
 &=& \sum_{v \in V}P_{u,v}f(v).
\end{eqnarray*} 
This proves harmonicity on $V_I$.
\proofend

\section{Voltages and Current Flows}\label{voltageCurrent}
Consider a network $G=(V,E,c)$ and let a pair of nodes $a$ and $z$ be known as the \emph{source} and \emph{sink} respectively. Treating $a,z$ as the only elements of a boundary set on the network, a function $W$ harmonic on $V \setminus \{a,z\}$ is known as a \emph{voltage}. 

For an edge $e=(u,v)$, denote by $\overrightarrow{e}=(\overrightarrow{u,v})= (\overleftarrow{v,u})$ an orientation of the edge from $u$ to $v$. Furthermore, if $\overrightarrow{e}=(\overrightarrow{u,v})$ then let $\overleftarrow{e}=(\overleftarrow{u,v}) = (\overrightarrow{v,u})$. A \emph{flow}  $\varphi: \overrightarrow{E}\rightarrow \mathbb{R}$ where $\overrightarrow{E} = \{\overrightarrow{e} : e \in E\}\cup  \{\overleftarrow{e}: e \in E\}$ is a function on oriented edges which is antisymmetric, meaning that $\varphi(\overrightarrow{e}) = -\varphi(\overleftarrow{e})$. For a flow $\varphi$, define the \emph{divergence} of $\varphi$ at a node $u$ to be 
\[
\operatorname{div}\varphi(u) = \sum_{\overrightarrow{e}=(\overrightarrow{u,v})\in \overrightarrow{E}}\varphi(\overrightarrow{e}).
\] 
Observe, for a flow $\varphi$, 
\begin{equation}
\sum_{u \in V}\operatorname{div}\varphi(u) = \sum_{u \in V}\sum_{\overrightarrow{e}=(\overrightarrow{u,v})\in \overrightarrow{E}}\varphi(\overrightarrow{e}) = \sum_{e \in E}\varphi(\overrightarrow{e})+\varphi(\overleftarrow{e}) = 0. \label{divIsZero}
\end{equation}

We term as \emph{a flow from $a$ to $z$} a flow $\varphi$ satisfying 
{\begin{enumerate}
	\item \emph{Kirchhoff's node law}
	\[
		\operatorname{div}\varphi(v) = 0 \, \, \, \text{ for all } v \notin \{a,z\},
	\]
	\item $\operatorname{div}\varphi(a) \geq 0$.
\end{enumerate}}
The \emph{strength} of a flow $\varphi$ from $a$ to $z$ is defined to be $\|\varphi\| = \operatorname{div}\varphi(a)$ and a \emph{unit flow} from $a$ to $z$ is a flow from $a$ to $z$ with strength $1$. Observe also that by (\ref{divIsZero}) $\operatorname{div}\varphi(a) = -\operatorname{div}\varphi(z)$.

Given a voltage $W$ on the network, the \emph{current flow} $I$ associated with $W$ is defined on oriented edges $\overrightarrow{e} = (\overrightarrow{u,v})$ by the following relation, known as \emph{Ohm's Law}:
\begin{equation}
I(\overrightarrow{e})=\frac{W(u)-W(v)}{r(e)} = c(e)(W(u)-W(v))\label{currentFlowEq}
\end{equation}
Conductances (resistances) are defined for an edge with no regard to orientation, so in (\ref{currentFlowEq}) we have used for notational convenience that $c(e) = c(\overrightarrow{e}) = c(\overleftarrow{e})$. We will continue to use this.

Let $G$ be a network and for some chosen boundary points $V_B \subset V$ let $f: V \rightarrow \mathbb{R}$ harmonic on the internal nodes $V_I = V \setminus V_B$. For a transformation of the form $\mathcal{T}(x) = \alpha x + \beta$, applying $\mathcal{T}$ to $f$ on the boundary points we get a new set of boundary node values given by $\mathcal{T}(f(v))$ for any $v \in V_B$. At the same time, $\mathcal{T}(f)$ is a solution to the Dirichlet problem for the new boundary node values, and so by Theorem \ref{uniqueharmonictheorem}, it is the only solution. 
\\

Now let $I_{W}$ be the current flow by a voltage $W$ on $G$ (with some chosen source and sink) and $I_{\mathcal{T}(W)}$ the current flow from the transformation $\mathcal{T}$ applied to $W$. It can be seen from the definition of current flow that $I_{\mathcal{T}(W)}= \alpha \cdot I_{W}$. In particular, this means that current flow is invariant with respect to $\beta$. Thus, assuming constant edge conductances, current flow is determined entirely by $\Delta_W=W(a)-W(z)$. We may therefore denote the current flow determined by $\Delta_W$ by $I_{\Delta_W}$. Observe $I_{0} = 0$ and $I_{\alpha \cdot \Delta_W} = \alpha \cdot I_{\Delta_W}$. Thus if, for a given $G$, any finite, non-zero current flow exists, then $\|I_{\Delta_W}\|$ as a function of $\Delta_W$ is a bijective mapping from $\mathbb{R}$ to $\mathbb{R}$. In particular, if any finite, non-zero current flow exists, then a unit current flow exists and is unique.

We show that $I$ is a flow from $a$ to $z$ when $W(a) \geq W(z)$. Firstly, consider any node $u \notin \{a,z\}$: 
\begin{eqnarray*}
\operatorname{div}I(u) = \sum_{\substack{\overrightarrow{e}=(\overrightarrow{u,v})\\ \overrightarrow{e} \in \overrightarrow{E}}}I(\overrightarrow{e})&=& \sum_{\substack{e=(u,v)\\ e \in E}}I(e)\\ 
&=& \sum_{\substack{e=(u,v)\\ e \in E}}c(e)(W(u)-W(v))\\
&=& W(u)\sum_{\substack{e=(u,v)\\ e \in E}}c(e) - c(u)\sum_{\substack{e=(u,v)\\ e \in E}}\frac{c(e)}{c(u)}W(v)\\
&=& W(u)c(u) - c(u)W(u) = 0 
\end{eqnarray*} 
The last line follows because 
\[
 \sum_{\substack{e=(u,v)\\ e \in E}}c(e) = c(u)
\]
by definition
and 
\[
\sum_{\substack{e=(u,v)\\ e \in E}}\frac{c(e)}{c(u)}W(v) = W(u)
\]
by harmonicity.

Now if $W(a) \geq W(z)$ then by Lemma \ref{maxminharmlemma}, $W(a) \geq W(v)$ for any $v \in V$. Therefore,
\[
\operatorname{div}I(a) = \sum_{\substack{e=(a,v)\\ e \in E}}c(e)(W(a)-W(v)) \geq 0.
\]
Thus, having proved both conditions of the definition, it is shown that the current flow is a flow from $a$ to $z$. Furthermore, since setting $\Delta_W = W(a)-W(z)=1$ will give some current flow $I_{1}$, setting $W(a)-W(z)=1/\|I_{1}\|$ will give a unit current flow. 

The significance of the unit current flow will become clear in the discussion of effective resistance. 

\section{Effective Resistance}
We begin with the definition
\\

\begin{definition}
Let $G=(V,E,c)$ be a network, and let $a,z$ be a pair of nodes in the network. Let $W$ be any voltage with $a,z$ treated as source and sink respectively and with $W(a) \geq W(z)$. Using the notation of Section \ref{voltageCurrent}, the \emph{effective resistance} between $a$ and $z$, denoted by $R(a,z)$ is defined as
\[
R(a,z)=\frac{\Delta_W}{\|I_{\Delta_W}\|}
\] 
\end{definition}
Clearly, for this definition to be proper, the ratio has to be invariant with respect to voltages, and indeed it was shown in section \ref{voltageCurrent} that $I_{\alpha \cdot \Delta_W} = \alpha \cdot I_{\Delta_W}$, thus preserving the ratio.
\\

It is important to note that the \emph{resistance} $r(u,v)$ of an edge $(u,v)$ is different to the \emph{effective resistance} $R(u,v)$ between the vertices $u,v$. Resistance $r(u,v)$ is $1/c(u,v)$, the inverse of conductance, which is part of the definition of the network $G=(V,E,c)$, and is the weighting function $c$ defined on an \emph{edge}. Effective resistance, on the other hand, is a property of the network, but not explicitly given in the tuple $(V,E,c)$, and it is defined between \emph{a pair of vertices}. 
\\

\begin{theorem}[\textbf{\cite{DoyleSnell} or \cite{LPW}}]\label{ReffMetric}
Effective resistance forms a metric on the nodes of a network $G=(V,E,c)$, that is, 
\textbf{(1)} $R(v,v)=0$ for any $v \in V$
\textbf{(2)} $R(u,v)\geq 0$ for any vertices $u,v \in V$
\textbf{(3)} $R(u,v) = R(v,u)$ for any vertices $u,v \in V$ 
\textbf{(4)} $R(u,w) \leq R(u,v) + R(v,w)$ for any vertices $u,v,w \in V$ (\emph{triangle inequality}).
\end{theorem}

Define the \emph{energy} $\mathcal{E}(\varphi)$  of a flow $\varphi$ on a network $G=(V,E,c)$ as 
\begin{equation}
\mathcal{E}(\varphi) = \sum_{e \in E}[\varphi(e)]^2r(e)\label{energydef}
\end{equation}

Note, the sum in (\ref{energydef}) is over unoriented edges, so each edge $e$ is considered only once. Because flow is antisymmetric by definition, the term $[\varphi(e)]^2$ is unambiguous.

The following theorem is useful in using current flows to approximate effective resistance. We shall see such an application in section \ref{Effective Resistance Lemmas}. For a proof, see, for example \cite{LPW}.
\\

\begin{theorem}[\textbf{Thomson's Principle}]\label{Thomson's Principle}
For any network $G=(V,E,c)$ and any pair of vertices $u,v \in V$,
\begin{equation}
R(u,v) = \min \{\mathcal{E}(\varphi) : \varphi \text{ is a unit flow from } u \text{ to } v\}. \label{ThomsonEq}
\end{equation}
The unit current flow is the unique $\varphi$ that gives the minimum element of the above set.
\end{theorem}

\subsection{Rayleigh's Monotonicity Law, Cutting \& Shorting}\label{RayleighStuff}
Rayleigh's Monotonicity Law, as well as the related Cutting and Shorting Laws, are intuitive principles that play important roles in our work. They are very useful means of making statements about bounds on effective resistance in a network when the network is somehow altered. With minor alterations of notation, we quote \cite{LPW} Theorem 9.12, including proof.
\\

\begin{theorem}[\textbf{Rayleigh's Monotonicity Law}]\label{RayMonoLaw}
If $G=(V,E)$ is a network and $c,c'$ are two different weightings of the network such that $r(e)\leq r'(e)$ for all $e \in E$, (recall, $r(e) = 1/c(e)$), then for any $u,v \in V$,
\[
R(u,v) \leq R'(u,v)
\]
where $R(u,v)$ is the effective resistance between $u$ and $v$ under the weighting $c$ (or $r$), and $R'(u,v)$ under weighting $c'$ (or $r'$).
\end{theorem} 
\proofstart
Note that $\min_{\varphi}\sum_{e \in E}[\varphi(e)]^2r(e) \leq \min_{\varphi}\sum_{e \in E}[\varphi(e)]^2r'(e)$ and apply Thomson's Principle (Theorem \ref{Thomson's Principle}).
\proofend
\\

\begin{lemma}[\textbf{Cutting Law}]\label{CuttingLaw}
Removing an edge $e$ from a network cannot decrease the effective resistance between any vertices in the network. 
\end{lemma}
\proofstart
Replace $e$ with an edge of infinite resistance (zero conductance) and apply Rayleigh's Monotonicity Law.
\proofend 
\\

\begin{lemma}[\textbf{Shorting Law}] \label{ShortingLaw}
To short a pair of vertices $u,v$ in a network $G$, replace $u$ and $v$ with a single vertex $w$ and do the following with the edges: Replace each edge $(u,x)$ or $(v,x)$ where $x \notin \{u,v\}$ with an edge $(w,x)$. Replace each edge $(u,v)$ with a loop $(w,w)$. Replace each loop $(u,u)$ or $(v,v)$ with a loop $(w,w)$. A new edge has the same conductance as the edge it replaced. Let $G'$ denote the network after this operation, and let $R$ and $R'$ represent effective resistance in $G$ and $G'$ respectively. Then, for a pair of vertices $a,z, \notin \{u,v\}$, $R'(a,z) \leq R(a,z)$, $R'(a,w) \leq R(a,u)$ and $R'(a,w) \leq R(a,v)$.
\end{lemma}
\proofstart 
Consider nodes $a,z \notin \{u,v\}$ of the network $G=(V, E, c)$, and let $G'$ denote $G$ after a Shorting operation on $u,v$. For any flow $\varphi$ from $a$ to $z$ in $G$, we can define a flow $\varphi'$ in $G'$ as follows: For any $e =(x,y)$ such that $x,y \notin \{u,v\}$, $\varphi'(\overrightarrow{e}) =\varphi(\overrightarrow{e})$. For $e=(u,u)$ or $e=(v,v)$, let $e' = (w,w)$ be the edge that replaced $e$, and have $\varphi'(\overrightarrow{e'}) =\varphi(\overrightarrow{e})$ (orientating the loops arbitrarily). For $e=(u,v)$ let $e' = (w,w)$ be the edge that replaced $e$, and have $\varphi'(\overrightarrow{e'}) =\varphi(\overrightarrow{e})$ (again, orientating the loops arbitrarily). It is easily seen that $\varphi'$ is a valid unit flow from $a$ to $z$, and that the energy on each edge is the same for both $\varphi$ and $\varphi'$. It follows that $\mathcal{E}(\varphi')$ = $\mathcal{E}(\varphi)$. Hence by Thomson's Principle, if $\varphi = \varphi_{min}$, the unit current flow from $a$ to $z$ in $G$, then $R'(a,z)\leq \mathcal{E}(\varphi')$ = $\mathcal{E}(\varphi) = R(a,z)$
\ignore{
\begin{eqnarray*}
R'(a,z) &=& \min \{\mathcal{E}(\varphi') : \varphi' \text{ is a unit flow from } a \text{ to } z \text{ in }G'\}\\
 &\leq& \min \{\mathcal{E}(\varphi) : \varphi \text{ is a unit flow from } a \text{ to } z \text{ in }G\}\\ 
 &=& R(a,z).
\end{eqnarray*}}
A similar argument can be made for $R'(a,w) \leq R(a,u)$ and $R'(a,w) \leq R(a,v)$.
\proofend

Sometimes the Shorting Law is defined as putting a zero-resistance edge between $u,v$, but since zero-resistance (infinite-conductance) edges are not defined in our presentation, we refer to the act of ``putting a zero-resistance edge'' between a pair of vertices as a metaphor for shorting as defined above.

\subsection{Commute Time Identity}
The following theorem, first given in \cite{Chandra} is a fundamental tool in our analysis of random walks on graphs, in chapter \ref{ch:The Cover Time of Cartesian Product Graphs}. The proof is not difficult but we omit it because a presentation would be lengthy. We refer to \cite{Chandra} or \cite{LPW} for a proof.
\\

\begin{theorem}[\cite{Chandra}]\label{Th:commtetimetheorem}
Let $G=(V,E,c)$ be a network. Then for a pair of vertices $u,v \in V$
\begin{equation*}
\com[u,v] = c(G)R(u,v).
\end{equation*}
(The reader is reminded that $c(G)=\sum_{v \in V}c(v) = 2\sum_{e \in E}c(e)$, as defined in (\ref{c(G)})).
\end{theorem}

\section{Parallel and Series Laws}
The parallel and series laws are rules that establish equivalences between certain structures in a network. They are useful for reducing a network $G$ to a different form $G'$, where the latter may be more convenient to analyse. We quote from \cite{LPW}, with minor modifications for consistency in notation.
\\

\begin{lemma}[\textbf{Parallel Law}]
Conductances in parallel add.
\end{lemma}
Suppose edges $e_1$ and $e_2$, with conductances $c(e_1)$ and $c(e_2)$ respectively, share vertices $u$ and $v$ as endpoints. Then $e_1$ and $e_2$ can be replaced with a single edge $e$ with $c(e) = c(e_1)+c(e_2)$, without affecting the rest of the network. All voltages and currents in $G\setminus \{e_1,e_2\}$ are unchanged and the current $I(\overrightarrow{e}) = I(\overrightarrow{e_1}) + I(\overrightarrow{e_2})$. For a proof, check Ohm's and Kirchhoff's laws with $I(\overrightarrow{e}) = I(\overrightarrow{e_1}) + I(\overrightarrow{e_2})$.
\\

\begin{lemma}[\textbf{Series Law}]
Resistances in series add.
\end{lemma}
If $v \in V\setminus \{a,z\}$, where $a$ and $z$ are source and sink, is a node of degree $2$ with neighbours $v_1$ and $v_2$, the edges $(v_1, v)$ and $(v, v_2)$ can be replaced with a single edge $(v_1,v_2)$ with resistance $r(v_1, v_2) = r(v_1, v) + r(v, v_2)$. All potentials and currents in $G\setminus \{v\}$ remain the same and the current that flows from $v_1$ to $v_2$ is  $I(\overrightarrow{v_1,v_2})=I(\overrightarrow{v_1,v}) = I(\overrightarrow{v,v_2})$. For a proof, check Ohm's Law and Kirchhoff's Law with $I(\overrightarrow{v_1,v_2})=I(\overrightarrow{v_1,v}) = I(\overrightarrow{v,v_2})$.


\chapter{Techniques and Results for Hitting and Cover Times}
\label{Techniques and Results for Hitting and Cover Times}
In this chapter we present some of the techniques for bounding hitting and cover times, as well as particular results. We start with section \ref{preciseCalc} where we give calculations for hitting and cover times of some particular graph structures. The graph classes that are the subject of the next section are important for subsequent chapters, and as examples, they serve to convey some of the techniques used to precisely calculate hitting and cover time.  This allows comparisons to be drawn with more general techniques and bounds. In general, it is quite difficult to calculate precise cover times for all but a few classes of graphs; the examples given in section \ref{preciseCalc} are amongst the simplest and most common structures studied in the literature. 

We refer the reader to section \ref{graphdefnsExamples} for reminders on definitions of graph structures and section \ref{Random Walks on Graphs} for a reminder of relevant notations and definitions relating to random walks on graphs.

We only deal with connected, undirected graphs in this chapter.

The $n$'th \emph{harmonic number}, $h(n) = \sum_{i=1}^{n}\frac{1}{i}$ is a recurring quantity, and a short hand proves useful. Note $h(n) = \ln n + \gamma + O(1/n)$ where $\gamma \approx 0.577$ (see, e.g.\, \cite{ConwayGuy}). Thus even for relatively small values of $n$, $\ln n$ is a close approximation for $h_n$.

\section{Precise Calculations for Particular Structures}\label{preciseCalc}
We deal with three particular classes of graphs: complete graphs, paths and cycle. These results are given in \cite{LovaszSurvey}, amongst others.

\subsection{Complete Graph}
\begin{theorem}\label{theorem:completecovertime}
$K_n = (V,E)$ be the complete graph on $n$ vertices.
\begin{description}
	\item[(i)] $\fretT[u] = n$ for any $u \in V$.
	\item[(ii)] $\hit[u,v] = n-1$ if $u\neq v$, for any pair of vertices $u,v\in V$.
	\item[(iii)] $\cov[K_n] = (n-1)h(n-1)$.
\end{description}
\end{theorem}
\proofstart
\textbf{(i)} $|E|=\frac{n(n-1)}{2}$ and $d(u)=n-1$ for any $u\in V$. Now use Theorem \ref{Th:RandomWalkStationDist} and  Theorem \ref{ergodicTheorem}.

\textbf{(ii)} The walk has transition probability $P_{u,v} = \frac{1}{n-1}$ if $u\neq v$ and $P_{u,u} = 0$. Thus, $\hit[u,v]$ is the expectation of a geometric random variable with parameter $ \frac{1}{n-1}$, i.e.,\ $\hit[u,v] = n-1$.

\textbf{(iii)} Let $T(r)$ be the expected number of steps until $r$ distinct vertices have been visited by the walk for the first time. By the symmetry of the graph, $T(r)$ will be invariant with respect to starting vertex. Suppose the walk starts at some vertex $u\in V$. Then $T(1)=0$. Since it will move to a new vertex in the next step, $T(2)-T(1)=1$. Suppose it visits the $r$'th distinct vertex at some time $t$, where $r<n$. Then there are $n-r$ unvisited vertices, and the time to visit any vertex from this set is a geometric random variable with probability of success $\frac{n-r}{n-1}$, the expectation of which is $\frac{n-1}{n-r}=T(r+1)-T(r)$. Thus by linearity of expectation,
\begin{eqnarray*}
\cov[K_n] &=& (T(n)-T(n-1)) + (T(n-1)-T(n-2)) + \ldots + (T(2)-T(1))+ T(1)\\
&=& \sum_{r=1}^{n-1}\frac{n-1}{n-r}=(n-1)h(n-1)
\end{eqnarray*}
\proofend 

Thus, $\cov[K_n]\sim n\ln n$.

\subsection{Path}\label{section:pathcovertime}
Without loss of generality, we label the vertices of the $n$-vertex path graph, $P_{n}$ by the set $[0,n-1]$, where labels are given in order starting from one end, i.e., $P_n = (0,1,2,\ldots,n-1)$.
\\
 
\begin{theorem}\label{theorem:pathCover}
Let $P_{n}$ be the path graph on $n$ vertices. 
\begin{description}
	\item[(i)] $\fretT[0] = \fretT[n-1]= 2(n-1)$ and $\fretT[i]=n-1$ for $0<i<n-1$. 
	\item[(ii)] When $i< j$, $\hit[i,j] = j^2-i^2$. In particular, $\hit[0,n-1] = (n-1)^2$. 
	\item[(iii)] $\cov[P_{n}] = \left\{ 
  \begin{array}{l l}
    \frac{5(n-1)^2}{4} & \quad \text{if $n$ odd}\\
   \frac{5(n-1)^2}{4}-\frac{1}{4} & \quad \text{if $n$ even}
  \end{array} \right.$ 
\end{description}
\end{theorem}

\proofstart
\textbf{(i)} $d(0) = d(n-1) = 1$, and $d(i)=2$ when $0<i<n-1$. Furthermore, $|E| = n-1$. Now use Theorem \ref{Th:RandomWalkStationDist} and  Theorem \ref{ergodicTheorem}.

\textbf{(ii)} By part (i), $\fretT[n-1] = 2(n-1)$, but we also know that $\fretT[n-1]= 1+\hit[n-2,n-1]$, because when the walk is on vertex $n-1$, it has no choice but to move to vertex $n-2$ in the next step. So we have $\hit[n-2,n-1] = 2(n-1)-1$. Furthermore, $\hit[r-1,r] = 2r-1$ for any $0<r\leq n-1$, since this is the same as $\hit[r-1,r]$ on  $P_{r+1}$. We have, by linearity of expectation,
\[
\hit[i,j] = \sum_{r=i+1}^j\hit[r-1, r] = \sum_{r=i+1}^j 2r-1 = (j+i+1)(j-i) - (j-i) = j^2-i^2.
\]

\textbf{(iii)} Our analysis will be more notationally convenient with $n+1$ vertices than with $n$ vertices. Let $a_r$ denote the expected time to reach either one of the ends, when starting on some vertex $r\in [0,n]$. The variables $a_r$ satisfy the following system of equations:
\begin{equation}
  a_r = \left\{ 
  \begin{array}{l l}
    0 & \quad \text{if }r=0\\
    1+\frac{1}{2}a_{r-1} + \frac{1}{2}a_{r+1} & \quad \text{if } 0< r <n\\
    0 & \quad \text{if }r=n\\
  \end{array} \right.\label{syseq}
\end{equation}
A solution is $a_r = r(n-r)$: $a_0 = 0(n-0) = 0 = n(n-n) = a_n$, and 
\begin{align*}
&a_r - 1 - \frac{1}{2}a_{r-1} - \frac{1}{2}a_{r+1}\\
&\quad = r(n-r)- 1 - \frac{1}{2}[(r-1)(n-(r-1)) + (r+1)(n-(r+1))]\\
\ignore{&\quad =r(n-r) - 1 -\frac{1}{2}[(r-1)n - (r-1)^2 + (r+1)n - (r+1)^2]\\
&\quad=r(n-r) - 1 -\frac{1}{2}[2rn - (r^2-2r+1+r^2+2r+1)]\\
&\quad=rn-r^2 - 1 - rn + r^2+1\\}
& \quad= 0.
\end{align*}
Furthermore, this is the only solution. This can be seen by studying the system of equations written as a matrix equation, and determining that the equations are linearly independent, but a more elegant method relies on the principle that harmonic functions achieve their maximum and minimum on the boundary, as expressed in Lemma \ref{maxminharmlemma}: Suppose that there is some other solution $a'_r$ set the system of equations (\ref{syseq}). Consider $f_r = a_r-a'_r$. We have $f_0 = a_0 - a'_0 = 0$, and $f_n = a_n - a'_n = 0$. Furthermore, for $0<r<n$,
\begin{eqnarray*}
f_r &=& a_r-a'_r \\
		&=& 1+\frac{1}{2}a_{r-1} + \frac{1}{2}a_{r+1}-\left(1+\frac{1}{2}a'_{r-1} + \frac{1}{2}a'_{r+1}\right)\\
 		&=& \frac{1}{2}(a_{r-1}-a'_{r-1}) + \frac{1}{2}(a_{r+1}-a'_{r+1})\\
 		&=& \frac{1}{2}f_{r-1} + \frac{1}{2}f_{r+1}.
\end{eqnarray*}
Thus, $f_r$ is harmonic on $[1,n-1]$. Since $f_r=0$ on the boundary vertices $0,n$, by Lemma \ref{maxminharmlemma} $f_r=0$ for $r \in [1,n-1]$, and so $a'=a$.

To cover $P_{n+1}$ when starting at vertex $r\in [0,n]$, the walk needs to reach either one of the ends, then make its way to the other. Thus
\[
\cov_r[P_{n+1}] = a_r+ \hit[0,n] =  r(n-r) + n^2.
\]
When $n$ is even (i.e, the path $P_{n+1}$ has an odd number of vertices), then $\cov_r[P_{n+1}]$ is maximised at $r=\frac{n}{2}$, in which case, 
\begin{equation*}
\cov_r[P_{n+1}] = \cov[P_{n+1}] = \frac{5n^2}{4}. 
\end{equation*}
When $n$ is odd, $\cov_r[P_{n+1}]$ is maximised at $r=\left\lfloor \frac{n}{2} \right\rfloor$ or $r=\left\lceil \frac{n}{2} \right\rceil$ , in which case,
\begin{eqnarray*}
\cov_r[P_{n+1}] = \cov[P_{n+1}] &=& \left\lceil\frac{n}{2}\right\rceil\left(n-\left\lceil\frac{n}{2}\right\rceil\right) + n^2\\
&=&\left(\frac{n}{2}+\frac{1}{2}\right)\left(\frac{n}{2}-\frac{1}{2}\right)+n^2\\
&=&\frac{5n^2}{4}-\frac{1}{4}.
\end{eqnarray*}
\proofend

\subsection{Cycle}\label{section:cyclecovertime}

\begin{theorem}\label{theorem:cycleCover}
Let $\mathbb{Z}_n$ be the cycle graph on $n$ vertices. 
\begin{description}
	\item[(i)] $\fretT[u]=n$ for any vertex $u$.
	\item[(ii)] For a pair of vertices $u,v$ distance $r\leq n/2$ from each other on $\mathbb{Z}_n$, $\hit[u,v] = r(n-r)$,	
	\item[(iii)] $\cov[\mathbb{Z}_n] = \frac{n(n-1)}{2}$.
\end{description}
\end{theorem}

\proofstart
\textbf{(i)} The vertices all have the same degree, so by symmetry, and in conjunction with Theorem \ref{Th:RandomWalkStationDist}, $\pi_u$ must be the same for all $u\in V$, that is, $\pi_u = 1/n$. Now apply Theorem \ref{ergodicTheorem}.

\textbf{(ii)} We use the principles of the proof of Theorem \ref{theorem:pathCover}. We assume the vertices of $\mathbb{Z}_n$ are labelled with $[0,n-1]$ in order around the cycle. Hence, vertex $0$, for example, would have vertices $1$ and $n-1$ as neighbours. We wish to calculate $\hit[r,0]$. Observe that there are two paths from $r$ to $0$; one path is $(r, r-1, \ldots, 0)$. On this path the distance between vertices $r$ and $0$ is $r$. The other path is $(r, r+1, \ldots, n-1, 0)$. The distance between $r$ and $0$ on this path is $(n-r)$. By the same principles as the proof of Theorem \ref{theorem:pathCover}, we calculate $\hit[r,0]$ by equating it with the expected time it takes a walk to reach $0$ or $n$ of a path graph $(0,1,2,\ldots,n)$, when the walk starts at vertex $r$. As calculated for Theorem \ref{theorem:pathCover}, this is $r(n-r)$.

(iii) To determine the cover time, observe that at any point during the walk, the set of vertices that have been visited will be contiguous on the cycle; there will be an ``arc'' (path) of visited vertices, and another of unvisited. Let $T(r)$ be expected time it takes a walk starting at some vertex to visit $r$ vertices of $\mathbb{Z}_n$. Given that the walk does indeed start on some vertex, we have $T(1)=0$. After it moves for the first time, it visits a new vertex, thus giving $T(2)-T(1)=1$. By linearity of expectation, 

\begin{equation}
\cov[\mathbb{Z}_n] = (T(n)-T(n-1))+ (T(n-1)-T(n-2)) + \ldots + (T(2)-T(1))+ T(1).\label{eq:cyclcovtime}
\end{equation}

Suppose the walk has just visited the $r$'th new vertex, where $r<n$. Without loss of generality, we can label that vertex $r$, and further label the arc of already visited neighbours $(r-1, r-2, \ldots, 1)$, in order. The arc of unvisited vertices is labelled $(r+1, r+2, \ldots, n-1, 0)$, such that $(r,r+1)$ and $(0,1)$ are edges on the cycle. Thus, the next time the walk visits a new vertex, it will be either the vertex labelled $r+1$ or the vertex labelled $0$ in the current labelling. Hence, $T(r+1)-T(r)$ is the same as the expected time it takes a walk on a path graph $(0,1,2,\ldots,r+1)$ to reach $0$ or $r+1$, when it starts on vertex $r$. As calculated above, this is $r(r+1-r) = r$. Equation (\ref{eq:cyclcovtime}) can thus be calculated as 
\[
\cov[\mathbb{Z}_n] = \sum_{r=0}^{n-1}r = \frac{n(n-1)}{2}.
\]
\proofend

\section{General Bounds and Methods}\label{General Bounds and Methods}
In this section, we detail two general approaches for bounding cover times: the \emph{spanning tree technique}, and the \emph{Matthews's technique}. Despite the simplicity of the techniques, they can often yield bounds that are within constant factors of the actual cover time. Both methods can be applied to a graph under question, but it is often the case that one is more suited, i.e, yields tighter bounds - than the other for a particular graph.  In both cases, the effectiveness of the technique is dependent on finding suitable bounds on hitting times between vertices (or sets of vertices), as well as the way in which the technique is applied.

\subsection{Upper Bound: Spanning Tree and First Return Time}\label{spantree}
Let  $G=(V,E)$ be an undirected, unweighted, simple, connected graph. Let $n=|V|$ and $m=|E|$. One way to upper bound the cover time of $G$ is to choose some sequence of vertices $\sigma=(v_0, v_1, \ldots, v_r)$ such that every vertex in $V$ is in $\sigma$, and sum the hitting time from one vertex to another in the sequence; that is,
\[
\cov[G] \leq \hit[v_0, v_1] + \hit[v_1, v_2] + \ldots + \hit[v_{r-1}, v_r].
\]
\\

\begin{proposition}\label{PropEulerTour}
For a tree $T =(V,E)$ we can generate a walk (sequence of edge transitions on $T$) $\sigma=(v_0, v_1, \ldots, v_{2|V|-2})$ such that each edge of $T$ is traversed once in each direction.
\end{proposition}

The sequence $\sigma$ of Proposition \ref{PropEulerTour} contains every vertex of $V$. It can, in fact, be generated be the \emph{depth first search} (DFS) algorithm started at some vertex $v=v_0$. We shall use DFS again in chapter \ref{ch:The Cover Time of Cartesian Product Graphs}. See, e.g.,\ \cite{Kleinberg} for a discussion of the algorithm.
\\

The following theorem and proof are given given in \cite{mitz}. The argument itself goes back to \cite{AKLLR}.
\\

\begin{theorem}\label{Thm:spanningTreeFirstRet}
Let  $G=(V,E)$ be an undirected, unweighted, simple, connected graph, and let $n=|V|$ and $m=|E|$. 
\[
\cov[G]<4mn.
\]
\end{theorem}
\proofstart

For any $v\in V$, using Theorem \ref{Th:RandomWalkStationDist} and  Theorem \ref{ergodicTheorem}, we have 
\begin{equation}
\fretT[v] =\frac{2m}{d(v)}.\label{eq:hitndksjnks}
\end{equation}
But we also know that

\begin{equation}
\fretT[v] = \sum_{u\in N(v)}\frac{1}{d(v)}\left(1+\hit[u,v]\right)\label{eq:hitjdanfsf}
\end{equation}
where $N(v)$ is the neighbour set of $v$.
Thus, equating (\ref{eq:hitndksjnks}) and (\ref{eq:hitjdanfsf}),
\[
\frac{2m}{d(v)} = \frac{1}{d(v)}\sum_{u\in N(v)}1+\hit[u,v]
\]
and so $\hit[u,v]<2m $.

Let $T =(V,E_T)$ be some spanning tree of $G$. Let $\sigma=(v_0, v_1, \ldots, v_{2n-2})$ be a walk as described in Proposition \ref{PropEulerTour}. Since each vertex of $V$ occurs in $\sigma$, we have
\[
\cov[G] \leq \sum_{i=0}^{2n-3} \hit[v_{i}, v_{i+1}] \leq  2m(2n-2)<4mn.
\]  
\proofend

\subsection{Upper Bound: Minimum Effective Resistance Spanning Tree}\label{MERSTArg}
In section \ref{spantree}, we computed an upper bound on the sum of commute times of a spanning tree of $G$. We can generalise this to trees that \emph{span} $G$, that is, include all the vertices of $G$, but who's edges are not necessarily contained in $G$.

Let $G$ be a connected, undirected graph. If $G$ is unweighted, assign unit weights (conductances) to the edges of $G$. Thus, $G=(V,E,c)$.
\\

\begin{definition}
Define the complete graph $K = (V, E', \rho)$ where the weighting function $\rho : E' \rightarrow \mathbb{R^+}$, and $\rho((u,v)) = R(u,v)$ where $R(u,v)$ is the effective resistance between $u$ and $v$ in $G$. Let $\mathcal{T} = \{T \text{ is a spanning tree of $K$}\}$, and let $T_* \in \mathcal{T}$ be such that $w(T_*) \leq w(T)$ for any $T\in \mathcal{T}$  (recall $w(T)$ is the total of the edge weights of $T$, as defined by equation \ref{w(G)}). We call $T_*$ the \emph{minimum effective resistance spanning tree} of $G$.
\end{definition}
\vspace{\baselineskip}
\begin{theorem}\label{MERSTcovBound}
Let $G$ be a connected, undirected graph. If $G$ is unweighted, assign unit weights (conductances) to the edges of $G$. Thus, $G=(V,E,c)$.
\begin{equation*}
\cov[G] \leq c(G)w(T_*),
\end{equation*}
where $T_*$ is the minimum effective resistance spanning tree of $G$.
\end{theorem} 
\proofstart
Using Theorem \ref{Th:commtetimetheorem}, we have for any $u,v\in V$, 
\[
\com[u,v] = c(G)R(u,v).
\] 
So
\begin{eqnarray}
c(G)w(T_*) = c(G)\sum_{(u,v)\in E'}R(u,v) = \sum_{(u,v)\in E'}\com[u,v]. \label{eq:isnjhcbj}
\end{eqnarray}
Now we continue using the same ideas of the proof of Theorem \ref{Thm:spanningTreeFirstRet}: we apply Proposition \ref{PropEulerTour} to generate a sequence $\sigma$ that transition each edge of $T_*$ once in each direction, thereby visiting every vertex of $G$. The RHS of the second equality of (\ref{eq:isnjhcbj}), $\sum_{(u,v)\in E'}\com[u,v]$, is the sum of hitting times for the sequence $\sigma$.
\proofend

For all but a few simple examples, it can be difficult to determine $w(T_*)$. However, bounds on effective resistances can often be determined using the various tools of electrical network theory; for example, through the use of flows and Thomson's principle (Theorem \ref{Thomson's Principle}), and other tools such as Rayleigh's laws, cutting and shorting laws, etc.

\subsection{Upper Bound: Matthews' Technique}\label{Upperbound: Matthews' Technique}
\begin{theorem}[\textbf{Matthews' upper bound}, \cite{Matthews}]
\label{MatLemma1} For a graph $G=(V,E)$,
\begin{equation}
\cov[G] \leq \hit^*[G]h(n), \label{mattUpper}
\end{equation}
where $\hit^*[G] = \max_{u,v\in V}\hit[u,v]$ and $h(n)$ is the $n$'th harmonic number $\sum_{i=1}^n\frac{1}{n}$.
\end{theorem}
We refer the reader to, e.g.,\ \cite{LPW} for a proof. The proof is not difficult, but is fairly lengthy.

The power of the method is two-fold. Firstly, one needs only to bound $\hit^*[G]$, which can be facilitated through electrical network theory as well as consideration of the structure of $G$. Secondly, it applies also to weighted graphs (note there is no restriction of being unweighted in the statement of the theorem).

To use electrical network theory, we can use the commute time identity of Theorem \ref{Th:commtetimetheorem}, and use the commute time as an upper bound for hitting time. 
In this case, (\ref{mattUpper}) is expressible as  $\cov[G] \leq c(G)R^*(G)h_n$, where $R^*(G)$ is the maximum effective resistance between any pair of vertices in $G$.

Despite the simplicity of the inequality, the method can yield bounds on cover time that are within a constant factor of the precise value. This is always the case if $\hit^*[G]$ can be shown to be $O(n)$, since this gives a cover time of $O(n\log n)$, and as we shall see in section \ref{FeigeSection}, cover times of graphs are $\Omega(n\log n)$. 

One example of an application of Theorem \ref{MatLemma1} that gives good bounds is on the complete graph. As was established in Theorem \ref{theorem:completecovertime}, has $\hit[u,v] = n-1$ for any pair $u,v$. The resulting bound is very close to the precise result. 

\subsubsection{Matthews' Technique for a Subset}
The inequality (\ref{mattUpper}) bounds the cover time of all the vertices of the graph, but it applies equally to a subset of the vertices $V'\subseteq V(G)$: 
\\

\begin{theorem}[\textbf{Matthews' bound, subset version}]\label{subsMatTheorem}
Let 
\[
\hit_{G}^*[V']=\max\{\hit[u,v] :  u,v \in V'\},
\] 
where $\hit[u,v]$ is the hitting time from $u$ to $v$ in $G$. For a random walk on $G$ starting at some vertex $v\in V'$, denote by $\cov_v[V']$ the expected time to visit all the vertices of $V'$. Then
\begin{equation}
\cov_v[V'] \leq \hit_{G}^*[V']h(|V'|), \label{subsmattUpper}
\end{equation}
\end{theorem}
The notation $\cov[V']$ shall mean $\max_{v\in V'}\cov_v[V']$. 

\subsection{Lower Bound: Matthews' Technique}
There is also a lower bound version of Matthews' technique. A proof is given in \cite{LPW}.
\\
\begin{theorem}[\textbf{Matthews' lower bound} \cite{Matthews}]\label{mattLower}
For the graph $G=(V,E)$, 
\[
\cov[G] \geq \max_{A \subseteq V}\hit_*[A]h(|A|-1).
\]
where $\hit_*[A] = \min_{u,v\in A, u\neq v}\hit[u,v]$.
\end{theorem}
When using Theorem \ref{mattLower}, one needs to be careful with the choice of $A$; the end of a path has hitting time $1$ to its neighbour. Thus, if both of these vertices are included in $A$, the bound is no better than $\log n$, and of course, cover time is at least $n$.

A refinement on Theorem \ref{mattLower} was given in \cite{ZuLower}. We shall quote it as in \cite{Feige2}, since this is simpler notation.
\\

\begin{lemma}[\cite{ZuLower}]\label{ZuLowerLemma}
Let $S$ be a subset of vertices of $G$ and let $t$ be such that for all $v\in S$, at most $b$ of the vertices $u \in S$ satisfy $\hit[v,u]<t$. Then 
\[
\cov_v[G]\geq t(\log(|S|/b)-2).
\]
\end{lemma}

\section{General Cover Time Bounds}
\subsection{Asymptotic General Bounds}\label{FeigeSection}
In two seminal papers on the subject of cover time, Feige gave tight asymptotic bounds on the cover time. As usual, the logarithm is base-$e$ unless otherwise stated.
\\

\begin{theorem}[\cite{Feige2}]\label{FeigeLower}
For any graph $G$ on $n$ vertices and any starting vertex $u$
\[
\cov_u[G] \geq (1+o(1))n\log n
\]
\end{theorem}
This lower bound is exhibited by $K_n$ the complete graph on $n$ vertices, as was demonstrated in Theorem \ref{theorem:completecovertime}. 
It is proven using Lemma \ref{ZuLowerLemma}. It is shown that at least one of the two conditions hold for any graph.
\begin{enumerate}
	\item There are two vertices $u$ and $v$ such that $\hit[u,v] \geq n\log n$ and $\hit[v,u] \geq n\log n$.
	\item The assumptions of Lemma \ref{ZuLowerLemma} hold with parameters  $|S|>n/(\log_2 n)^c$, $b<(\log_2 n)^c$ and $t\geq n(1-c/\log_2 n)$ for some constant $c$ independent of $n$.
\end{enumerate}
The upper bound is as follows
\\

\begin{theorem}[\cite{Feige1}]\label{FeigeUpper}
For any connected graph $G$ on $n$ vertices
\begin{equation}
\mathbf{CyCOV}[G] \leq (1+o(1))\frac{4}{27}n^3 \label{cycCoverUB}
\end{equation}
\end{theorem}
The quantity $\mathbf{CyCOV}[G]$ is the \emph{cyclic cover time}, which is the expected time it takes to visit all the vertices of the graph in a  specified cyclic order, minimised over all cyclic orders. Clearly, $\cov[G]< \mathbf{CyCOV}[G]$. 

Consider the \emph{lollipop graph}, which is a path of length $n/3$ connected to a complete graph of $2n/3$ vertices. Let $u$ be the vertex that connects the clique to the path, and $v$ be the vertex at the other end of the path. It can be determined that $\hit[u,v]=\com[u,v]-\hit[v,u]=(1+o(1))\frac{4}{27}n^3$. This can be seen by applying Theorem \ref{Th:commtetimetheorem} to get $\com[u,v]$ and Theorem \ref{theorem:pathCover} to get $\hit[v,u]$.  This demonstrates that the asymptotic values of maximum hitting time and cyclic cover time (and therefore also cover time) can be equal (up to lower order terms) even when the cyclic cover time is maximal (up to lower order terms). 

Theorem \ref{FeigeUpper} is proved by a contradiction argument on the minimum effective resistance spanning tree $T_*$ (see section \ref{MERSTArg}). A trade off is demonstrated between $w(T_*)$, the sum of effective resistances of edges of $T_*$,  and $m$, the number of edges of the graph. This bounds the cover time as per Theorem \ref{MERSTcovBound}.

\chapter{The Cover Time of Cartesian Product Graphs}\label{ch:The Cover Time of Cartesian Product Graphs}
In this chapter, we study the cover time of random walks on on the Cartesian product $F$ of two graphs $G$ and $H$. In doing so, we develop a relation between the cover time of $F$ and the cover times of $G$ and $H$. When one of $G$ or $H$ is in some sense larger than the other, its cover time dominates, and can become of the same order as the cover time of the product as a whole. Our main theorem effectively gives conditions for when this holds. The probabilistic technique which we introduce, based on a quantity called the \emph{blanket time}, is more general and may be of independent interest, as might some of the lemmas developed in this chapter. The electrical network metaphor is one of the principle tools used in our analysis. 

$G$ and $H$ are assumed to be finite (as all graphs in this thesis are), undirected, unweighted, simple and connected.

\section{Cartesian Product of Graphs: Definition, Properties, Examples}
\subsection{Definition}

\begin{definition}
\label{DefCartProd} Let $G=(V_G, E_G)$ and $H = (V_H, E_H)$ be simple, connected, undirected graphs. The \emph{Cartesian product}, $G \Box H$ of $G$ and $H$ is the graph $F=(V_F, E_F)$ such
that
\begin{description}
    \item[(i)] $V_F = V_G\times V_H$
    \item[(ii)] $((a,x), (b,y)) \in E_F$ if and only if either
    \begin{enumerate}
        \item $(a,b) \in E_G$ and $x=y$, or
        \item $a=b$ and $(x,y) \in E_H$
    \end{enumerate}
\end{description}
\end{definition}

We call $G$ and $H$ the \emph{factors} of $F$, and we say that $G$ and $H$ are \emph{multiplied} together.

We can think of $F = G \Box H$ in terms of the following
construction: We make a copy of one of the graphs, say $G$, once for
each vertex of the other, $H$. Denote the copy of $G$ corresponding to
vertex $x \in V_H$ by $G_x$. Let $a_x$ denote a vertex in $G_x$ corresponding to $a \in V_G$. If there is an edge $(x,y) \in E_H$, then add an edge $(a_x, a_y)$ to the construction.

\textbf{Notation}
For a graph $\Gamma = (V_{\Gamma}, E_{\Gamma})$, denote by 
\begin{description}
	\item[(i)]$n_{\Gamma}$ the number of vertices $|V_{\Gamma}|$, and 
	\item[(ii)]$m_{\Gamma}$ the number of edges $|E_{\Gamma}|$. 
\end{description}	
In addition we will use use the notation $N$ and $M$ to stand for $n_F$ and $m_F$ respectively.

\subsection{Properties}\label{properties}
\textbf{Commutativity of the Cartesian Product Operation}
For a pair of graphs $G$ and $H,$ $G \Box H$ is isomorphic to $H \Box G$; that is, if vertex labels are ignored, the graphs are identical. Note, however, that by (i) of Definition \ref{DefCartProd}, the two different orders on the product operation do produce different labellings.

\textbf{Vertices and Edges of the Product Graph}
The number of vertices and edges of a Cartesian product is related to the vertices and edges of its factors as follows:

\begin{description}
	\item[(i)] $N =n_Gn_H$.
	\item[(ii)] $M = n_Gm_H + n_Hm_G$.
\end{description}

(i) follows from the properties of Cartesian product of two sets, and (i) of Definition \ref{DefCartProd}. To see (ii), we have by (ii)--1 of Definition \ref{DefCartProd} the following: For a vertex $x \in V_H$, there is, for each $(a,b)\in E_G$, an edge $((a,x), (b,x))\in E_F$. That is, we have the set  
\[
S_{H,x} = \{((a,x), (b,x)) : (a,b)\in E_G\} 
\]
Similarly, we have by (ii)--2 of Definition \ref{DefCartProd} the following: For a vertex $a \in V_G$, there is, for each $(x,y)\in E_H$, an edge $((a,x), (a,y))\in E_F$. That is, we have the set  
\[
S_{G,a} = \{((a,x), (a,y)) : (x,y)\in E_H\} 
\]
Thus, 
\[
E_F = \bigcup_{x \in V_H}S_{H,x} \cup \bigcup_{a \in V_G}S_{G,a}.
\]
Now $|S_{H,x}| = m_G$ for all $x \in V_H$, and $|S_{G,a}| = m_H$ for all $a \in V_G$, and since the sets are all disjoint, we have 
\begin{eqnarray*}
M = |E_F| &=&  \left|\bigcup_{x \in V_H}S_{H,x} \cup \bigcup_{a \in V_G}S_{G,a}\right| \\
&=& \left|\bigcup_{x \in V_H}S_{H,x}\right| + \left|\bigcup_{a \in V_G}S_{G,a}\right|\\
&=& \sum_{x \in V_H}m_G + \sum_{a \in V_G}m_H\\
&=& n_Hm_G + n_Gm_H
\end{eqnarray*}

\textbf{Associativity of the Cartesian Product and a Generalisation to an Arbitrary Number of Factors}
We can extend the definition of the Cartesian product to an arbitrary number of factors: $(G_1 \Box G_2) \Box G_3)\ldots...)\Box G_r$. If we always represent the resulting product vertices by $r$-tuples and edges by pairs of $r$-tuples, then any bracketing in which a bracket contains a product of two graphs (either of which may be a product itself) will give the same product, that is, the Cartesian product is associative. We can therefore represent it unambiguously by $F=G_1\Box G_2 \ldots \Box G_r$. This assumes the order of the operands is kept the same. If the order is permuted, the tuple representing the vertex labelling will be permuted in the same way, but the two permutations will produce isomorphic products.

For a natural number $d$, we denote by $G^d$ the $d$'th Cartesian power, that is, 
$G^d = G$ when $d=1$ and $G^d = G^{d-1}\Box G$ when $d>1$.

\subsection{Examples}
We give examples of Cartesian product of graphs, some of which are important to the proofs of this chapter.
First, we remind the reader of some specific classes of graphs, and define new ones: Let $P_n$ denote the \emph{$n$-path}, the path graph of $n$ vertices.
Let $\mathbb{Z}_n$ represent the \emph{$n$-cycle}, the cycle graph with $n$ vertices.
The Cartesian product of a pair of paths, $P_p \Box P_q$ is a $p \times q$ \emph{rectangular grid}, and when $p=q=n$, is a $n \times n$ \emph{grid}, or \emph{lattice}. The  product of a pair of cycles $\mathbb{Z}_p \Box \mathbb{Z}_q$ is a \emph{toroid}, and when  $p=q$, is a \emph{torus}. Both grids and toroids can be generalised to higher powers in the obvious way to give $d-$dimensional grids and toroids respectively, where $d$ is the number of paths or cycles multiplied together, respectively. 

\begin{figure}[htb]
\begin{center}
\leavevmode
\includegraphics{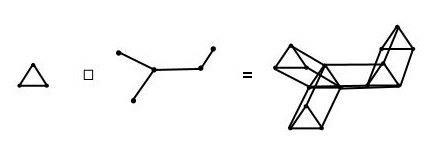}
\end{center}
\caption{Cartesian product of a triangle with a tree.}
\label{fig:CartProdTriangleTree}
\end{figure}

To give another - somewhat more arbitrary - example, a pictorial representation of the product of a triangle graph with a tree is given in Figure \ref{fig:CartProdTriangleTree}.

\section{Blanket Time}
\label{blanketTime}
We introduce here a notion that is related to the cover time, and is an important part of the main theorem and the proof technique we use.
\\

\begin{definition}[\cite{Winkler}]
\label{DefBlanket}
For a random walk $\mathcal{W}_u$ on a graph $G=(V,E)$ starting at some vertex $u \in V$, and $\delta \in [0,1)$, define the the random variable
\begin{equation}
B_{\delta, u}[G] = \min\{t: \forall v \in V, N_v(t)>\delta\pi_v t\}, \label{blanketRV}
\end{equation}
where $N_v(t)$ is the number of times $\mathcal{W}_u$ has visited $v$ by time $t$ and $\pi_v$ is the stationary probability of
vertex $v$. The \emph{blanket time} is 
\[
\blank_{\delta}[G] = \max_{u \in V}\E[B_{\delta, u}[G]].
\]
\end{definition}

The following was recently proved in \cite{Ding}.
\\

\begin{theorem}[\cite{Ding}]
\label{BlanketTheorem}
For any graph $G$, and any $\delta \in (0,1)$, we have
\begin{equation}
    \blank_{\delta}[G]\leq \kappa(\delta)\cov[G]
\end{equation}
Where the constant $\kappa(\delta)$ depends only on
$\delta$.
\end{theorem}

We define the following
\\

\begin{definition}[\textbf{Blanket-Cover Time}]\label{Blanket Cover def}
For a random walk $\mathcal{W}_u$ on a graph $G=(V,E)$ starting at some vertex $u \in V$, define the the random variable 
\[
\beta_u[G] =\min\{t: \forall v \in V\, N_v(t)\geq\pi_v\cov[G]\},
\] 
where $N_v(t)$ is the number of times $\mathcal{W}_u$ has visited $v$ by time $t$ and $\pi_v$ is the stationary probability of
vertex $v$. The \emph{blanket-cover time} is the quantity 
\[
\bcov[G] =\max_{u \in V}\E[\beta_u[G]].
\] 
\end{definition}

Thus the blanket-cover time of a graph is the expected first time at each vertex $v$ is visited at least $\pi_v\cov[G]$ times -
which we shall refer to as the \emph{blanket-cover criterion}.

In the paper that introduced the blanket time, \cite{Winkler}, the following equivalence was asserted, which we conjecture to be true.
\\

\begin{conjecture}\label{Blanket Cover Criterion Equiv}
 $\bcov[G] = O(\cov[G])$.
\end{conjecture}

In the same paper, this equivalence was proved for paths and cycles. However, we have not found a proof for the more general case. It can be shown without much difficulty that $\bcov[G] = O((\cov[G])^2)$. Using the following lemma, we can improve upon this.
\\

\begin{lemma}[\cite{LovaszBlanket}]\label{LovaszBlanketLem}
Let $i$ and $j$ be two vertices and $k\geq 1$. Let $W_k$ be the number of times $j$ had been visited
when $i$ was visited the $k$-th time. Then for every $\varepsilon>0$,
\[
\Pr\left(W_k<(1-\varepsilon)\frac{\pi_j}{\pi_i}k\right)\leq \exp \left(\frac{-\varepsilon^2k}{4\pi_i\com[i,j]}\right).
\]
\end{lemma}
We use it thus:
\\

\begin{lemma}
\[
\bcov[G] = O\left((\log n)\com^*[G]\right) 
\]
where $\com^*[G]=\max_{u,v \in V(G)}\com[u,v]$. 
\end{lemma}
\proofstart
At time $t$ some vertex $i$ must have been visited at least $\pi_it$ times, otherwise we would get $t= \sum_{v\in V}N_v(t) < \sum_{v\in V}\pi_vt = t$, where $N_v(t)$ is the number of times $v$ has been visited by time $t$.

We let the walk run for $\tau =A(\log n)\com^*[G]$ steps where $A$ is a large constant. Some vertex $i$ will have been visited at least $\pi_i\tau$ times. Now we use Lemma \ref{LovaszBlanketLem} with $k=\pi_i\tau$. Then for any $j$,
\begin{eqnarray*}
\Pr\left(W_k<(1-\varepsilon)\frac{\pi_j}{\pi_i}k\right)&\leq& \exp \left(\frac{-\varepsilon^2k}{4\pi_i\com[i,j]}\right)\\ 
&\leq& \exp \left(\frac{-\varepsilon^2A\log n}{4}\right)\\
&\leq& 1/n^c
\end{eqnarray*}
for some constant $c>1$. Hence with probability at most $1/n^{c-1}$ the walk has failed to visit each vertex $j$ at least $\pi_j\cov[G]$ times (by Matthews' bound, Theorem \ref{MatLemma1}). We repeat the process until success. The expected number of attempts is $1+O(n^{1-c})$. 

\proofend
\section{Relating the Cover Time of the Cartesian Product to Properties of its Factors}

\textbf{Notation} For a graph $\Gamma$, denote by: $\delta_{\Gamma}$ the minimum degree; $\theta_{\Gamma}$ the average degree; $\Delta_{\Gamma}$ the maximum degree; $D_{\Gamma}$ the diameter.  

The main theorem of this chapter is the following.
\\

\begin{theorem}
\label{main}
Let $F= (V_F, E_F) = G \Box H$ where $G=(V_G, E_G)$ and $H = (V_H, E_H)$ are simple, connected, unweighted, undirected graphs. 
We have
\begin{equation}
\cov[F] \geq \max\left\{\left(1+\frac{\delta_G}{\Delta_H}\right) \cov[H], \left(1+\frac{\delta_H}{\Delta_G}\right) \cov[G]\right\}. \label{LBIneq}
\end{equation}
Suppose further that $n_H \geq D_G+1$, then
\begin{equation}
\cov[F] \leq K\left(\left(1+\frac{\Delta_G}{\delta_H}\right)\bcov[H] + \frac{Mm_Gm_Hn_H\ell^2}{\cov[H]D_G}\right)\label{UBIneq}
\end{equation}
where $M =|E_F|= n_Gm_H + n_Hm_G$, $\ell=\log (D_G+1) \log(n_GD_G)$ and $K$ is some universal constant. 
\end{theorem}

The main part of the work is the derivation of (\ref{UBIneq}); the inequality (\ref{LBIneq}) is relatively straightforward to derive. Note, by the commutativity of the Cartesian product, $G$ and $H$ in the may be swapped in \eqref{UBIneq}, subject to the condition $n_G \geq D_H+1$.

Theorem \ref{main} extends much work done on the particular case of the two-dimensional toroid on $n^2$ vertices, that is, $\mathbb{Z}^2_n = \mathbb{Z}_n \Box \mathbb{Z}_n$, culminating in a result of \cite{Dembo}, which gives a tight asymptotic result for the cover time of $\mathbb{Z}^2_n$ as $n \rightarrow \infty$. Theorem \ref{main} also extends work done in \cite{Jonasson} on powers $G^d$ of general graphs $G$, which gives upper bounds for the cover time of powers of graphs. Specifically, it shows $\cov[G^2] = O(\theta_G N\log^2 N)$ and
for $d \geq 3$, $\cov[G^d] = O(\theta_G N\log N)$. Here $N=n^d$, is the number of vertices in the product, and $\theta_G=2|E|/n$ is the average degree of $G$. A formal statement of the theorem is given in section \ref{RelatedWork} and further comparisons made in section \ref{Cover Time: Examples and Comparisons}.

To prove the Theorem \ref{main}, we present a framework to bound the cover time of a random walk on a graph which works by dividing the graph up into (possibly overlapping) regions, analysing the behaviour of the walk when locally observed on those regions, and then composing the analysis of all the regions over the whole graph. Thus the analysis of the whole graph is reduced to the analysis of outcomes on local regions and subsequent compositions of those outcomes. This framework can be applied more generally than Cartesian products.

Some of the lemmas we use may be of independent interest. In particular, Lemmas \ref{lg1x} and \ref{lg2x} provide bounds on effective resistances of graph products that extend well-known and commonly used bounds for the $n \times n$ grid. 

The lower bound in Theorem \ref{main} implies that $\cov[G\Box H] \geq \cov[H]$ (and $\cov[G\Box H] \geq \cov[G]$), and the upper bound can be
viewed as providing conditions sufficient for $\cov[G\Box H] = O(\bcov[H])$ (or $\cov[G\Box H] = O(\bcov[G])$). For example, since paths and cycles have $\bcov[G] = \Theta(\cov[G])$, then $\cov[\mathbb{Z}_p \Box \mathbb{Z}_q] = \Theta(\cov[\mathbb{Z}_q])= \Theta(q^2)$ subject to the condition $p\log^4 p= O(q)$. Thus for this
example, the lower and upper bounds in Theorem \ref{main} are within a constant
factor.

Before we discuss the proof of Theorem \ref{main} and the framework use to produce it, we discuss related work, and give examples of the application of the theorem to demonstrate how it extends that work.

\section{Related Work}
\label{RelatedWork} A $d$-dimensional torus on $N=n^d$ vertices is
the $d$'th power of an $n$-cycle, $\mathbb{Z}^d_n$. The behaviour of
random walks on this structure is well studied.
\\
 
\begin{theorem}[see, e.g., \cite{LPW}]
\begin{description}
	\item[(i)] $\cov[\mathbb{Z}^2_n] = \Theta(N\log^2
N)$.
	\item[(ii)] $\cov[\mathbb{Z}^d_n] = \Theta(N \log N)$ when
$d \geq 3$.
\end{description}
\end{theorem}
In fact, there is a precise asymptotic value for the $2$-dimensional case.
\\

\begin{theorem}[\cite{Dembo}]\label{Thm:Dembo}
 $\cov[\mathbb{Z}^2_n] \sim
\frac{1}{\pi}N \log^2 N$. 
\end{theorem}

The following result of \cite{Jonasson} gives bounds on the cover time for powers of more general graphs:
\\

\begin{theorem}[\cite{Jonasson}, Theorem 1.2]
\label{JonassonTheorem1.2} { Let $G=(V,E)$ be any connected, finite graph on
$n$ vertices with $\theta_G = 2|E|/n$. Let $d \geq 2$ be an integer
and let $N=n^d$. For $d=2$, $\cov[G^d] = O(\theta_G N\log^2 N)$ and
for $d \geq 3$, $\cov[G^d] = O(\theta_G N\log N)$. These bounds are
tight. }
\end{theorem}

\cite{Jonasson} does not address products of graphs that are different, nor does it seem that the proof techniques
used could be directly extended to deal with such cases. Our proof techniques are different, but both this work and 
\cite{Jonasson} make use of electrical network theory and analysis of subgraphs of the product that are isomorphic to the square grid $P_k\Box P_k$.

A  number of theorems and lemmas related to random walks
and effective resistance between pairs of vertices in graph products
are given in \cite{Bollobas}.
To give the reader a flavour we quote Theorem 1 of that paper, which
is useful as a lemma implicitly in this paper and in the proof of
\cite{Jonasson} Theorem 1.2 to justify the intuition that the
effective resistance is maximised between opposite corners of the
square lattice.
\\

\begin{lemma}[\cite{Bollobas}, Theorem 1]
\label{BollobasTheorem1} { Let $P_n$ be an $n$-vertex path with
endpoints $x$ and $y$. Let $G$ be a graph and let $a$ and $b$ be any
two distinct vertices of G. Consider the graph $G \times P_n$. The
effective resistance R((a,x),(b,v)) is maximised over vertices $v$
of $P_n$ at $v=y$. }
\end{lemma}
For $P^2_n$ this is used twice: 
\[
R((0,0),(r,s)) \leq R((0,0),(n-1,s)) \leq R((0,0),(n-1,n-1)).
\]

\section{Cover Time: Examples and Comparisons}\label{Cover Time: Examples and Comparisons}
In this section, we shall apply Theorem \ref{main} to some examples and make comparisons to established results. 

\subsection{Two-dimensional Torus}
 We shall apply the upper bound of Theorem \ref{main}, 
to the $2$-d torus, $\mathbb{Z}^2_n$: 
\begin{description}
	\item[(i)] $G=H = \mathbb{Z}_n$;
	\item[(ii)] $\Delta_{\mathbb{Z}_n} = \delta_{\mathbb{Z}_n} = 2$;
	\item[(iii)] $m_{\mathbb{Z}_n} = n_{\mathbb{Z}_n} = n$;
	\item[(iv)] $D_{\mathbb{Z}_n} = \lfloor\frac{n}{2}\rfloor$;
	\item[(v)] Thus $M = 2m_{\mathbb{Z}_n}n_{\mathbb{Z}_n} = 2n^2$, and
	\item[(vi)] $\ell=\log (D_G+1) \log(n_GD_G) = \log (\lfloor\frac{n}{2}\rfloor+1) \log(n\lfloor\frac{n}{2}\rfloor)$.
	\item[(vii)] By Theorem \ref{theorem:cycleCover}, $\cov[\mathbb{Z}_n] = \frac{n(n-1)}{2}$.
	\item[{(viii)}] $\bcov[G] = \Theta(\cov[G])$.
\end{description}

\begin{eqnarray*}
\cov[F] &\leq& K\left(\left(1+\frac{\Delta_G}{\delta_H}\right)\bcov[H] + \frac{Mm_Gm_Hn_H\ell^2}{\cov[H]D_G}\right)\\
&=& O\left(n^2 + \frac{2n^2nnn\ell^2}{n^2n}\right)\\
&=&O\left(n^2\ell^2\right)\\
&=&O\left(n^2\log^4n\right)\\
&=&O\left(N\log^4N\right).
\end{eqnarray*}
This is a factor $\log^2N$ out of the actual value $\frac{1}{\pi}N \log^2 N$ of Theorem \ref{Thm:Dembo}. Theorem \ref{JonassonTheorem1.2} gives 
$O(N\log^2N)$ bound. 

\subsection{Two-dimensional Toroid with a Dominating Factor}
Theorem \ref{main} does not cope well with squares, for which Theorem \ref{JonassonTheorem1.2} provides strong bounds. Instead it is more effectively applied to cases where there is some degree of asymmetry between the factors. The previous example bounded $\cov[\mathbb{Z}_p\Box \mathbb{Z}_q]$ for the case where $p=q$. If, however, $p\log^4 p= O(q)$, then we get a stronger result. 

\begin{description}
	\item[(i)] $G \equiv \mathbb{Z}_p$;
	\item[(ii)] $H \equiv \mathbb{Z}_q$
	\item[(iii)]$\Delta_G = \Delta_H = \delta_G = \delta_H=2$;
	\item[(iv)] $\bcov[H]=\Theta(\cov[H])$.
	\item[(v)]  $m_G =n_G = p$;
	\item[(vi)]  $m_H=n_H = q$;
	\item[(vii)] $D_G = \lfloor\frac{p}{2}\rfloor$.
	\item[(viii)] Thus $M=2pq$, and
	\item[(ix)]  $\ell = \log(\lfloor\frac{p}{2}\rfloor+1)\log(p\lfloor\frac{p}{2}\rfloor)$.
	\item[(x)]  By Theorem \ref{theorem:cycleCover}, $\cov[\mathbb{Z}_p] = \frac{p(p-1)}{2}$ and $\cov[\mathbb{Z}_q] = \frac{q(q-1)}{2}$.
\end{description} 

Thus, 
\begin{eqnarray*}
\cov[F] &\leq& K\left(\left(1+\frac{\Delta_G}{\delta_H}\right)\bcov[H] + \frac{Mm_Gm_Hn_H\ell^2}{\cov[H]D_G}\right)\\
&=& O\left(q^2 + \frac{2pqpqq\ell^2}{q^2p}\right)\\
&=&O\left(q^2 + pq\log^4p\right)\\
&=&O(q^2)
\end{eqnarray*}
if $p\log^4 p= O(q)$. 

Comparing this to the lower bound of Theorem \ref{main},

\begin{equation*}
\cov[F] \geq \max\left\{\left(\frac{\delta_G}{\Delta_H}+1\right) \cov[H], \left(\frac{\delta_H}{\Delta_G}+1\right) \cov[G]\right\} 
\end{equation*}
which implies

\begin{equation*}
\cov[F] =\Omega \left(\left(\frac{\delta_G}{\Delta_H}+1\right) \cov[H]\right)  = \Omega(q^2).
\end{equation*}
Thus, Theorem \ref{main} gives upper and lower bounds within a constant a multiple for this example. That is, it tells us $\cov[\mathbb{Z}_p \Box \mathbb{Z}_q] = \Theta(\cov[\mathbb{Z}_q])= \Theta(q^2)$ subject to the condition $p\log^4 p= O(q)$. Looking at it another way, it gives conditions for when the cover time of the product $F=G\Box H$ is within a constant multiple of the cover time of one of it's factors. We describe that factor as the \emph{dominating factor}.

\section{Preliminaries}
\subsection{Some Notation}\label{Notation}
For clarity, and because a
vertex $u$ may be considered in two different graphs, we may use
$d_G(u)$ to explicitly denote the degree of $u$ in graph $G$.

$h(n)$ denotes the $n$'th harmonic number, that is,  $h(n) =
\sum_{i=1}^n1/i$. Note $h(n) = \log n + \gamma + O(1/n)$ where
$\gamma \approx 0.577$. All logarithms in this chapter are base-$e$.

In the notation $(.,y)$, the `.' is a place holder for some unspecified element, which may be different from one tuple to another. For example, if we refer to two vertices $(.,a), (.,b) \in G\Box H[S]$, the first elements of the tuples may or may not be the same, but $(.,a)$, for example, refers to a \emph{particular} vertex, not a set of vertices $\{(x,a) : a\in V(G)\}$.

\subsection{The Square Grid}
 
The $k \times k$ grid graph $P^2_k$, where $P_k$ is the $k$-path, plays an important role in our work.
We shall analyse random walks on subgraphs isomorphic to this
structure. It is well known in the literature (see, e.g. \cite{DoyleSnell}, \cite{LPW}) 
that for any pair of vertices $u,v \in V(P^2_k)$, we have $R(u,v) \leq
C\log k$ where $C$ is some universal constant.
We shall quote part of \cite{Jonasson} Lemma 3.1 in our notation and refer
the reader to the proof there.
\\

\begin{lemma}[\cite{Jonasson}, Lemma 3.1(a)]\label{squaregrideffreslemma}
Let $u$ and $v$ be any two vertices of $P^2_k$. Then $R(u,v) <
8h(k)$, where $h(k)$ is the $k$'th harmonic number.
\end{lemma}

\section{Locally Observed Random Walk}
\label{DefnSubgraphWalk} Let $G=(V,E)$ be a connected, unweighted
(equiv., uniformly weighted) graph. Let $S \subset V$ and let $G[S]$
be the subgraph of $G$ induced by $S$.  Let $B=\{v \in S: \exists x
\not \in S, (v,x) \in E\}$. Call $B$ the \emph{boundary} of $S$, and the
vertices of $V \setminus S$ \emph{exterior vertices}. If $v \in S$ then $d_G(v)$ (the
degree of $v$ in $G$) is partitioned into $d(v,in)=|N(v,in)|=|N(v)
\cap S|$ and $d(v,out)=|N(v,out)|=|N(v) \cap (V\setminus S)|$, (\emph{inside} and
\emph{outside} degree). Here $N(v)$ denotes the neighbour set of $v$.

Let $u,v \in B$. Say that $u,v$ are \emph{exterior-connected}  if there is
a $(u,v)$-path $u ,x_1,...x_k,v$ where $x_i \in V \setminus S, k \geq 1$. Thus
all vertices of the path except $u,v$ are exterior, and the path
contains at least one exterior vertex. Let $A(B)=\{ (u,v):$ $u,v$
are exterior-connected $\}$. Note $A(B)$ may include self-loops.

Call edges  of $G[S]$ \emph{interior}, edges of $A(B)$ \emph{exterior}. We say
that a walk $\omega =(u ,x_1,...x_k,v)$ on $G$ is an exterior walk
if $u,v \in S$ and $x_i \notin S$, $1 \leq i \leq k$.

We derive a weighted multi-graph $H$ from $G$ and $S$ as follows:
$V(H)=S$, $E(H)= E(G[S]) \cup A(B)$. Note if $u,v \in B$ and $(u,v) \in E$
then $(u,v) \in E(G[S])$, and if, furthermore, $u,v$ are exterior connected,
then $(u,v) \in A(B)$ \emph{and these edges are distinct}, hence, $H$ may not only
have self-loops but also parallel edges, i.e., $E(H)$ is a multiset.

Associate with an orientation
$(\vec{u,v})$ of an edge $(u,v) \in A(B)$ the set of all exterior
walks $\omega =(u ,x_1,...x_k,v)$, $k \geq 1$ that start at $u$ and
end at $v$, and associate with each such walk the value $p(\omega)
= 1/(d_G(u)d_G(x_1)...d_G(x_k))$ (note, the $d(x_i)$ is not
ambiguous, since $x_i \notin E(H)$, but we leave the `$G$' subscript
in for clarity). This is precisely the probability that the walk
$\omega$ is taken by a simple random walk on $G$ starting at $u$.
Let
\begin{equation}
p_H(\vec{u,v})= \sum_{k \geq
1}\sum_{\omega=(u,x_1...x_k,v)}p(\omega),
\end{equation}
where the sum is over all exterior walks $\omega$.

We set the edge conductances (weights) of $H$ as follows: If $e$ is an interior
edge, $c(e)=1$. If it is an exterior edge $e=(u,v)$ define $c(e)$ as
\begin{equation}
c(e)= d_G(u)p_H(\vec{u,v}) = \sum_{k \geq
1}\sum_{\omega=(u,x_1...x_k,v)}\frac{1}{d_G(x_1)...d_G(x_k)} =
d_G(v)p_H(\vec{v,u})
\end{equation}

Thus the edge weight is consistent. A weighted
random walk on $H$ is thus a finite reversible Markov chain with all
the associated properties that this entails.
\\

\begin{definition}
\label{localObs} The weighted graph $H$ derived from $(G,S)$ is
termed \emph{the local observation of $G$ at $S$}, or \emph{$G$
locally observed at $S$}. We shall denote it as $H=Loc(G,S)$.
\end{definition}

The intuition in the above is that we wish to observe a random walk
$\mathcal{W}(G)$ on a subset $S$ of the vertices. When
$\mathcal{W}(G)$ makes an external transition at the border $B$, we
cease observing and resume observing if/when it returns to the
border. It will thus appear to have transitioned a virtual edge between the
vertex it left off and the one it returned on. It will therefore
appear to be a weighted random walk on $H$. This equivalence is
formalised thus
\\

\begin{definition}
\label{localObsWalk}
Let $G$ be a graph and $S \subset V(G)$. For an (unweighted) random walk $\mathcal{W}(G)$ 
on $G$ starting at $x_0 \in S$, derive the Markov chain $\mathcal{M}(G,S)$ on the states 
of $S$ as follows: (i) $\mathcal{M}(G,S)$ starts on $x_0$ (ii) If $\mathcal{W}(G)$
makes a transition through an internal edge $(u,v)$ then so does $\mathcal{M}(G,S)$ 
(iii)If $\mathcal{W}(G)$ takes an exterior walk $\omega=(u,x_1...x_k,v)$ then $\mathcal{M}(G,S)$
remains at $u$ until the walk is complete and subsequently transitions to $v$.
We call $\mathcal{M}(G,S)$ \emph{the local observation of $\mathcal{W}(G)$ at $S$}, or 
\emph{$\mathcal{W}(G)$ locally observed at $S$}.
\end{definition}

\begin{lemma}
\label{equiv}
For a walk $\mathcal{W}(G)$ and a set $S \subset V(G)$, the local observation of $\mathcal{W}(G)$ at $S$, $\mathcal{M}(G,S)$
is equivalent to the weighted random walk $\mathcal{W}(H)$ where $H=Loc(G,S)$.  
\end{lemma}
\proofstart
The states are clearly the same so it remains to show that the transition probability $P_{\mathcal{M}}(u,v)$ from $u$ to $v$ in $\mathcal{M}(G,S)$ is the same as $P_{\mathcal{W}(H)}(u,v)$ in $\mathcal{W}(H)$.
Recall that $B$ is the border of the induced subgraph $G[S]$. If $u \notin B$ then an edge $(u,v)\in E(H)$ is internal and so has unit conductance in $H$, as it does in $G$. Furthermore, for an internal edge $e$, $e\in E(H)$ if and only if $e\in E(G)$, thus $d_H(u) = d_G(u)$ when $u \notin B$. Therefore $P_{\mathcal{W}(H)}(u,v) = 1/d_H(u) = 1/d_G(u) = P_{\mathcal{M}}(u,v)$.

Now suppose $u \in B$. Let $E(u)$ denote the set of all edges incident with $u$ in $H$ and recall $A(B)$ 
above is the set of exterior edges. 
The total conductance (weight) of the exterior edges at $u$ is
\begin{eqnarray*}
\sum_{e \in E(u) \cap A(B)} c_H(e)&=&\sum_{x \in N(u, out)}\sum_{v \in B}
\Pr(\text{walk from } x \text{ returns to } B \text{ at } v)\\
&=&\sum_{x \in N(u, out)} 1\\
&=& d(u,out).
\end{eqnarray*}
(Note the $H$ subscript in $c_H(e)$ above is redundant since exterior edges are only defined for $H$, 
but we leave it for clarity).

Thus for
$u \in B$
\begin{eqnarray*}
c_H(u)=\sum_{e \in E(u)}c_H(e)&=& \sum_{e \in E(u) \cap G[S]} 1 +
\sum_{e \in E(u) \cap A(B)} c_H(e)\\
&=&d(u,in)+d(u,out)\\
&=&d_G(u)
\end{eqnarray*}

Now 
\begin{eqnarray}
P_{\mathcal{M}}(u,v) = \mathbf{1}_{\{(u,v) \in G[S]\}}\frac{1}{d_G(u)} 
+\sum_{k \geq 1}\sum_{\omega=(u,x_1...x_k,v)}\frac{1}{d_G(u)d_G(x_1)...d_G(x_k)}  
\end{eqnarray}
where the sum is over all exterior walks $\omega$.
Thus 
\begin{equation}
P_{\mathcal{M}}(u,v) = \mathbf{1}_{\{(u,v) \in G[S]\}}\frac{1}{d_G(u)} + p_H(\vec{u,v})
\end{equation}

\begin{eqnarray}
P_{\mathcal{W}(H)}(u,v) &=& \frac{1}{c_H(u)}\left[\mathbf{1}_{\{(u,v) \in G[S]\}} +  \mathbf{1}_{\{(u,v) \in A(B)\}}c_H(u,v)\right]\\
  			 &=& \frac{1}{d_G(u)}\left[\mathbf{1}_{\{(u,v) \in G[S]\}} +  \mathbf{1}_{\{(u,v) \in A(B)\}}d_G(u)p_H(\vec{u,v})\right]\\
  			 &=& \mathbf{1}_{\{(u,v) \in G[S]\}}\frac{1}{d_G(u)} +  \mathbf{1}_{\{(u,v) \in A(B)\}}p_H(\vec{u,v})\\
  			 &=& P_{\mathcal{M}}(u,v)
\end{eqnarray}
\proofend

\section{Effective Resistance Lemmas}\label{Effective Resistance Lemmas}

For the upper bound of Theorem \ref{main}, we require the following lemmas. 
\\

\begin{lemma}\label{subgraphEffResLemma}
Let $G$ be an undirected graph.  Let $G'\subseteq G$, be any subgraph such that such that $V(G') = V(G)$. For any $u,v \in V(G)$, 
\[
R(u,v)\leq R'(u,v)
\]
where $R(u,v)$ is the effective resistance between $u$ and $v$ in $G$ and $R'(u,v)$ similarly in $G'$. 
\end{lemma}
\proofstart
Since $V(G') = V(G)$, $G'$ can be obtained from $G$ by only removing edges. The lemma follows by the Cutting Law (Lemma \ref{CuttingLaw}).
\proofend

Denote by $R_{max}(G)$ the maximum effective
resistance between any pair of vertices in a graph $G$. 
\\

\begin{lemma}
\label{lg1x} 
For a graph $G$ and tree $T$, $R_{max}(G\Box T) < 4R_{max}(G\Box P_r)$ 
where $|V(T)| \leq r \leq 2|V(T)|$ and $P_r$ is the path on $r$ vertices.
\end{lemma}
\proofstart

Note first the following: 
\begin{description}
\item[(i)] By the parallel law, an edge $(a,b)$
of unit resistance can be replaced with two parallel edges between
$a,b$, each of resistance $2$. 
\item[(ii)] By the shorting law, a vertex
$a$ can be replaced with two vertices $a_1,a_2$ with a
zero-resistance edge between them and the ends of edges incident on
$a$ distributed arbitrarily between $a_1$ and $a_2$. 
\item[(iii)] By the same principle of the cutting law, this edge can be broken without decreasing effective resistance between any pair of vertices.
\end{description}
Transformations \textbf{(i)} and \textbf{(ii)} do not alter the effective resistance $R(u,v)$ between a pair of vertices $u,v$ in the network. For any vertex $u\notin\{a_1, a_2, a\}$, $R(u,a_1)=R(u,a_2)$ and these are equal to $R(u,a)$ before the operation. 

Points \textbf{(ii)} and \textbf{(iii)} require elaboration. In this thesis, we do not define zero-resistance (infinite conductance) edges. As stated in section \ref{RayleighStuff}, to say that a zero-resistance edge is placed between $a_1$ and $a_2$, is another way of referring to shorting as defined in Lemma \ref{ShortingLaw}. It would seem then, that \textbf{(ii)}, in fact says nothing. However, it serves as a useful short hand for talking about operations on the graph when used in conjunction with \textbf{(iii)}. If \textbf{(ii)} and \textbf{(iii)} are always used together, that is, if a zero-resistance edge created from \textbf{(ii)} is always cut by \textbf{(iii)}, then this is equivalent to the reverse of process of shorting two vertices $a_1$ and $a_2$ into $a_3$, as per Lemma \ref{ShortingLaw}. Hence, these two operations together are sound.

We continue thus:
\begin{enumerate}
\item Let $F=G\Box T$. Let each edge of $F$ have unit resistance. In what follows, we shall modify $F$, but shall continue to refer to the modified graphs as $F$. 

\item Starting from some vertex $v$ in $T$, perform a depth-first search (DFS) of $T$ stopping at the first return to $v$ after all vertices in
$T$ have been visited. Each edge
of $T$ is traversed twice; once in each orientation. Each vertex $x$ will be 
visited $d(x)$ times. 

\item Let $(e_i)$ be
the sequence of oriented edges generated by the search. The idea is to use
$(e_i)$ to construct a transformation from $F=G\Box T$ to $G\Box P_r$. 
From $(e_i)$, we derive
another sequence $(a_i)$, which is generated by following $(e_i)$ and if we have edges $e_i, e_{i+1}$
with $e_i = (a,b)$, $e_{i+1}=(b,c)$ such that it is neither the first time nor the last time $b$ is visited
in the DFS, then we replace $e_i,e_{i+1}$
with $(a,c)$. We term such an operation an \emph{aggregation}. Observe that in the sequence $(a_i)$, all leaf vertices of $T$ appear only once (just as in $(e_i)$), and a non-leaf vertex appears twice.

\item By \textbf{(i)} above, we can replace each (unit resistance) edge in $F$ by a pair of parallel edges each of 
resistance $2$. 

\item For a pair of parallel edges in the $T$ dimension, arbitrarily label one of them with an 
orientation, and label the other with the opposite orientation. Note, orientations are only an aid to the proof, and are not a flow restriction.
We therefore see that $(e_i)$ can be interpreted
as a sequence of these parallel oriented edges. 

\item We further modify $F$ using $(a_i)$: If $(a,b),(b,c)$ was aggregated to $(a,c)$,
then replace each pair of \textbf{oriented} edges $((x,a),(x,b))$ and $((x,b),(x,c))$ in $F$ with an oriented edge $((x,a),(x,c))$.
The resistances of $((x,a),(x,b))$ and $((x,b),(x,c))$ were $r((x,a),(x,b))=2$ and $r((x,b),(x,c))=2$. Set the resistance $r((x,a),(x,c))=r((x,a),(x,b))+r((x,b),(x,c)) = 4$. 

\item The above operation is the same as restricting flow through $((x,a),(x,b))$ and $((x,b),(x,c))$ to only going from one to the other at vertex $(x,b)$, without the possibility of going through other edges. The infimum of the energies of this subset of flows is at least the infimum of the energies of the previous set and so by Thomson's principle, the effective resistance cannot be decreased by this operation.

\item For each copy $G_i$ of $G$ in $F$ excluding those that correspond to a leaf of $T$, we can create a ``twin'' copy $G_i'$. 
Associate with each vertex $x \in V(F)$ (except those excluded) a newly-created twin vertex $x'$ with no incident edges. Thus, $V(G_i)$ has a twin set $V(G_i')$, though the latter has no edges yet. 

\item Recall the parallel edges created initially from all the edges of $F$; we did not manipulate those in the $G$ dimension, but we do so now:
redistribute half of the parallel edges of $G_i$ in the $G$ dimension to the set of twin vertices $V(G_i')$ so as to make $G_i'$ a copy of $G$ (isomorphic to it). Now put a zero-resistance edge between $x$ and $x'$. By \textbf{(ii)}, effective resistance is unchanged by this operation.
 
\item We now redistribute the oriented parallel edges in the $T$ dimension so as to respect the sequence $(a_i)$. We do this as follows: follow the sequence $(a_i)$ by traversing edges in their orientation. Consider the following event: In the sequence $(a_i)$ there is an element $a_j=(a,b)$  and $b$ has appeared in some element $a_i$ such that $i<j$. Then $a_j$ is the second time that $b$ has occurred in the sequence. Now change each edges $((x,a),(x,b))\in F$ to $((x,a),(x,b)')$. If $b=v$, then stop; otherwise, $a_j$ is followed by $a_{j+1}=(b,c)$, for some $c\in V(T)$. In this case, also change all $((x,b),(x,c))\in F$ to $((x,b)',(x,c))$. Continue in the same manner to the end of the sequence $(a_i)$.

\item We then remove the zero-resistance edges between each pair of twin vertices, and by \textbf{(iii)}, this cannot decrease the effective 
resistance. 
\end{enumerate}

Using the sequence $(a_i)$ to trace a path of copies of 
$G$, we see that the resulting structure is isomorphic to $G \Box P_r$. Since the aggregation process only aggregates edges that pass through a previously seen vertex, $r$ is at least $|V(T)|$. Also, because each edge is traversed at most once in each direction, $r$ is at most $2|V(T)|$. Each edge has resistance at most $4$, and so the lemma follows. 
\proofend 
\\

\begin{lemma}
\label{lg2x} 
For graphs $G,H$ suppose $D_G+1 \leq n_H \leq \alpha (D_G+1)$, for some $\alpha$. Then $R_{max}(G \Box H) < \zeta \alpha \log(D_G+1)$, where $\zeta$ is some universal constant.
\end{lemma}

\proofstart
Let $(a,x), (b,y)$ be any two vertices in $G \Box H$. Let $D$ be some diametric path of $G$. Let $\langle a,D\rangle$ 
represent the shortest path from $a$ to $D$ in $G$ (which may trivially be $a$ if it is on $D$). Similarly with
$\langle b,D\rangle$. Let $T_D = D \cup \langle a,D\rangle \cup \langle b,D\rangle$. Let $k=D_G+1$. Note $k \leq |V(T_D)| \leq 3k$. Now let $T_H$ be any
spanning tree of $H$. Applying Lemma \ref{lg1x} twice we have
\begin{equation*} 
R_{max}(T_D \Box T_H) < 4 R_{max}(T_D \Box P_s) < 16R_{max}(P_r\Box P_s)
\end{equation*}
where $k\leq r \leq 6k$ and $k\leq s \leq 2\alpha k$. Considering a series of connected $P^2_k$ subgraphs and using Lemma \ref{squaregrideffreslemma} and the triangle inequality for effective resistance, we have $R_{max}(P_r\Box P_s) \leq 16(6+2\alpha)8h(k)$, where $h(k)$ is the $k$'th harmonic number. Since $T_D \Box T_H \subseteq G \Box H$, the lemma follows by Lemma \ref{subgraphEffResLemma}.
\proofend

A diametric path $D$ is involved in the proof of Lemma \ref{lg2x} because the use of $D$ means that the dimension of $P_r$ is effectively maximised, and we can break up the grid $P_r\Box P_s$ roughly into $k\times k$ square grids, each with maximum effective resistance $O(\log k) = O(\log D_G)$. If, for example, the shortest path between $a$ and $b$ is used, the product $P_r\Box P_s$ may have $r$ much smaller than $s$, looking like a long thin grid, which may have a high effective resistance. 

\ignore{
Lemma \ref{lg2x} gives us an upper bound of $\zeta\log (D_G+1)$ for the effective resistance in a block (definition below),
which in turn allows us to bound the maximum hitting time within a block, and therefore the cover time of the block via Matthews'
technique.
}

\section{A General bound}
In this section, we prove Theorem \ref{main}, starting with the lower bound.

\subsection{Lower Bound}
The following is a partial restatement of Theorem \ref{main} for the lower bound (inequality \eqref{LBIneq}).
\\

\textbf{Theorem \ref{main} (partial restatement)} \emph{
Let $F= (V_F, E_F) = G \Box H$ where $G=(V_G, E_G)$ and $H = (V_H, E_H)$ are simple, connected, unweighted, undirected graphs. 
We have
\begin{equation}
\cov[F] \geq \max\left\{\left(1+\frac{\delta_G}{\Delta_H}\right) \cov[H], \left(1+\frac{\delta_H}{\Delta_G}\right) \cov[G]\right\}.
\end{equation}
}

\proofstart
In order for the walk $\mathcal{W}$ to cover $F$, it needs to have covered the $H$ dimension of $F$. That is, each copy of $G$ in $F$ needs to have been visited at least once. The probability of a transition in the $H$ dimension is distributed as a geometric random variable with success probability at most  $\frac{\Delta_H}{\Delta_H+\delta_G}$. Thus, the expectation of the number of steps of $\mathcal{W}$ per transition in the $H$ dimension is at least $\frac{\Delta_H+\delta_G}{\Delta_H}$. Transitions of $\mathcal{W}$ in the $H$ dimension are independent of the location of $\mathcal{W}$ in the $G$ dimension, and have the same distribution (in the $H$ dimension) as a walk on $H$. This proves 
\[
\cov[F] \geq \left(1+\frac{\delta_G}{\Delta_H}\right) \cov[H].
\]
By commutativity, 
\[
\cov[F] \geq \left(1+\frac{\delta_H}{\Delta_G}\right) \cov[G].
\]

\proofend

\subsection{Upper Bound}
The following proves the upper bound in Theorem \ref{main}. It is envisaged that theorem is used with the idea in mind that $G$ is small relative to $H$, and so the cover time
of the product is essentially dominated by the cover time of $H$. We give a partial restatement of the theorem for the upper bound (inequality \eqref{UBIneq}).
\\

\textbf{Theorem \ref{main} (partial restatement)} \emph{Let $F= (V_F, E_F) = G \Box H$ where $G=(V_G, E_G)$ and $H = (V_H, E_H)$ are simple, connected, unweighted and undirected. Suppose further that $n_H \geq D_G+1$.
\begin{equation}
\cov[F] \leq K\left(\left(1+\frac{\Delta_G}{\delta_H}\right)\bcov[H] + \frac{Mm_Gm_Hn_H\ell^2}{\cov[H]D_G}\right)
\end{equation}
where $\ell=\log (D_G+1) \log(n_GD_G)$ and $K$ is some universal constant.
}

\proofstart
Let $k = D_G+1$. We group the vertices of $H$ into sets such that for any set $S$ and the subgraph of $H$ induced by $S$, $H[S]$: 
(i)$|S| \geq k$, (ii)$H[S]$ is connected, (iii) The diameter of $H[S]$ is at most $4k$. We do this through the following decomposition algorithm on $H$: Choose some arbitrary vertex $v\in V(H)$ as the root, and using a breadth-first search (BFS) on $H$, 
descend from $v$ at most distance $k$. The resulting tree $T(v)\subseteq H$ will have diameter at most $2k$. For each leaf $l$ of $T(v)$,
continue the BFS using $l$ as a root. If $T(l)$ has fewer than $k$ vertices, append it to $T(v)$. If not, recurse on the leaves of $T(l)$.
The set of vertices of each tree thus formed satisfies the three conditions above. The root is part of a new set, unless it has been appended to another tree.

In the product $F$ we refer to copies of $G$
as \emph{columns}. In $F$ we have 
a natural association of each column with the set $S \subseteq V(H)$ defined above. We define $Block[S] =(G\Box H[S])$.

[Refer to section \ref{Notation} for a reminder of the notation $(.,y)$]. For any two vertices $(.,a), (.,b) \in G\Box H[S]$ there exists a tree $T\langle a,b\rangle$ subgraph of the tree $T$ in $H$ that generated $S$ such that $a$ and $b$ are connected in $T\langle a,b\rangle$
and $k \leq |V(T\langle a,b\rangle)| \leq 4k$. Then using Lemmas \ref{lg2x} and \ref{subgraphEffResLemma}, we can upper bound the effective resistance 
$R((.,a), (.,b))$ in $B=Block[S]$,  

\begin{equation}
\label{R_{max}(B)}
R_{max}(B) \leq 4\zeta\log (D_G+1).
\end{equation}

Furthermore, if $B'=Loc(F,V(B))$ ($Loc$ is defined in Definition \ref{localObs}), then $B \subseteq B'$ so by Lemma \ref{subgraphEffResLemma},

\begin{equation}
\label{R_{max}(B')}
R_{max}(B') \leq 4\zeta\log (D_G+1).
\end{equation}

We  use the following two-phase approach to bound the cover time of $F= G\Box H$.
\begin{description}
  \item[Phase 1] Perform a random walk $\mathcal{W}(F)$ on $F$ until the blanket-cover criterion is satisfied for the $H$ dimension.
  \item[Phase 2] Starting from the end of phase 1, perform a random walk on $F$ until all vertices of $F$ not visited in phase 1 are visited.
\end{description}
Phase 1 can be thought of in the following way: We couple $\mathcal{W}(F)$ with a walk $\mathcal{W}(H)$ such that (i) if $\mathcal{W}(F)$ starts at $(.,x)$, then $\mathcal{W}(H)$ starts at $x$, and (ii) $\mathcal{W}(H)$ moves to a new vertex $y$ from a vertex $x$ when and only when $\mathcal{W}(F)$ moves from $(.,x)$ to $(.,y)$. This coupled process runs until $\mathcal{W}(H)$ satisfies the blanket-cover criteria for $H$, i.e., when each vertex $v \in V(H)$ has been visited at least $\pi(v)\cov[H]$ times. An implication is that the corresponding column $G_v$ in $F$ will have been visited at least that many times.

Having grouped $F$ into blocks, we analyse the outcome of phase 1 by relating $\mathcal{W}(F)$ to the local observation on each block. A particular
block $B$ will have some vertices unvisited by $\mathcal{W}(F)$ if and only if $\mathcal{W}(F)$ locally observed on $B$ fails to visit all vertices.
We refer to such a block as \emph{failed}. Consider the weighted random walk $\mathcal{W}(B')$ on $B' = Loc(F,V(B))$. This has the same distribution as $\mathcal{W}(F)$ locally observed on $B$. Hence, we bound the probability of $\mathcal{W}(F)$ failing to cover $B$ by bounding the probability that
$\mathcal{W}(B')$ fails to cover $B'$. Done for all blocks, we can bound the expected time it takes phase 2 to cover the failed blocks. We think of phase 1 as doing most of the ``work'', and phase 2 as a ``mopping up'' phase. Mopping up a block in phase 2 is costly, but if there are few of them, the overall cost is within a small factor of phase 1. 

We bound $\Pr(\mathcal{W}(B') \text{ fails})$ by exploiting the fact that $\mathcal{W}(B')$ will have made some minimal number of transitions $t$. This is guaranteed because phase 1 terminates only when $\mathcal{W}(H)$ has satisfied the blanket-cover criterion on $H$. If $\kappa$ counts the number of steps of a walk $\mathcal{W}(B')$ until $B'$ is covered, then
\begin{equation}
\Pr(\mathcal{W}(B') \text{ fails to cover} B') \leq \Pr(\kappa > t) \leq \frac{\E[\kappa]}{t} \label{probfailMArkov}
\end{equation}
by Markov's inequality.
\\

\begin{definition}
For graphs $I = J \Box K$, and $S \subseteq V(I)$, denote by $S.K$ the
\emph{projection of $S$ on to $K$}, that is, $S.K = \{v \in K: (.,v)
\in  S\}$.
\end{definition}

For a weighted graph $G$, recall that $c(G)$ is the twice the sum of the conductances (weights) of all
edges of $G$ (refer to section \ref{weighted graphs}). By the definition of $G\Box H[S]$ and section \ref{DefnSubgraphWalk},
\begin{eqnarray}
    c(B') \leq m_G|V(B).H| +n_G\sum_{u \in V(B).H}d(u) \label{c(B')UpperBound}
\end{eqnarray}

Using \eqref{R_{max}(B')} and
Theorem \ref{Th:commtetimetheorem} we
 therefore have for any $u,v \in V(B')$, 
$\com[u,v] \leq Kc(B')\log (D_G+1)$ for some universal constant $K$. (In what follows $K$ will change, but we shall keep the same
symbol, with an understanding that what we finish with is a universal constant).
Hence, by Matthews' Technique (Theorem \ref{MatLemma1}), 
\begin{equation}
\cov[B']\leq Kc(B')\log (D_G+1) \log (|V(B')|).
\end{equation}

For a block $B$, the number
of transitions on the $H$ dimension - and therefore the number of
transitions on $B$ -  as demanded by the blanket-cover criterion is at least

\begin{eqnarray}
    \tau = \sum_{u \in V(B).H}\pi_H(u)\cov[H] =  \frac{\cov[H]}{2m_H}\sum_{u \in V(B).H}d_H(u), \label{tTime}
\end{eqnarray}
where $\pi_H(u)$ and $d_H(u)$ denote the stationary probability and degree of $u$ in $H$.

Now

\begin{eqnarray}
    \Pr(\mathcal{W}(F) \text{ fails on }B)
     &=& \Pr(\mathcal{W}(B') \text{ fails on }B') \nonumber \\
     &\leq& Kc(B')\log (D_G+1) \log (|V(B')|)/\tau, \label{xxx}
\end{eqnarray}
as per \eqref{probfailMArkov}. For convenience, we left $l_B=\log (D_G+1) \log (|V(B)|)$ (recall $V(B)=V(B')$). Hence, using \eqref{tTime} with \eqref{c(B')UpperBound} and \eqref{xxx},
\begin{eqnarray}
    \Pr(\mathcal{W}(F) \text{ fails on }B) &\leq& \frac{Kl_Bm_H}{\cov[H]}\frac{m_G|V(B).H| +n_G\sum_{u \in V(B).H}d(u)}{\sum_{u \in V(B).H}d_H(u)}\nonumber \\
    &=&  \frac{Kl_Bm_H}{\cov[H]}\left(n_G+\frac{m_G|V(B).H|}{\sum_{u \in V(B).H}d_H(u)}\right)\label{xxxy}
\end{eqnarray}

Phase 2 consists of movement between failed blocks,
and covering a failed block it has arrived at. The total
block-to-block movement is upper bounded by the time is takes to
cover the $H$ dimension of $F$ (in other words, for each column to
have been visited at least once). We denote this by $\cov_F[H]$.
Let $\cov_F[B]$ denote the cover time of the set of vertices of a block $B$ by the walk $\mathcal{W}(F)$. Let the random variables $\phi_1$ and  $\phi_2$ represent the time it takes to complete phase 1 and phase 2 respectively. 

\begin{eqnarray*}
    \E [\phi_2] \leq \cov_F[H] +\sum_{B \in F} \Pr(\mathcal{W}(F) \text{ fails on }B)\cov_F[B].
\end{eqnarray*}

For $\mathcal{W}(H)$, the random variable $\beta_H =
\min\{t:(\forall v)N_v(t)\geq\pi(v)\cov[H]\}$ counts the time it takes
to satisfy the blanket-cover criterion on $H$. 

The expected number of movements on $F$ per movement on the $H$ dimension is at most $(\Delta_G + \delta _H)/\delta_H$.
Therefore, 
\[
\E [\phi_1] \leq \frac{\Delta_G + \delta _H}{\delta_H} \E[\beta_H] = \frac{\Delta_G + \delta _H}{\delta_H}\bcov[H].
\]
Similarly, 
\[
\cov_F[H] \leq \frac{\Delta_G + \delta _H}{\delta_H}\cov[H].
\]

Using \eqref{R_{max}(B)}, Lemma \ref{subgraphEffResLemma}, and Theorems \ref{Th:commtetimetheorem} and \ref{subsMatTheorem} on $B$, we have 
\begin{equation}
\cov_F[B]\leq K'c(F)l_B \label{cov_F[B]UB}
\end{equation}
where $c(F)= 2|E(F)| = 2M$. 

Hence,
\begin{eqnarray*}
\cov[F] &\leq& \E [\phi_1] + \E [\phi_2]\\ 
				&\leq& K\frac{\Delta_G + \delta _H}{\delta_H}\bcov[H] +\sum_{B \in F} \Pr(\mathcal{W}(F) \text{ fails on }B)\cov_F[B]. 
\end{eqnarray*}

We have, using \eqref{xxxy} and \eqref{cov_F[B]UB},
\begin{equation}
\sum_{B \in F} \Pr(\mathcal{W}(F) \text{ fails on }B)\cov_F[B]
 \leq K\frac{Mm_H}{\cov[H]}\sum_{B \in F}\left(n_G+\frac{m_G|V(B).H|}{\sum_{u \in V(B).H}d_H(u)}\right)l_B^2. \label{fjsndjsn}
\end{equation}

Since $\sum_{u \in V(B).H}d(u) \geq |V(B).H|$, the outer summation in \eqref{fjsndjsn} can be bounded thus
\begin{eqnarray}
\sum_{B \in F}\left(n_G+\frac{m_G|V(B).H|}{\sum_{u \in V(B).H}d_H(u)}\right)l_B^2 \leq m_G\log^2(D_G+1)\sum_{B \in F}\log^2(|V(B)|). \label{oisjfoijo}
\end{eqnarray}
Since each block $B \in F$ has at least $D_G+1$ columns, we can upper bound the sum in the RHS of \eqref{oisjfoijo} by assuming all blocks have this minimum. The number of such blocks in $F$ will be $|V(H)|/(D_G+1)$, each block having $(D_G+1)n_G$ vertices. Hence

\begin{equation}
\sum_{B \in F}\log(|V(B)|)^2 \leq \frac{n_H}{D_G}\log^2(n_G(D_G+1)). \label{ifjsohgifkjj} 
\end{equation}

Putting together \eqref{fjsndjsn}, \eqref{oisjfoijo} and \eqref{ifjsohgifkjj}, we get

\begin{equation*}
\sum_{B \in F} \Pr(\mathcal{W}(F) \text{ fails on }B)\cov_F[B] \leq K\frac{Mm_Gm_Hn_H\ell^2}{\cov[H]D_G}
\end{equation*}

where $\ell = \log (D_G+1)\log (n_GD_G)$.

\proofend

\chapter{Random Graphs of a Given Degree Sequence}\label{CovDS}

In this chapter we study the asymptotic cover time of random graphs that have a prescribed 
{\em degree sequence}, that is, for a graph with vertex set $V = \{1,2,\ldots,n\}$ with $n \to \infty$, we have a sequence $\textbf{d} = (d_1, d_2,\ldots, d_n)$ of positive integers where $d_i$ is the degree of vertex $i$ and we wish to determine the cover time of a graph picked uniformly at random from the set of all connected simple graphs of degree sequence $\textbf{d}$. Denote by \(\mathcal{G}(\textbf{d})\) the space of all such graphs with a uniform distribution on them. Thus, we study the cover time of a graph $G$ picked from \(\mathcal{G}(\textbf{d})\). We may relax our terminology slightly and say that a graph $G$ is picked \uar\ from the set \(\mathcal{G}(\textbf{d})\), even though \(\mathcal{G}(\textbf{d})\) is not a set but a set with an associated distribution on the elements.

It must be noted that not every sequence of $n$ positive integers will allow for a simple graph, nor even a graph. We know, for example, for the sequence to be graphical, the sum of the degrees has to be even, i.e., $\sum_{i=1}^{n}d_i = 2m$ for some natural number $m$. With this condition, a graph having the prescribed sequence may exist, but may not be simple. For example, a $2$-vertex graph with degree sequence $(1,3)$ will force multiple edges and/or loops. Nevertheless, the study of the types of graphs can be facilitated by means of an intimately related random process known as the \emph{configuration model}, which is explained in section \ref{TheConfigurationModelSection}, and is treated in many sources, including, e.g., \cite{Janson}.

We reiterate some definitions: a statement $\mathcal{P}(n)$ parameterised on an integer $n$ holds \emph{with high probability} (\whp) if $\Pr(\mathcal{P}(n) \text{ is true})\rightarrow 1$ as $n \rightarrow \infty$; the notation $f(n) \sim g(n)$ means $f(n)/g(n) \rightarrow 1$ as $n \rightarrow \infty$.

Our study of the cover time imposes certain technical restrictions on the degree sequence under consideration in addition to the $\sum_{i=1}^{n}d_i = 2m$ restriction already mentioned. A degree sequence \textbf{d} satisfying these restrictions, which are described in section \ref{seca}, is called \emph{nice}.

We denote by $\theta$ the average vertex degree, i.e., $\theta=2m/n$,  and by $d$  the \emph{effective minimum degree}. The latter is a fixed positive integer, and the first entry in the ordered degree sequence which occurs $\Theta(n)$ times, meaning that for any $d'<d$, there are o(n) vertices with degree $d'$. The significance of the effective minimum degree is discussed in \ref{seca}.
\\

\begin{theorem} \label{CovTh}
Let  \(G\) be chosen \uar\ from \(\mathcal{G}(\textbf{d})\), where \textbf{d} is nice.
Then \whp,
\begin{equation}
\cov[G] \sim \frac{d-1}{d-2}\frac{\theta}{d}\;n\log n. \label{Covertime}
\end{equation}
\end{theorem}

The logarithms are base-$e$, as they are in the rest of this chapter, unless stated otherwise.

We note that if \(d\sim\th\), i.e., the graph is pseudo-regular, then
\[
\cov[G] \sim\frac{d-1}{d-2}\;n\log n,
\]
which is the same asymptotic limit for random \(d\)-regular graphs given in \cite{CFreg}.

\section{Random Graphs: Models and Cover Time}\label{Random Graphs: Models and Cover Time}
The study of random graphs goes back to \cite{ER} and \cite{Gilbert} and has since become a highly active field of research. Those two models have proved to be a highly fertile ground in which to develop ideas, tools and techniques for the study of random graphs. Despite the fact that these two models were introduced by different groups of authors, both of them have come to be referred to as \emph{Erd\H{o}s--R\'{e}nyi} (E--R) random graphs, after the authors of \cite{ER}.

Subsequent models attempted to deal with the shortcomings of {Erd\H{o}s-Renyi in accurately capturing the structural properties of ``real world'' networks. This has been an area where mathematical theory and empirical study have been mutually beneficial to each other; theory has served to deepen understanding of ``real-world'' networks, and empirical research has generated data that has led to the developments of theory. Section \ref{Other Models} gives examples.

The classical work on random graphs is \cite{Bellabook}. Another popular example is \cite{Janson}, but there are many others.

In the next section, we describe some models of random graphs as well as cover time results on them. We start with the two original ones of \cite{ER} and \cite{Gilbert}. It should be noted, that in all graphs we consider, vertices are labelled, and thus two graphs which may be indistinguishable without labelling, will be different objects with labelling.

Before doing so, we define the term graph space.  A \emph{graph space} $\mathcal{G}$ is a set of graphs together with an assignment of probability $p(G)$ to each graph $G\in \mathcal{G}$ such that $\sum_{G\in \mathcal{G}}p(G)=1$, i.e.\, it is a probability distribution. When we say we pick a random graph $G\in \mathcal{G}$, we mean we are picking it from the set with probability $p(G)$. 

\subsection{Erd\H{o}s--R\'{e}nyi}
In the model of \cite{ER}, the graph space $\mathcal{G}(n,m)$ is the set of all graphs with $n$ (labelled) vertices and $m$ edges, together with a uniform distribution on the set. Hence, a graph $G$ is picked uniformly at random from all graphs on $n$ vertices and $m$ edges.  In the model of \cite{Gilbert}, there are $n$ (labelled) vertices and each of the $\binom{n}{2}$ possible edges exists with some fixed probability $p$, independently of the others. This graph space, denoted by $\mathcal{G}(n,p)$ contains every graph on $n$ vertices, but the distribution on them is not uniform. For a particular graph $G$ on $n$ vertices, if the edge set $E(G)$ is such that $|E(G)|=m$, then the probability of $G$ being picked from $\mathcal{G}(n,p)$ is 
\[
p^m(1-p)^{\binom{n}{2}-m}.
\]
Note, that whilst these graphs will always be simple, they may not be connected. Conditions for connectivity were studied with the introduction of the models, and continued to be thereafter, becoming a major area of focus for research on these models. 

A cover time result for $\mathcal{G}(n,p)$ was given in \cite{JonassonCover}:
\\

\begin{theorem}[\cite{JonassonCover}]
For $G\in \mathcal{G}(n,p)$, \whp,
\begin{description}
	\item [(i)] If $\frac{np}{\log n} \rightarrow \infty$ then $\cov[G] = (1+o(1))n\log n$.
	\item [(ii)] If $c>1$ is a constant and $np=c\log n$ then $\cov[G]>(1+\alpha)n\log n$ for some constant $\alpha= \alpha(c)$.
\end{description}
\end{theorem}
This result was then strengthened in \cite{CooperER}:
\\

\begin{theorem}[\cite{CooperER}]
Suppose that $np = c\log n = \log n + \omega$ where $\omega = (c-1)\log n \rightarrow \infty$ and $c=O(1)$. If $G\in \mathcal{G}(n,p)$, 
then \whp, 
\[
\cov[G] \sim c \log{\left(\frac{c}{c-1}\right)}n\log n.
\]
\end{theorem}

\subsection{Random Regular}
Let $r$ be a positive integer. A \emph{random $r$-regular} graph $G$ on $n$ vertices is a graph picked \uar\ from the set of of all $r$-regular graphs on $n$ vertices. The graph space is denoted by $\mathcal{G}(n,r)$.

The cover time for random regular graphs was studied in \cite{CFreg}:
\\

\begin{theorem}[\cite{CFreg}]
Let $r\geq 3$ be a constant. For $G\in \mathcal{G}(n,r)$, \whp,
\[
\cov[G] \sim \frac{r-1}{r-2}n\log n.
\]
\end{theorem}

\subsection{Other Models}\label{Other Models}
The cover time of a particular generative model of the \emph{preferential attachment graph} is studied in \cite{CGweb}. In this model, at each time step, a new vertex $v$ is added to the graph, and a fixed number $m$ edges are randomly added between $v$ and the existing vertices. The probability of adding to a vertex $u$ is in proportion to the degree of $u$ at that time in the process. The cover time was determined to be asymptotically equal to $\frac{2m}{m-1}n\log n$, where $n$ is the final number of vertices in the graph.

This model was suggested by \cite{Barabasi} as a means of generating graphs with the \emph{scale-free} property, which is the name given in the same paper to graphs having a power-law degree distribution. This is one in which $P(k)\propto k^{-\gamma}$ for some constant $\gamma$, where $P(k)$, the fraction of vertices with degree $k$. Such a property is thought to exist in many ``real world'' networks, such as the WWW and actor collaboration networks \cite{Barabasi}, and the Internet \cite{Faloutsos} \footnote{The model of \cite{Barabasi} was analytically determined by \cite{BollPref} to have $P(k)\propto k~^{-3}$ for all $k\leq n^{1/15}$ where $n$ is the final number of vertices in the model. This closely matched simulation results of \cite{Barabasi} and \cite{Barabasi2}
which gave values for $\gamma$ of $2.9 \pm 0.1$. As a comparison, experimental studies for the WWW \cite{Albert} suggests $\gamma\approx 2.1$ and $\gamma\approx 2.45$ for the in-degree and out-degree respectively. Similarly, experimental studies for the Internet in \cite{Faloutsos} suggest $\gamma$ between $2.15$ and $2.20$.}.

In a \emph{random geometric graph} (\cite{Penrose}), the $n$ vertices are scattered \uar\ on (some subset of) a $d$-dimensional space where $d\geq 2$, and an edge is placed between vertices $u$ and $v$ if the Euclidean distance between them is at most some fixed constant $r$, often called the \emph{radius}. This type of random graph for $d=2$ has been used as a model of wireless ad-hoc and sensor networks (\cite{Gupta}, \cite{chenAdhoc}, \cite{Broutin}) where the radius represents the radio communication range of devices that are placed randomly on the plane. Two recent papers deal with the cover time of random geometric graphs: \cite{ColinGeometric} and \cite{chenGeometric}.

\section{Mixing Time, Eigenvalues and Conductance}
In this section we discuss parameters that are related to random walks on graphs. The concepts introduced here play a fundamental role in the proof Theorem \ref{CovTh}. 

In chapter \ref{Theory of Markov Chains and Random Walks}, we discussed the convergence of a random walk $\mathcal{W}$ on a graph $G$ to a unique stationary distribution $\pi$. The formal statement was made in Theorem \ref{Th:RandomWalkStationDist}, which in turn was based on Theorem \ref{ergodicTheorem} for Markov chains more general than random walks on undirected graphs.

Although these theorems asserted that a walk would converge to a stationary distribution (given certain conditions), there was no mention of how \emph{quickly} this convergence happens. That is, there was no mention of how close to stationarity was the distribution of the random walk $\mathcal{W}$ after some number of steps $t$, nor how many steps were required to get close to stationarity (for some well-defined meaning of ``close''). There are a number of related definitions of ``closeness'' between one distribution and another, and also a number of related definitions of ``quickness'' for convergence of distributions, but they are similar enough to one another that they all reflect the fundamental behaviour of a walk on the graph in roughly the same way. The rate of convergence of the walk is called the \emph{mixing rate}, and the time it takes the walk to get close to stationarity is called the \emph{mixing time}. A walk is \emph{rapidly mixing} if the mixing time is somehow small compared to the size of the graph - say, polylogarithmic in the number of vertices. We shall define these notions in precise terms (see, e.g., \cite{LPW} or \cite{LovaszSurvey}), but first we shall introduce the role of eigenvalues in the study of random walks on graphs.

\subsection{Theory and Application of the Spectra of Random Walks}
Recall that for the matrix $\mathbf{P}$ of transition probabilities of a random walk on a graph $G$, the $t$-step probabilities are given by $\mathbf{P}^t$. This suggests that the set of eigenvalues of $\mathbf{P}$ - the \emph{spectrum} - and their associated eigenvectors may have some important role. Indeed, we have already seen one particularly important left-eigenvector, the stationary distribution $\pi$: $\pi \mathbf{P} = \pi$. Note further that   $\mathbf{P1}=\mathbf{1}$, where $\mathbf{1}$ is the column vector of $n$ elements with every element $1$. 

The spectral theory of the transition matrix allows us to prove a convergence in distribution, and as we shall see, also give bounds on how close the distribution at time $t$ is to the limit give some starting distribution $\mathbf{p}$. This is given in terms of a distance between $\mathbf{p}\mathbf{P}^t$ and $\pi$, bounded by a function of eigenvalues and $t$. 

We follow the presentation given in \cite{LovaszSurvey}. Recall the definition of the adjacency matrix given in section \ref{graphdefns}. Let $\mathbf{A}$ be the adjacency matrix and $\mathbf{P}$ be the transition matrix for a connected, simple, undirected and unweighted graph $G$. Let $n=|V(G)|$, and without loss of generality, let us assume the vertices of $G$ are labelled $1$ to $n$. Let $\mathbf{D}$ be the diagonal $n\times n$ matrix such that $\mathbf{D}_{i,i} = 1/d(i)$, the degree of vertex $i$. Observe $\mathbf{P} = \mathbf{DA}$.

Although $\mathbf{A}$ is symmetric, $\mathbf{P}$ will not be unless $G$ is a regular graph. In order to use the powerful tools of Spectral Theory, we require a matrix to be symmetric. We therefore use the related matrix $\mathbf{N}=\mathbf{D}^{1/2}\mathbf{A}\mathbf{D}^{1/2}= \mathbf{D}^{-1/2}\mathbf{P}\mathbf{D}^{1/2}$, which is, in fact, symmetric. 
\\

\begin{proposition}\label{SpectralProp}
For a real, symmetric $n\times n$ matrix $\mathbf{M}$, by the Spectral Theorem, 
\begin{enumerate}
	\item Eigenvectors of $\mathbf{M}$ with different eigenvalues are orthogonal; eigenvectors with the same eigenvalue need not be.
	\item $\mathbf{M}$ has a full orthonormal basis of eigenvectors $\mathbf{v_1}, \mathbf{v_2}, \ldots, \mathbf{v_n}$, with corresponding eigenvalues $\lambda_1, \lambda_2, \ldots, \lambda_n$. All eigenvalues and eigenvectors are real.
	\item $\mathbf{M}$ is diagonalisable:
				\[
					\mathbf{M}= \mathbf{E}\mathbf{\Lambda} \mathbf{E}^T
				\]
				where the columns of $\mathbf{E}$ are the orthonormal basis $\mathbf{v_1}, \mathbf{v_2}, \ldots, \mathbf{v_n}$, and $\mathbf{\Lambda}$ is a diagonal matrix with entries corresponding to the eigenvalues of the columns of $\mathbf{E}$ (in corresponding order). Thus, $\mathbf{M}$ can be expressed in the following form:
				\[
					\mathbf{M} = \sum_{i=1}^n\lambda_i\mathbf{v_i}\mathbf{v_i}^T.
				\]
\end{enumerate}  
\end{proposition}
See, for example, \cite{hornmatrix} for details.

Since $\mathbf{N}$ is real and symmetric, then by Proposition \ref{SpectralProp} it has the form
\[
	\mathbf{N} = \sum_{i=1}^n\lambda_i\mathbf{v}_i\mathbf{v}_i^T,
\] 
where $\lambda_1\geq \lambda_2\geq \ldots \lambda_n$ and the $v_i$ form an orthonormal set. Consider the column vector $\mathbf{w}$ with $\mathbf{w}_i=\sqrt{d(i)}$, where $d(i)$ is the degree of vertex $i$. Observe
\begin{eqnarray*}
\mathbf{Nw} &=& \mathbf{D}^{-1/2}\mathbf{P}\mathbf{D}^{1/2}\mathbf{w}\\ 
&=& \mathbf{D}^{-1/2}\mathbf{P1}\\ 
&=& \mathbf{D}^{-1/2}\mathbf{1}\\ 
&=& \mathbf{w}.
\end{eqnarray*}
Thus $\mathbf{w}$ is an eigenvector of $\mathbf{N}$ with size $\|w\| = \sqrt{\sum_{i=1}^nd(i)} = \sqrt{2m}$, where $m=|E(G)|$ is the number of edges in $G$. Since this eigenvector is positive, by the Perron--Frobenius theorem, the eigenvalue associated with it is unique and strictly larger than the second largest eigenvalue. Furthermore, it is at least the absolute size of the smallest eigenvalue (which may be negative). That is, for the eigenvalues  $\lambda_1, \lambda_2, \ldots, \lambda_n$ of $\mathbf{N}$,
\begin{equation}
1=\lambda_1> \lambda_2\geq \ldots \geq \lambda_n \geq -1 \qquad \text{ and }\qquad \lambda_1 \geq |\lambda_n|. \label{PerronFrob}
\end{equation}
It therefore follows that $\mathbf{v_1}=\mathbf{w}/\|\mathbf{w}\|$, i.e.,\  $\mathbf{v_1}[i] = \sqrt{d(i)/2m} =\sqrt{\pi_i}$.

For a proof of the following proposition, see, for example \cite{Lovaszcombbook}.
\\

\begin{proposition}
If $G$ is non-bipartite, then $\lambda_n>-1$
\end{proposition}

The following is from \cite{LovaszSurvey} (with modifications for consistency of notation).
\[
\mathbf{P}^t = \mathbf{D}^{1/2}\mathbf{N}^t\mathbf{D}^{-1/2} = \sum_{k=1}^n\lambda_k^t\mathbf{D}^{1/2}\mathbf{v_k}\mathbf{v_k}^T\mathbf{D}^{-1/2} = \mathbf{Q} + \sum_{k=2}^n\lambda_k^t\mathbf{D}^{1/2}\mathbf{v_k}\mathbf{v_k}^T\mathbf{D}^{-1/2}
\]
where $\mathbf{Q}_{i,j} = \pi_j$. That is,
\begin{equation}
p_{i,j}^{(t)}=\mathbf{P}_{i,j}^t = \pi_j + \sum_{k=2}^n\lambda_k^t\mathbf{v_k}[i]\mathbf{v_k}[j]\sqrt{\frac{d(j)}{d(i)}}. \label{specconveq}
\end{equation}
If $G$ is not bipartite, then $|\lambda_k|<1$ for $2\leq k \leq n$ and so
\[
p_{i,j}^{(t)} \rightarrow \pi_j \text{ as } t \rightarrow \infty. 
\]
This proves Theorem \ref{Th:RandomWalkStationDist}, and moreover, demonstrates how the eigenvalues are related to the speed of convergence of the walk to the stationary distribution. 

Thus, from equation (\ref{specconveq}), we can see that for any $i,j$  
\begin{equation}
|p_{i,j}^{(t)} -\pi_j| \leq  \sqrt{\frac{d(j)}{d(i)}}\lambda^t_*. \label{specdistbound}
\end{equation}
where $\lambda_* = \max\{\lambda_2, |\lambda_n|\}$,  since $\lambda_2 \geq \lambda_k$ for $2\leq k \leq n$. The quantity $\lambda_1 - \lambda_2 = 1-\lambda_2$ is called the  \emph{spectral gap} and the quantity $\lambda_1 - \lambda_* = 1-\lambda_*$ is called the \emph{absolute spectral gap}, and bounding these quantities is a common means of bounding mixing time.

\subsection{Conductance}\label{conductancesectionprelims}
The conductance we refer to in this section and for the rest of this chapter is \emph{not} the same concept as that used in conjunction with electrical network theory, where it refers to weight of an edge in a network. It will be suitable for our purposes to first state the definition in terms of Markov chains, and then reduce it for random walks on graphs. Definitions are given in e.g., \cite{Durrett}, \cite{LPW} and \cite{LovaszSurvey}.
\\

\begin{definition}[\textbf{Conductance}]\label{MarkovchainCondDef}
Let $\mathcal{M}$ be an irreducible, aperiodic Markov chain on some state space $\Omega$. Let the stationary distribution of $\mathcal{M}$ be  $\pi$ with $\pi(x)$ denoting the stationary probability of $x\in \Omega$. Let $\mathbf{P}$ be the transition matrix for $\mathcal{M}$. For $x,y\in \Omega$ let $Q(x,y)=\pi(x)\mathbf{P}[x,y]$ and for sets $A, B \subseteq \Omega$, let $Q(A,B) =\sum_{x \in A, y\in B}Q(x,y)$. The \emph{conductance} of $\mathcal{M}$ is the quantity
\begin{equation}
\Phi = \Phi(\mathcal{M})=\min_{\substack{S \subseteq \Omega\\ \pi(S)\leq 1/2}}\frac{Q(S,\overline S)}{\pi(S)} \label{MarkovCond}
\end{equation}
where $\pi(S) = \sum_{x\in S}\pi(x)$, and  $\overline S = \Omega \setminus S$.
\end{definition}

For an unweighted simple graph $G=(V,E)$ with $n=|V|$, $m=|E|$, 
\[
 Q(i,j) = \pi(i)\mathbf{P}[i,j] = 
  \begin{cases}
   \frac{d(i)}{2m}\frac{1}{d(i)} = \frac{1}{2m} & \text{if } (i,j)\in E  \\
   0       & \text{if } (i,j)\notin E
  \end{cases}
\]
Thus 
\[
Q(S,\overline S) = \frac{E(S:\overline S)/2m}{d(S)/2m}
\]
where $E(S:\overline S)$ denotes the number of edges with one end in $S$ and the other in $\overline S$ and $d(S)=\sum_{i\in S} d(i)$. 

Hence
\begin{equation}
\Phi = \Phi(G)=\min_{S \subseteq V: \pi(S)\leq 1/2}\frac{E(S:\overline S)}{d(S)}, \label{graphCond}
\end{equation}
where $\pi(S) = \sum_{i\in S}\p(i) = d(S)/2m$ is the probability of a random walk on $G$ being in $S$ when it is in the stationary distribution.

To glean some intuition behind equation (\ref{graphCond}), observe that in the stationary distribution, each of the $2m$ total orientations of the edges has the same probability of being transitioned. The quantity $E(S:\overline S)$ counts the number of orientations out of $S$, and $d(S)$ counts the total number of orientations that start in a vertex of $S$. Thus $E(S:\overline S)/d(S)$ gives the probability, in the stationary distribution, of moving out of $S$ at a given step given that the walk was in $S$. Or said in another way, when in the stationary distribution, it is the frequency of moving out of $S$, divided by the frequency of being in $S$. Intuitively, therefore, it would appear that a higher conductance might imply more rapid mixing of a random walk. Indeed, this is the case, as can be seen from this important and useful theorem, which was independently proved by \cite{JerrumSinclair} and \cite{LawlerSokal}:
\\

\begin{theorem}[\cite{JerrumSinclair}]\label{JSCondTh}
Let $\lambda_2$ be the second largest eigenvalue of a reversible, aperiodic transition matrix $\mathbf{P}$. Then
\begin{equation}
\frac{\Phi^2}{2} \leq 1-\lambda_2 \leq 2\Phi \label{JSCondineq}
\end{equation}
\end{theorem}

To be able to use Theorem \ref{JSCondTh} on a graph $G$, it needs to be non--bipartite so that it is aperiodic (see section \ref{randwalksundirectedgraphs}, Lemma \ref{nonBipartiteErgo}). Furthermore, to use it in conjunction with inequality (\ref{specdistbound}), we need to make sure that $\lambda_* = \lambda_2$. Both of these problems can be solved if we make the random walk \emph{lazy}. This means replacing the transition matrix $\mathbf{P}$ with $\mathbf{L}=\frac{1}{2}\mathbf{P}+\frac{1}{2}\mathbf{I}$, where $\mathbf{I}$ is the identity matrix. This introduces a loop probability of $1/2$. It means that asymptotically, the cover time becomes precisely twice as large. Introducing this loop probability makes the graph non-bipartite, but moreover, making it have probability (at least) $1/2$ ensures that all eigenvalues are non-negative. This is easy to see: If $\mathbf{x}$ is an eigenvector of $\mathbf{P}$ with eigenvalue $\lambda$, then
\[
\mathbf{L}\mathbf{x}=\left(\frac{1}{2}\mathbf{P}+\frac{1}{2}\mathbf{I}\right)\mathbf{x} = \frac{1}{2}\left(\mathbf{P}\mathbf{x}+\mathbf{I}\mathbf{x}\right) =
 \frac{1}{2}\left(\lambda \mathbf{x}+\mathbf{x}\right) = \frac{1}{2}\left(\lambda +1\right)\mathbf{x}
\]
and $\lambda+1\geq 0$ since $\lambda\geq -1$ by the Perron-Frobenius theorem (see (\ref{PerronFrob})). Furthermore, since $\mathbf{x}$ was an arbitrary eigenvector of $\mathbf{P}$, all the eigenvectors of $\mathbf{P}$ are eigenvectors of $\mathbf{L}$. This implies $\mathbf{P}$ and $\mathbf{L}$ have the same eigenvectors $\mathbf{v_1}, \mathbf{v_2}, \ldots, \mathbf{v_n}$, and the eigenvector $\mathbf{v_i}$ with eigenvalue $\lambda_i$ under $\mathbf{P}$ has eigenvalue $\frac{1}{2}\left(\lambda_i +1\right)$ under $\mathbf{L}$.

Thus, for the lazy walk $\mathbf{L}$, $\lambda_2\geq |\lambda_n|$ and so $\lambda_* = \lambda_2$. Thus, using inequality (\ref{JSCondineq}) in conjunction with inequality (\ref{specdistbound}), we have, for any $i,j$ 
\begin{equation}
|p_{i,j}^{(t)} -\pi_j| \leq  \sqrt{\frac{d(j)}{d(i)}}\lambda^t_*  \leq  \sqrt{\frac{d(j)}{d(i)}}\left(1-\frac{\Phi^2}{2}\right)^t\label{distcondbound}
\end{equation}

Our reason for expressing a probability distance bound in terms of conductance is that it is sometimes far easier to bound conductance than the eigenvalues of a graph. This is the case with the graph space we analyse in this chapter. In section \ref{ConductanceLowerbound}, we shall prove that, the conductance of a random graph with given degree sequence is at least $1/100$, \whp. This will imply that a random walk on such a graph is rapidly mixing \whp, and this is crucial to the proof of Theorem \ref{CovTh}. The meaning of ``rapidly mixing'' in this context shall be made precise in due course.

\section{Random Graphs of a Given Degree Sequence: Structural Aspects}
\subsection{The Configuration Model}
\label{TheConfigurationModelSection}
The configuration model is a random process that has proved useful for studying graphs with prescribed degree sequences. We assume the necessary condition that the sum of the degrees is even, $\sum_{i=1}^{n}d_i = 2m$ and with each vertex $i \in V$, we associate $d_i$ half-edges, or \emph{stubs}, which we consider distinguishable. Starting with an arbitrary stub from amongst the $2m$, we choose another stub in the graph \uar\ and pair the two. We repeat this process, taking an arbitrary stub and pairing it with another chosen \uar\ from the remaining $2m-3$. We continue this way until all stubs have been paired. Since there are an even number of stubs in total, the process must terminate successfully. The set of pairings that results is called a \emph{configuration}. A configuration $C$ maps to a graph $G(C)$ on the same vertex set and with each stub pairing considered to constitute an edge in $G(C)$. Note that there will be multiple configurations mapping to the same graph (i.e., the mapping from configurations to graphs is many--to--one).

In the above process since every other stub was picked at random with no regard for which vertex it was associated with, we may have loops and/or multiple edges in the resulting graph. It is quite clear that because the stubs are distinguishable, each possible configuration has the same probability, with the number of configurations being 
\[
(2m-1)!! = (2m-1)(2m-3)...(1) = \frac{(2m)!}{2^mm!}.
\]
(If we think of stubs as vertices themselves, we may consider this process as picking a perfect matching on the stubs uniformly at random from all possible matchings). For a degree sequence $\mathbf{d}$, we shall write $\mathbb{CM}(\mathbf{d})$ for the configuration space for $\mathbf{d}$. That is, $\mathbb{CM}(\mathbf{d})$ is the set all possible $\frac{(2m)!}{2^mm!}$ configurations on the set of vertices $V$ with degree sequence $\mathbf{d}$, and a uniform distribution on them. Thus, each configuration $C\in \mathbb{CM}(\mathbf{d})$ is picked with probability $\frac{2^mm!}{(2m)!}$.

The number of configurations  mapping to a particular graph with the degree sequence in question is not uniform in general, that is, for graphs $G, H$, the sets $\{C \in \mathbb{CM}(\mathbf{d}) : G(C)=G\}$ and $\{C \in \mathbb{CM}(\mathbf{d}) : G(C)=H\}$ may be of different cardinality. However, each simple graph with the prescribed degree sequence $\textbf{d}$ corresponds to $\Pi_{i=1}^{n}(d_i)!$ configurations. Thus, conditioning on the outcome of the process producing a simple graph, we have a uniform distribution.

The point of the configuration model is that it tends to be much easier to analyse and prove statements with than a direct analysis of simple random graphs of a given degree sequence. However this is of little use if statements proved in the paradigm of the configuration model cannot be carried over to statements about the graph space of interest. We use the following principle.
\\

\begin{proposition}\label{configToSimpleProp}
Suppose that for a configuration space $\mathbb{CM}(\mathbf{d})$ on $n$ vertices there is a function 
\[
f(n)\leq \Pr(C \text{ picked $\uar$ from $\mathbb{CM}(\mathbf{d})$ is simple})
\]
and suppose that a statement $\mathcal{P}(\mathbf{d})$  proved in the configuration model (and thus parameterised on $\mathbf{d}$) is false with some probability at most $g(n)$. If $g(n)=o(f(n))$, then $\mathcal{P}(\mathbf{d})$ holds \whp\ when we condition on the configuration $C$ drawn from  $\mathbb{CM}(\mathbf{d})$ being simple (and thus mapping to a simple graph).
\end{proposition}

Each simple graph $G$ is mapped to by $\Pi_{i=1}^{n}(d_i)!$ configurations, but the actual probability of a simple graph - that is, the total probability of the subspace of configurations which map to simple graphs, is difficult to determine precisely. Estimates are given, e.g.\ in \cite{McW},  but these require restrictions on the degree sequence. Further details will be given in section \ref{seca}.

It should be noted that \cite{Durrett} gives a similar result, demonstrating that there exists a constant that is a lower bound for the conductance. The analysis in turn relies on results from \cite{Gkants2003}.

\subsection{Conductance: A Constant Lower Bound}\label{ConductanceLowerbound}

Let \(\textbf{d} = (d_1, d_2, ..., d_n)\) be a sequence of natural numbers \ignore{for which the sum \(2M = \sum_{i=1}^{n}d_i\) is even}and let \(G=(V,E)\) be a graph of $n$ vertices chosen \uar\ from the family $\mathcal{G}(\mathbf{d})$ of all simple graphs with degree sequence \(\textbf{d}\), i.e., such that $d_i$ denotes the degree of vertex $i$.

We make the following assumptions about the degree sequence:
\begin{description}
	\item[(i)] $\sum_{i=1}^{n}d_i = 2m$ where $m$ is a natural number.
	\item[(ii)] The minimum degree $\delta \geq 3$.
	\item[(iii)] The average degree $\theta  = 2m/n \leq n^{\zeta}$, where \(0<\zeta < 1/3\) is a constant.	 
\end{description}

We work in the configuration model and make the following further assumption:

\begin{description}
	\item[(iv)] $\mathbf{d}$ is further restricted in such a way that statements which fail with probability at most $n^{-\Omega(1)}$ in the configuration model hold \whp\ when the model is conditioned on mapping to a simple graph (in other words, we can apply Proposition \ref{configToSimpleProp}).	 
\end{description}

As per the statement of Theorem \ref{CovTh}, in this chapter the cover time is analysed for \emph{nice} sequences. As will be seen in section \ref{seca}, a nice sequence $\mathbf{d}$ has the property that $\Pr(C \in \mathbb{CM}(\mathbf{d}) \text{ is simple}) \geq e^{-o(\log n)}$, implying that if a statement is demonstrated to fail in the configuration model with probability at most $n^{-\Omega(1)}$, then by Proposition \ref{configToSimpleProp}, it holds \whp\ when we condition for simple graphs.

The statements in the proof of Theorem \ref{condLowerboundTheorem} do indeed fail with probability at most $n^{-\Omega(1)}$, and so for the purposes of Theorem \ref{CovTh}, the proof of Theorem \ref{condLowerboundTheorem} is valid. Nevertheless, some of the conditions of nice sequences are specified for the sake of the analysis of cover time rather than conductance, and so we leave the details to section \ref{seca} and present a proof for Theorem \ref{condLowerboundTheorem} that depends on the more general assumptions \textbf{(i)}-\textbf{(iv)}.
\\

\begin{theorem}\label{condLowerboundTheorem}
Subject to assumptions \textbf{(i)}-\textbf{(iv)}, for a graph $G \in \mathcal{G}(\mathbf{d})$, $\Phi(G)>1/100$ \whp. \label{condtheoremjkndfsk}
\end{theorem}
Before we proceed with the proof, we note an immediate corollary:
\\

\begin{corollary}
$G \in \mathcal{G}(\mathbf{d})$ is connected \whp.
\end{corollary}

\proofstart
For a set \(S \subseteq V\) let \(d(S) = \sum_{v\in S}d(v)\) and for a configuration $C$,  let \(E_C(V_1:V_2)\) denote the number of edges with one end in \(V_1\) and the other in \(V_2\) in $C$. Let $\mathcal{E}_C(S) = E_C(S:\overline S)/d(S)$. Let $\pi(S) = d(S)/2m$.

We work in the configuration model. Our general approach is to show that when a configuration $C\in \mathbb{CM}(\mathbf{d})$ is picked, the set 
\[
\#(C)=\{S \subseteq V : \pi(S) \leq 1/2,\, \mathcal{E}_C(S)\leq 1/100\}
\]
is empty with probability at least $1-n^{-\Omega(1)}$. Having done this, the theorem follows by condition \textbf{(iv)}.

Throughout we make use of the following result of Stirling (see, e.g., \cite{Fe}):
\begin{equation}
n! = \sqrt{2\pi n}\left(\frac{n}{e}\right)^ne^{\lambda_n} \label{stirlingprecise}
\end{equation}
where 
\begin{equation}
\frac{1}{12n+1}<\lambda_n<\frac{1}{12n}. \label{lambda_n}.
\end{equation}
\eqref{lambda_n} implies $1<e^{\lambda_n} <1.1$ and that $e^{\lambda_n} =1+O(1/n)$. For notational convenience, we shall omit the correcting factor. 

A useful application of \eqref{stirlingprecise} is for fractions of the form \((2k)!/{k!}\) whence we get  $\sqrt{2}\left(\frac{4k}{e}\right)^k$. 

Let $X=|\#(C)|$ when $C \in \mathbb{CM}(\mathbf{d})$ is picked. Let \(\beta = 99/100\), \(\varepsilon = 1-\beta \), $\mathcal{F}(2K) = \frac{(2K)!}{K!2^{K}}$
and
\begin{equation}
\mathcal{H}(S) = \binom{d(S)}{\lceil \beta d(S)\rceil^*}\frac{\mathcal{F}(\lceil \beta d(S)\rceil^*)\mathcal{F}(2m-\lceil \beta d(S)\rceil^*)}{\mathcal{F}(2m)}, \label{H(S)}
\end{equation}
where $\lceil \beta d(S)\rceil^*$ is the smallest even integer greater than or equal to $\beta d(S)$. For notational convenience, we omit the ceiling and $*$ symbols and note that doing so can only incur (small) constant correcting factors that will not affect the results.

$\mathcal{H}(S)$ is an upper bound on the probability that a particular set of vertices $S$ will have $\mathcal{E}_C(S)\leq 0.01$ when $C$ is picked. 

By linearity of expectation
\begin{equation}
\E[X] \leq  \sum_{\substack{S\subseteq V,\\ \pi(S) \leq 1/2}}\mathcal{H}(S).\label{violationsum464557}
\end{equation}

Letting $\partial =d(S)$, the RHS of \eqref{H(S)} can be expanded thus:

\begin{flalign}
\binom{\partial}{\beta \partial}\frac{\mathcal{F}(\beta \partial)\mathcal{F}(2m-\beta \partial)}{\mathcal{F}(2m)}
& = \frac{(\partial)!}{(\beta \partial)!(\partial-\beta \partial)!} \frac{(\beta \partial)!}{(\beta \partial/2)!}\frac{(2m-\beta \partial)!}{(m-\beta \partial/2)!}\frac{m!}{(2m)!} \nonumber \\
& = \frac{\partial!}{(\varepsilon \partial)!(\beta \partial/2)!}\frac{(2m-\beta \partial)!}{(m-\beta \partial/2)!}\frac{m!}{(2m)!} \label{84518463}
\end{flalign}

By Stirling:
\begin{flalign}
\frac{\partial!}{(\varepsilon \partial)!(\beta \partial/2)!} & \approx \frac{\sqrt{2\pi \partial}\partial^{\partial}e^{-\partial}}{\sqrt{2\pi\varepsilon \partial}(\varepsilon \partial)^{\varepsilon \partial}e^{-\varepsilon \partial}\sqrt{\pi \beta \partial}(\beta \partial/2)^{\beta \partial/2}e^{-\beta \partial/2}} \nonumber \\
& =\frac{1}{\sqrt{\pi \varepsilon \beta \partial}} \frac{1} {e^{\beta \partial/2}}\frac{\partial^{\partial}}{(\varepsilon \partial)^{\varepsilon \partial}(\beta \partial/2)^{\beta \partial/2}} \label{7g87g5d}
\end{flalign}
and
\begin{flalign}
\frac{(2m-\beta \partial)!}{(m-\beta \partial/2)!} & \approx \sqrt{2}\left(\frac{4m-2\beta \partial}{e}\right)^{m-\beta \partial/2},\label{juykj5f} \\
\frac{m!}{(2m)!} & \approx \frac{1}{\sqrt{2}}\left(\frac{e}{4m}\right)^m. \label{54g8h7jf8a}
\end{flalign}
We substitute \eqref{7g87g5d},  \eqref{juykj5f} \eqref{54g8h7jf8a} and into \eqref{84518463}; observe there are seven factorial terms and we absorb the constant factor corrections from Stirling's approximation and the dropping of the ceiling and $*$ symbols into a constant $K_1$. So
\begin{flalign}
\mathcal{H}(S) 
&\leq K_1 \frac{e^{-\beta \partial/2}}{\sqrt{\pi \varepsilon \beta \partial}}\frac{\partial^{\partial}}{(\varepsilon \partial)^{\varepsilon \partial}(\beta \partial/2)^{\beta \partial/2}}\left(\frac{4m-2\beta \partial}{e}\right)^{m-\beta \partial/2}\left(\frac{e}{4m}\right)^m \nonumber \\
& =\frac{K}{\sqrt{\partial}}\frac{\partial^{\partial}}{(\varepsilon \partial)^{\varepsilon \partial}(\beta \partial/2)^{\beta \partial/2}}\frac{1}{(4m-2\beta \partial)^{\beta \partial/2}}\left(\frac{4m-2\beta \partial}{4m}\right)^{m}\nonumber \\ 
& = \frac{K}{\sqrt{\partial}}\left[\frac{\partial}{(\varepsilon \partial)^{\varepsilon} (\beta \partial/2)^{\beta /2}}\frac{1}{(4m-2\beta \partial)^{\beta /2}}\right]^{\partial}\left(\frac{2m-\beta \partial}{2m}\right)^{m}\nonumber \\ 
& = \frac{K}{\sqrt{\partial}}\left[\frac{1}{{\varepsilon}^{\varepsilon} {\beta}^{\beta /2}}\left(\frac{\partial}{2m-\beta \partial}\right)^{\beta /2} \right]^{\partial} \left(\frac{2m-\beta \partial}{2m}\right)^{m} \label{9r039ujfwowu98}
\end{flalign}
where $K$ is a constant.

Call a set \(S\) \emph{small} if \(d(S) \leq ({\theta n})^{1/4}\), otherwise call it \emph{large}. We handle the cases of small and large sets separately, and let the random variables $Y$ and $Z$ count for them respectively, i.e., $X=Y+Z$. We show that $E[Y]\leq n^{-\Omega(1)}$ and $E[Z]\leq n^{-\Omega(1)}$ and therefore, by Markov's inequality, $\Pr(X>0)\leq \E[Y]+\E[Z]\leq  n^{-\Omega(1)}$.

\subsubsection{Small Sets}

We bound the part of the sum of \eqref{violationsum464557} for small sets, and we partition it into sets of equal size (in the number of vertices). Since \(d(S)\leq ({\theta n})^{1/4}\), \(|S| = o(n)\). 

\begin{eqnarray*}
\E[Y] \leq \sum_{\substack{S\subseteq V,\\ d(S) \leq ({\theta n})^{1/4}}}\mathcal{H}(S) =\sum_{i=1}^{o(n)}\sum_{\substack{S\subseteq V,\\ d(S) \leq ({\theta n})^{1/4},\\|S|=i}}\mathcal{H}(S)
\leq  \sum_{i=1}^{o(n)}\binom{n}{i}\mathcal{H}(S_i)
\end{eqnarray*}

where \(S_i\) is the set of size \(i\) for which \(\mathcal{H}(S)\) is greatest over all sets \(S\) over the range of the sum with \(|S|=i\). 

Let \(S'\) with \(|S'|=s\) and \(\partial'=d(S')\) be the set over the range of the sum for which \(\binom{n}{i}\mathcal{H}(S_i)\) is greatest. Then 
\[
\sum_{i=1}^{o(n)}\binom{n}{i}\mathcal{H}(S_i) \leq o(n)\binom{n}{s}\mathcal{H}(S') \leq o(n)\left(\frac{en}{s}\right)^s\mathcal{H}(S'). 
\]
Using \eqref{9r039ujfwowu98},
\begin{eqnarray}
\E[Y] &\leq& o(n) \left(\frac{en}{s}\right)^s \frac{K}{\sqrt{\partial'}}\left[\frac{1}{{\varepsilon}^{\varepsilon} {\beta}^{\beta /2}}\left(\frac{\partial'}{2m-\beta \partial'}\right)^{\beta /2} \right]^{\partial'} \left(\frac{2m-\beta \partial'}{2m}\right)^{m} \nonumber\\
&=&o(n) \left(\frac{en}{s}\right)^s \left[\frac{1}{{\varepsilon}^{\varepsilon} {\beta}^{\beta /2}}\left(\frac{1}{2m/\partial'-\beta}\right)^{\beta /2} \right]^{\partial'}\label{ljnvoisg564548}
\end{eqnarray}

\(2m=\theta n\) and \(\partial' \leq (\theta n)^{1/4}\) so \(2m/\partial' \geq (\theta n)^{3/4}\), hence,
\[
\frac{1}{2M/\partial'-\beta} \leq \frac{1}{(\theta n)^{3/4}-\beta} \leq \frac{1.1}{(\theta n)^{3/4}}
\]
for large enough $n$. Hence,
\[
\frac{1}{{\varepsilon}^{\varepsilon} {\beta}^{\beta /2}}\left(\frac{1}{2m/\partial'-\beta}\right)^{\beta /2} \leq \left(\frac{1.1}{{\varepsilon}^{2\varepsilon/\beta}\beta(\theta n)^{3/4}}\right)^{\beta/2} \leq \left(\frac{1.3}{(\theta n)^{3/4}}\right)^{\beta/2}.
\]
Letting $\rho=\partial'/s$,
\begin{eqnarray}
\E[Y] &\leq& o(n)\left[\frac{en}{s}\left(\frac{2}{\theta n}\right)^{\frac{3}{4}\frac{\beta \rho}{2}} \right]^{s} \nonumber\\
 &\leq& \left[n^{1/s}en\left(\frac{2}{\theta n}\right)^{3\beta \rho/8} \right]^{s} \label{fsjfovpoakgvsgs} \\
 &\leq& \left[en^{1/s+1-3\beta \rho/8} \right]^{s}\label{f4e68g4g54hd}
\end{eqnarray} 
where \eqref{f4e68g4g54hd} follows from \eqref{fsjfovpoakgvsgs} because $\theta > 2$.

Now \(\rho \geq 3\) so \(1/s+1-3\beta \rho/8 \leq -0.01375\) for \(s\geq 10\). For sets of size \(s< 10\), we can do away with the \(o(n)\) term that multiplies the sum and replace it by a constant. This means that the exponent of \(n\) is \(1-3\beta \rho/8 \leq -0.11\). Thus, $\E[Y] \leq n^{-\Omega(1)}$.

\subsubsection{Large sets}
We now consider subsets $S$ for which $d(S) \geq (\theta n)^{1/4}$. Let $d(S) = \rho(S) c n = \alpha(S) \theta n$ where  \(\rho(S) = d(S)/|S|\) and \(0 < \alpha(S) < \frac{1}{2}\). Let $\partial=d(S)$, $\alpha=\alpha(S)$,  and note $K/\sqrt{\partial}<1$ for large enough $n$.  Hence, by \eqref{9r039ujfwowu98},
\begin{eqnarray}
\mathcal{H}(S)  &\leq& \left[\frac{1}{{\varepsilon}^{\varepsilon} {\beta}^{\beta /2}}\left(\frac{\partial}{2m-\beta \partial}\right)^{\beta /2} \right]^{\partial} \left(\frac{2m-\beta \partial}{2m}\right)^{m}\nonumber \\
&=& \left[\frac{1}{{\varepsilon}^{\varepsilon} {\beta}^{\beta /2}}\left(\frac{\alpha \theta n}{\theta n-\beta \alpha \theta n}\right)^{\beta /2} \right]^{\alpha \theta n} \left(\frac{\theta n-\beta \alpha \theta n}{\theta n}\right)^{\frac{\theta n}{2}}\nonumber \\
&=&\left[\frac{1}{{\varepsilon}^{\varepsilon} {\beta}^{\beta /2}}\left(\frac{\alpha}{1-\alpha\beta}\right)^{\beta /2} \right]^{\alpha \theta n} \left(1-\alpha\beta\right)^{\frac{\theta n}{2}}\nonumber\\
&=&\left[\frac{(\alpha\beta)^{\alpha\beta}(1-\alpha\beta)^{1-\alpha\beta}}{({\varepsilon}^{\varepsilon}{\beta}^{\beta})^{2\alpha}}\right]^{\frac{\theta n}{2}} = f(S).\label{f86s84f6ds8f46x}
\end{eqnarray}

We split the proof of the large sets into two parts: Those sets for which \(\alpha \leq 1/\theta\) and those for which \(1/{\theta}\leq \alpha\ \leq 1/2\).

\underline{$\alpha \leq 1/\theta$}\\
Let \(S_c' \in \mathcal{S}_c = \{S\subset V:|S| = cn\}\) be such that \(f(S_c') \geq f(S)\) for any \(S\in \mathcal{S}_c\).
For a constant $0<c<1$, define the random variable $Z_c = \sum_{S\in \mathcal{S}_c}\mathbf{1}_{\#(C)}(S)$, (the indicator random variable $\mathbf{1}_{\#(C)}(S)=1$ if and only if $S \in \#(C)$). 

Then 
\[
\E[Z_c] = \sum_{S\in \mathcal{S}_c}\mathcal{H}(S) \leq \binom{n}{cn}f(S_c')
\]

Applying Stirling's approximation to \(\binom{n}{cn}\) we have
\begin{flalign*}
\binom{n}{cn} & = \frac{n!}{(cn)!(n-cn)!}\\
& \approx \frac{\sqrt{2\pi n} n^ne^{-n}}{\sqrt{2\pi cn} (cn)^{cn}e^{-cn}\sqrt{2\pi (1-c)n}((1-c)n)^{(1-c)n}e^{-(1-c)n}}\\
& \leq \frac{K_2}{\sqrt{c(1-c)n}}\left(\frac{1}{c^c(1-c)^{1-c}}\right)^n
\end{flalign*}
where $K_2$ is some constant (which we shall assume absorbs the correcting factors $e^{\lambda_n}$ in the Stirling approximation).

Hence 
\begin{flalign*}
\E[Z_c] & \leq \frac{K_2}{\sqrt{c(1-c)n}}\left(\frac{1}{c^c(1-c)^{1-c}}\right)^n f(S_c')\\
&=\frac{K_2}{\sqrt{c(1-c)n}}\left(\left(\frac{(\alpha\beta)^{\alpha\beta}(1-\alpha\beta)^{1-\alpha\beta}}{({\varepsilon}^{\varepsilon}{\beta}^{\beta})^{2\alpha}}\right)^{\frac{\theta}{2}}\frac{1}{c^c(1-c)^{1-c}}\right)^n
\end{flalign*}
Consider the function 
\[g(x) = x^x(1-x)^{1-x} \mbox{ , } 0\leq x\leq 1/2
\]
\(g(0) = 1\) and the function is monotonically decreasing with minimum  \(g(1/2) = 1/2\).

Since \(\delta \geq 3\), \(\rho cn=\alpha \theta n\) implies \(c \leq \alpha\theta/3\). Now \(\alpha \leq 1/\theta\) implies \(\alpha\theta/3 < 1/2\), therefore \(g(c) \geq g(\alpha\theta/3)\).\\
Hence
\begin{flalign*}
\E[Z_c] & \leq \frac{K_2}{\sqrt{c(1-c)n}}\left(\left(\frac{(\alpha\beta)^{\alpha\beta}(1-\alpha\beta)^{1-\alpha\beta}}{({\varepsilon}^{\varepsilon}{\beta}^{\beta})^{2\alpha}}\right)^{\frac{\theta}{2}}\frac{1}{(\alpha\theta/3)^{\alpha\theta/3}(1-\alpha\theta/3)^{1-\alpha\theta/3}}\right)^n\\
& = \frac{K_2}{\sqrt{c(1-c)n}}\left(\frac{(\alpha\beta)^{\alpha\beta\theta/2}(1-\alpha\beta)^{1-\alpha\beta\theta/2}}{(\alpha\theta/3)^{\alpha\theta/3}(1-\alpha\theta/3)^{1-\alpha\theta/3}}\frac{(1-\alpha\beta)^{\theta/2-1}}{({\varepsilon}^{\varepsilon}{\beta}^{\beta})^{\alpha\theta}}\right)^n.
\end{flalign*}
Consider the function \(h(x,c) = (cx)^x(1-cx)^{1-x}\) where \(0\leq c\leq 1\).
\begin{flalign*}
\ln(h(x,c)) &= x\ln(cx) + (1-x)\ln(1-cx)\\
\frac{\partial}{\partial c}\mbox{ }\ln(h(x,c)) &= x\left(\frac{1}{c} - \frac{1-x}{1-cx}\right) = 0 \mbox{ at } c=1\\ 
\frac{{\partial}^2}{\partial c^2}\mbox{ }\ln(h(x,c)) &= -x\left(\frac{1}{c^2} + \frac{x(1-x)}{(1-cx)^2}\right) < 0
\end{flalign*}
Therefore \(h(x,c) < h(x,1) = g(x)\).

So \(h(\alpha\beta\theta/2, 2/\theta) < g(\alpha\beta\theta/2) < g(\alpha\theta/3)\).

Hence
\[
\E[Z_c] \leq \frac{K_2}{\sqrt{c(1-c)n}}\left(\frac{(1-\alpha\beta)^{\theta/2-1}}{({\varepsilon}^{\varepsilon}{\beta}^{\beta})^{\alpha\theta}}\right)^n.
\]
Now
\[
\frac{(1-\alpha\beta)^{\theta/2-1}}{({\varepsilon}^{\varepsilon}{\beta}^{\beta})^{\alpha\theta}} = \frac{1}{1-\alpha\beta} \left(\frac{(1-\alpha\beta)^{\frac{1}{2}}}{({\varepsilon}^{\varepsilon}{\beta}^{\beta})^{\alpha}}\right)^{\theta},
\]
and
\[
\frac{\partial}{\partial\theta}\left\{\frac{1}{1-\alpha\beta}\left(\frac{(1-\alpha\beta)^{\frac{1}{2}}}{({\varepsilon}^{\varepsilon}{\beta}^{\beta})^{\alpha}}\right)^{\theta}\right\} = \frac{1}{1-\alpha\beta} \left(\frac{(1-\alpha\beta)^{\frac{1}{2}}}{({\varepsilon}^{\varepsilon}{\beta}^{\beta})^{\alpha}}\right)^{\theta}\ln\left(\frac{(1-\alpha\beta)^{\frac{1}{2}}}{({\varepsilon}^{\varepsilon}{\beta}^{\beta})^{\alpha}}\right).
\]
Now
\[
\frac{(1-\alpha\beta)^{\frac{1}{2}}}{({\varepsilon}^{\varepsilon}{\beta}^{\beta})^{\alpha}} = 1 \mbox{ at } \alpha = 0
\] and
\begin{flalign*}
\frac{\partial}{\partial\alpha}(1-\alpha\beta) & = -\beta\\
\frac{\partial}{\partial\alpha}({\varepsilon}^{\varepsilon}{\beta}^{\beta})^{2\alpha} & = ({\varepsilon}^{\varepsilon}{\beta}^{\beta})^{2\alpha} \ln ({\varepsilon}^{\varepsilon}{\beta}^{\beta})^2 > \ln ({\varepsilon}^{\varepsilon}{\beta}^{\beta})^2.
\end{flalign*}
Consider
\begin{flalign*}
\frac{d}{d\beta}\left\{\ln ({\varepsilon}^{\varepsilon}{\beta}^{\beta})^2+\beta\right\} & = \frac{d}{d\beta}\left\{\ln ({(1-\beta)}^{1-\beta}{\beta}^{\beta})^2+\beta \right\}\\
& = 2\ln\left(\frac{\beta}{1-\beta}\right) + 1 > 0 \mbox{ for } \beta > \frac{1}{2}
\end{flalign*}
and \(\ln ({\varepsilon}^{\varepsilon}{\beta}^{\beta})^2 > -\beta \) when \(\beta = 0.99\), hence \(\frac{(1-\alpha\beta)^{\frac{1}{2}}}{({\varepsilon}^{\varepsilon}{\beta}^{\beta})^{\alpha}} < 1\) for \(\alpha > 0\). Therefore, 
\[
\frac{\partial}{\partial\theta} \left\{\frac{(1-\alpha\beta)^{\theta/2-1}}{({\varepsilon}^{\varepsilon}{\beta}^{\beta})^{\alpha\theta}}\right\} < 0
\]
and since \(\theta \geq \delta \geq 3\), we have that
\[
\frac{(1-\alpha\beta)^{\theta/2-1}}{({\varepsilon}^{\varepsilon}{\beta}^{\beta})^{\alpha\theta}} \leq \frac{(1-\alpha\beta)^{\frac{1}{2}}}{({\varepsilon}^{\varepsilon}{\beta}^{\beta})^{3\alpha}}\leq \frac{e^{-\alpha\beta/2}}{({\varepsilon}^{\varepsilon}{\beta}^{\beta})^{3\alpha}} < 0.73^{\alpha}
\]
Because \(\alpha\theta n \geq (\theta n)^{1/4}\) for large sets, we have that \(0.73^{\alpha} \leq 0.73^{(\theta n)^{-3/4}}\) and so
\[
\E[Z_c] \leq \frac{K_2}{\sqrt{c(1-c)n}}0.73^{\left(\frac{n}{{\theta}^3}\right)^{\frac{1}{4}}}.
\]
\(\frac{1}{n}\leq c\leq \frac{(n-1)}{n}\) so \(c(1-c) \geq \frac{n-1}{n^2}\), implying $\frac{1}{\sqrt{c(1-c)n}}\leq \sqrt{\frac{n}{(n-1)}}=1+o(1)$.

Therefore
\[
\E[Z_c] \leq (1+o(1))K_20.73^{\left(\frac{n}{{\theta}^3}\right)^{\frac{1}{4}}}.
\]

\(Z_c\) represents sets of size \(cn\) so multiplying the expression by \(n\) for an upper bound over all set sizes gives an upper bound of $O(n)0.73^{\left(\frac{n}{{\theta}^3}\right)^{\frac{1}{4}}}$. Furthermore, by our assumptions $\theta \leq n^{\zeta}$, where $0<\zeta < 1/3$ is a constant, and so \(\frac{n}{{\theta}^3} \geq n^{\zeta'}\) where $\zeta'$ is a positive constant. This implies $O(n)0.73^{\left(\frac{n}{{\theta}^3}\right)^{\frac{1}{4}}} \leq n^{-\Omega(1)}$.

\underline{$1/{\theta}\leq \alpha\ \leq 1/2$}\\
Let \(\Lambda = \{S \subset V :  1/{\theta}\leq  \alpha(S) = d(S)/{\theta n} \leq 1/2\}\). Let \(S_{\Lambda} \in \Lambda\) be such that  \(f(S_{\Lambda}) \geq f(S)\) for any \(S \in \Lambda\). Then 
\[
\E\left[\sum_{S \in \Lambda}\mathbf{1}_{\#(C)}(S)\right] \leq 2^nf(S_\Lambda)
\]
Using \eqref{f86s84f6ds8f46x} with $\alpha=\alpha(S_{\Lambda})$,
\begin{equation*}
f(S_{\Lambda}) =  \left(\frac{(\alpha\beta)^{\alpha\beta}(1-\alpha\beta)^{1-\alpha\beta}}{({\varepsilon}^{\varepsilon}{\beta}^{\beta})^{2\alpha}}\right)^{\frac{\theta n}{2}}. 
\end{equation*}
Let
\[
A(\alpha) = \frac{(\alpha\beta)^{\alpha\beta}(1-\alpha\beta)^{1-\alpha\beta}}{({\varepsilon}^{\varepsilon}{\beta}^{\beta})^{2\alpha}}
\] 
\begin{flalign*}
\ln(A(\alpha)) & = (\alpha\beta)\ln((\alpha\beta)) + (1-\alpha\beta)\ln(1-\alpha\beta) - 2{\alpha}\ln({\varepsilon}^{\varepsilon}{\beta}^{\beta})\\
\frac{\partial}{\partial\alpha}\ln(A(\alpha)) & = \beta \ln(\alpha\beta) - \beta\ln(1-\alpha\beta) - 2\ln({\varepsilon}^{\varepsilon}{\beta}^{\beta})
\end{flalign*}
Setting \(\frac{\partial}{\partial\alpha}\ln(A(\alpha)) = 0\) gives 
\[
\alpha = \frac{{\varepsilon}^{2\varepsilon/\beta}\beta}{1+{\varepsilon}^{2\varepsilon/\beta}{\beta}^2} = 0.477 \text{ to 3 d.p.}
\]

\[
\frac{{\partial}^2}{\partial{\alpha}^2}\ln(A(\alpha)) = \beta\left(\frac{1}{\alpha}+\frac{\beta}{1-\alpha\beta} \right) > 0
\]
therefore the stationary point is a minimum. Furthermore, $A(1/2)<A(1/3)$ so \(A(\alpha) \leq A(1/\theta)\). Hence,
\begin{eqnarray*}
2^nf(S_\Lambda) &=& (2(A(\alpha))^{\frac{\theta}{2}})^n \\
&\leq& (2(A(1/\theta))^{\frac{\theta}{2}})^n\\
&\equiv& \left[2\left(\frac{\left(\beta/\theta\right)^{\frac{\beta}{\theta}}\left(1-\beta/\theta\right)^{1-\frac{\beta}{\theta}}}{({\varepsilon}^{\varepsilon}{\beta}^{\beta})^{2/\theta}}\right)^{\frac{\theta}{2}}\right]^n\\
 &=& \left[2\left(\frac{(\beta/\theta)^{\frac{\beta}{2}}(1-\beta/\theta)^{\frac{1}{2}\left(\theta- \beta\right)}}{{\varepsilon}^{\varepsilon}{\beta}^{\beta}}\right)\right]^n.
\end{eqnarray*}
Let 
\[
T(\theta) = \left({\frac{\beta}{\theta}}\right)^{\beta}\left(1-\frac{\beta}{\theta}\right)^{\theta- \beta}
\]
Then
\[
\ln (T(\theta)) = \beta\ln\beta - (\theta - \beta)\ln(\theta - \beta) - \theta\ln\theta
\]
and
\[
\frac{\partial}{\partial\theta}\ln (T(\theta)) = \ln\left(\frac{\theta-\beta}{\theta}\right)<0.
\]
Then since \(\theta \geq 3\), we have that \(T(\theta) \leq T(3) \) and so
\[
2^nf(S_\Lambda) \leq \left[2\left(\frac{(\beta/3)^{\frac{\beta}{2}}(1-\beta/3)^{\frac{1}{2}\left(3- \beta\right)}}{{\varepsilon}^{\varepsilon}{\beta}^{\beta}}\right)\right]^n  \leq 0.8^n.
\]

This concludes the proof that $\E[Z]\leq n^{-\Omega(1)}$, and Theorem \ref{condtheoremjkndfsk} follows.
\proofend

\section{Assumptions About Degree Sequence}\label{seca}
We are studying the cover time of a graph $G$ picked uniformly at random \uar\ 
from the set \({\cal G}(\textbf{d})\) of simple graphs with vertex set \(V=[n]\) and degree sequence \(\textbf{d}
=(d_1,d_2,\ldots,d_n)\).
We make the following definitions:
Let \(V_j = \{i \in V: d_i = j\}\) and let \(n_j = |V_j|\).
Let \(\sum_{i=1}^{n}d_i = 2m\) and let
\(\th=2m/n\) be the average degree.

It seems reasonable to ask how the various entries in the  degree sequence affect the cover time.
In particular, how much does the cover time depend on the vertices of low degree, and how much on the average degree
of the graph? In fact, as in \cite{Coppersmith}, both parameters play a part, as is shown in Theorem \ref{CovTh}.

Let \(0<\alpha<1\) and \(0<\kappa<1/11\) be constants, and let \(d\) be a positive integer.
Let \(\g \rai\) with \(n\).
We suppose the degree sequence \textbf{d} satisfies
the following conditions:
\begin{description}
\parskip 1ex
\item[(i)] Average degree \(\th=o(\sqrt{\log n})\).
\item[(ii)]Minimum degree \(\delta \geq 3\).
\item[(iii)]For \(\delta \leq i < d\), \(n_i =O( n^{\kappa i/d})\).
\item[(iv)] \(n_d=\a n \pm o(n)\). We call \(d\) the {\em effective minimum degree}.
\item[(v)] Maximum degree \(\Delta =O(n^{\kappa(d-1)/d})\).
\item[(vi)] Upper tail size {\(\displaystyle \sum_{j=\g\th}^\D n_j=O(n^{\kappa(d-1)/d})\)}.
\end{description}

Immediately on fixing the degree sequence, \(\textbf{d}\),
some definitional problems arise, as  e.g. there may be
just a few low degree vertices spread over a wide range.
To get round this, we define an {\em absolute} minimum degree $\delta$, and an \emph{effective} minimum degree
\(d\). The effective minimum degree is the first entry in the sorted degree
sequence which occurs order \(n\) times. We fix the minimum degree at 3
to ensure the graph is connected (\whp). Between the minimum degree,
and the effective minimum degree, we place an upper bound on the
number \(n_i\) of vertices of degree \(i\). The bound we choose of \(n_i =O( n^{\kappa i/d})\),
is not as arbitrary as it looks. Certainly when \(\kappa=d/(d-1)\), the effective minimum
degree drops below \(d\), so clearly there is some \(\kappa < d/(d-1)\) which is critical.
Finally, we make some  constraints on the average degree and upper tail of the degree sequence,
to ensure simple graphs occur with high enough probability in the configuration model.
Thus we are left with the above list of conditions.

We call  a degree sequence  \textbf{d} which satisfies conditions (i)--(vi)  {\em nice}, and  apply the
same adjective to \(G\in {\cal G}(\textbf{d})\).

Examples of nice sequences/graphs are: Any $r$-regular graph (where $r$ is constant); The sequence which has at most $n^{1/20}$ vertices with degree at most $n^{1/20}$, and the rest are (constant) degree $d$. This will have average degree $\theta=(1+o(1))d$, but possibly high maximum degree; The sequence which has $\Theta(n^{1/20})$ vertices with degree $\Theta(n^{1/20})$, and the rest are equally divided between (constant) degree $d$ and degree $\Theta(\sqrt{\log n}/\log\log n)$. This will have a high average degree $\Theta(\sqrt{\log n}/\log\log n)$ as well as a large number of high degree vertices.

Recall from section \ref{TheConfigurationModelSection} we discussed the set of all configurations with degree sequence $\mathbf{d}$ and with uniform probability, known as the configuration space $\mathbb{CM}(\mathbf{d})$. That section discussed the probability of a configuration $C\in \mathbb{CM}(\mathbf{d})$ being simple. It was stated that each simple graph $G$ with degree sequence $\mathbf{d}$ had $\Pi_{i=1}^{n}(d_i)!$ configurations mapping to it, thus, conditioning on simplicity, the distribution across simple graphs was uniform. However, the actual probability of simplicity was not elaborated on. A result of \cite{McW} gives asymptotic values of the probability of simplicity, but requires certain constraints on $\mathbf{d}$.
  
Let \(\nu=\sum_i d_i(d_i-1)/(2m)\).
Assume that \(\D=o(m^{1/3})\) (as it will be for nice sequences). 
The probability that \(C\in\mathbb{CM}(\mathbf{d})\) maps to a simple graph is given by
\begin{equation}\label{psimple}
\Pr(C \text{ simple}) \sim e^{-\frac{\nu}{2}-\frac{\nu^2}{4}}.
\end{equation}

Observe that our assumptions (i)--(vi) that  \textbf{d}\ is nice imply that \(\n=o(\sqrt{\log n})\).
Indeed if \(\th=\sqrt{\log n}/\g^3\) where \(\g\to\infty\) then
\begin{equation*}\n\leq\frac{1}{\th n}\brac{\sum_{j=3}^{\g\th}n_jj^2+\sum_{j=\g\th}^{\D}n_jj^2}\leq
\frac{1}{\th n}\left(n\g^2\th^2+O\left(n^{3\kappa(d-1)/d}\right)\right)=o(\sqrt{\log n}).
\end{equation*}
All the \whp\ statements in the analysis below of the structure of random graphs fail with
probability
at most \(g(n) = n^{-\Omega(1)}\), whereas $\Pr(C \text{ simple})$ in \eqref{psimple}
is at least \(f(n) = e^{-o(\log n)}\). By Proposition \ref{configToSimpleProp}, this justifies our use of $\mathbb{CM}(\mathbf{d})$ to make statements about $\mathcal{G}(\mathbf{d})$.

\section{Estimating First Visit Probabilities}
\label{1stvisit}
\subsection{Convergence of the Random Walk}
\label{cover}
In this section $G$ denotes a fixed connected graph with $n$ vertices and $m$ edges.
A random walk $\cW_{u}$  is started from a vertex $u$.
Let $\cW_{u}(t)$ be the vertex
reached at step $t$ and let
$P_{u}^{(t)}(v)=\Pr(\cW_{u}(t)=v)$.

We assume that the random walk
$\cW_{u}$ on $G$ is ergodic i.e. $G$ is not bipartite. Thus, the
random walk $\cW_u$ has the steady state distribution $\pi$, where
$\pi_v=d(v)/(2m)$. Here $d(v)$ is the degree of vertex $v$.

Note $P_{u}^{(t)}(v)$ is the same as the quantity
$P^{(t)}_{u,v}$ defined in section \ref{MarkovChainFirstDef}. The former notation emphasises the point that the starting
vertex $u$ becomes irrelevant after the mixing time, and $P_{u}^{(t)}(v)$ becomes close to $\pi_v$.

\subsection{Generating Function Formulation}
We use the approach of \cite{CFreg}, \cite{CGweb}. 

Let $D(t)=\max_{u,x\in V}|P_{u}^{(t)}(x)-\pi_x| $, and
let $T$ be such that, for $t\geq T$
\begin{equation}\label{4}
\max_{u,x\in V}|P_{u}^{(t)}(x)-\pi_x| \leq n^{-3} .
\end{equation}
It follows from e.g. Aldous and Fill \cite{AlFi} that
$D(s+t)\leq
2D(s)D(t)$  and so for $k\geq 1$,
\begin{equation}\label{4a}
\max_{u,x\in V}|P_{u}^{(kT)}(x)-\pi_x| \leq
\frac{2^{k-1}}{n^{3k}}.
\end{equation}

Fix two vertices $u,v$.
Letting $h_t=P_{u}^{(t)}(v)$ for notational convenience, let
\begin{equation}\label{Hz}
H(z)=\sum_{t=T}^\infty h_tz^t
\end{equation}
generate $h_t$ for $t\geq T$. 

Next, considering the  walk $\cW_v$, starting
at $v$, let $r_t=\Pr(\cW_v(t)=v)$ be the probability  that this  walk
returns to $v$ at step $t = 0,1,...$.
Let
$$R(z)=\sum_{t=0}^\infty r_tz^t$$
generate $r_t$. Our definition of return involves $r_0=1$.

For $t\geq T$ let $f_t=f_t(u \rat v)$
be  the probability that the first visit of the walk $\cW_u$ to $v$
in the period $[T,T+1,\ldots]$ occurs at step $t$.
Let
$$F(z)=\sum_{t=T}^\infty f_tz^t$$
generate  $f_t$.
Then we have
\begin{equation}
\label{gfw} H(z)=F(z)R(z).
\end{equation}
Finally, for $R(z)$
 let
\begin{equation}
\label{Qs}
R_T(z)=\sum_{j=0}^{T-1} r_jz^j.
 \end{equation}

\subsection{First Visit Time Lemma: Single Vertex}
Let
\begin{equation}
\label{lamby}
\l=\frac{1}{KT}
\end{equation}
for some sufficiently large constant $K$.

The following lemma should be viewed in the context that $G$ is an $n$ vertex graph
which is part of a sequence of graphs with $n$ growing to infinity.
An almost identical lemma was first proved in \cite{CFreg}.
\\

\begin{lemma}
\label{L3}Suppose that
\begin{description}
\item[(a)]
For some constant $\psi>0$, we have
$$\min_{|z|\leq 1+\l}|R_T(z)|\geq \psi.$$
\item[(b)]$T^2\pi_v=o(1)$ and $T\pi_v=\Omega(n^{-2})$.
\end{description}
There exists
\begin{equation}\label{pv}
p_v=\frac{\pi_v}{R_T(1)(1+O(T\pi_v))},
\end{equation}
where $R_T(1)$ is from (\ref{Qs}), such that for all $t\geq T$,
\begin{equation}
\label{frat}
f_t(u \rat v)=(1+O(T\pi_v))\frac{p_v}{(1+p_v)^{t+1}}+
O(T\pi_ve^{-\l t/2}).
\end{equation}
\end{lemma}

We have the following corollary
\\

\begin{corollary}
\label{geom}
For $t\geq T$ let $\ul A_t(v)$ be the event that $\cW_u$
does not visit $v$ in steps $T,T+1,\ldots,t$.
Then, under the assumptions of Lemma \ref{L3},
\[
\Pr(\ul A_t(v))=
\frac{(1+O(T\p_v))}{(1+p_v)^{t}} +O(T^2\p_v e^{-\lambda t/2}).
\]
\end{corollary}

\proofstart
We use Lemma \ref{L3} and
$$\Pr(\ul A_t(v))= \sum_{\t > t} f_{\t}(u \rat v),$$
and note that $T^2 \pi_v=o(1)$.
\proofend

Rather than reproduce the proof of Lemma \ref{L3}, we instead present a simplified version which conveys some of the core ideas behind it in a more straightforward form. This simplified version gives a bound rather than a bound on $\Pr(\ul A_t(v))$ rather than the precise value given in Corollary \ref{geom}.

\textbf{NB} We use may use $R_v$ to denote $R_T(1)$, to emphasis the dependence of $R_T(1)$ upon $v$.

\subsection{First Visit Lemma: Simplification}\label{FVLSimp}
For the graph $G= (V,E)$ let $n=|V|$ and $m=|E|$.

\subsubsection{Visits to Single Vertices. }
Recall from section \ref{MarkovChainFirstDef} the definition of $\hit[\mathbf{p}, a]$ for a distribution $\mathbf{p}$ over states of a Markov chain:
\[
\hit[\mathbf{p}, a] = \sum_{i}\mathbf{p}_i\hit[i,a].
\]

The hitting time from the stationary distribution, $\hit[\pi, w] = \sum_{v\in V}\pi_v\hit[v,w]$, can be expressed
as $\hit[\pi, w] = Z_{w,w}/\pi_w$, where
\begin{equation} 
Z_{w,w} = \sum_{t=0}^{\infty} (P_{w}^{(t)}(w) - \pi_w), \label{Z_wwEquation}
\end{equation}
see e.g.\ \cite{AlFi}.

Let $P$ denote the transition matrix for a random walk on $G$,
and, for a walk $\cW_v$ starting from $v$ define
\begin{equation}
R_v(T) \; = \; \sum_{t=0}^{T-1}{P}_v^{(t)}(v).
\end{equation}
Thus $R_v$ is the expected number of returns made by $\cW_v$ to $v$
during $T$ steps, in the graph $G$. We note that $R_v \ge 1$, as
$P_v^{(0)}(v)=1$.

In subsequent discussions, we will assume $T$ is a mixing time
satisfying
\begin{equation}\label{mixing}
  |{P}_{u}^{(T)}(x)-{\pi}_x| \leq \frac{1}{n^3},
\end{equation}
for any vertices $u$ and $x$ in $G$. Let
$D(t)=\max_{u,x}|P_u^{(t)}(x)-\pi_x|$. As $\pi_x \ge 1/n^2$ for any
vertex of a simple graph, \eqref{mixing} implies that $D(T) \le
\pi_x$ for all $x \in V$.
\\

\begin{lemma}\label{Eval}
Let $T$ be a mixing time of a random walk $\cW_u$ on $G$ satisfying
\eqref{mixing}. Let vertex $v \in V$ be such that $T \cdot
\pi_v=o(1)$, and $\pi_v <1/2$, then
\begin{equation}
\hit[\pi,v] = \ooi \frac{R_v(T)}{\pi_v}.
\end{equation}
\end{lemma}
\proofstart

Let $D(t)=\max_{u,x}|P_u^{(t)}(x)-\pi_x|$. It follows from e.g.
\cite{AlFi} that $D(s+t) \le 2D(s)D(t)$. Hence, since $
\max_{u,x}|P_u^{(T)}(x)-\pi_x| \le \pi_v, $ then for each $k \ge 1$,
$ \max_{u,x}|P_u^{(kT)}(x)-\pi_x| \le (2\pi_v)^k. $ Thus
\begin{align}
Z_{v,v} & = \; \sum_{t=0}^{\infty} (P_{v}^{(t)}(v) - \pi_v) \notag \\
& \leq \;\sum_{t<T} (P_{v}^{(t)}(v) - \pi_v) + T \sum_{k\ge 1}
(2\pi_v)^k
 \notag
 \\ & = \; R_v(T) - T \cdot \pi_v + O(T \cdot \pi_v)
\\ &= R_v(T) \ooi.
\end{align}
The last inequality follows because $R_v(T) \ge 1$. \proofend

Let ${\ul A_t(v)}$ denote the event that $\cW_u$ does not visit $v$ in
steps $0,...,t$. We next derive a crude upper bound for
$\Pr({\ul A_t(v)})$ in terms of $\hit[\pi, v]$.
\\

\begin{lemma}\label{crude}
Let $T$ satisfy \eqref{mixing},  then
\[
\Pr({\ul A_t(v)}) \le  \exp\brac{\frac{-(1-o(1))\lfloor t/\tau_v\rfloor}{2}},
\]
where $\tau_v=T+2\hit[\pi, v]$.
\end{lemma}
\proofstart Let $\rho=\r(G,T, u)$ be the distribution of $\cW_u$ on $G$
after $T$, then
\[
\hit[\r, v]=\ooi \hit[\pi, v].
\]
Indeed,
\begin{eqnarray}
\hit[\r, v]&=&\sum_{w\in V} \r_w \hit[w, v]\nonumber \\
&=& \ooi\sum_{w\in V} \pi_w\hit[w, v]\nonumber \\
&=&\ooi \hit[\pi, v]\label{gggg6ash}
\end{eqnarray}

Let $h_{\r}(v)$ be the time to hit $v$ starting from $\r$. Then $\E[h_{\r}(v)]=\hit[\rho,v]$ so by Markov's inequality and using \eqref{gggg6ash},
\[
\Pr(h_{\r}(v)\ge 2 \hit[\pi, v]) \le \frac{\ooi }{2}.
\]

\ignore{Let $\t=T+2\hit[\pi, v]$. }By restarting the process at $\mathcal{W}(0)=u,
\mathcal{W}(\t_v), \mathcal{W}(2\t_v),...,\mathcal{W}(\rdown{t/\t_v}\t_v)$ we obtain
\[
\Pr({\ul A_t(v)}) \le  \bfrac{\ooi}{2}^{\rdown{t/\t_v}}.
\]
\proofend

\subsubsection{Visits to Vertex Sets}

We extend the results presented above to a subset of vertices. Let
$\emptyset\neq S\subseteq V$.

From graph $G$ we obtain a (multi)-graph $\G=\G_S$ by contracting
$S$ to a single vertex $\g$. Note that we retain  multiple edges and
loops in $\G_S$, so that $d(S)=d(\g)$ and that $m=|E(G)|=|E(\G)|$.
Let $\hat{\pi}$ be the stationary distribution of a random walk on
$\G$. If $v \not \in S$ then $\pi_v=\hat \pi_v$ and  also that
$\pi_S=\hat \pi_{\g}$.

Any walk $\cW_v$ from $v\not\in S$ to $S$ in $G$ with internal
vertices not in $S$ corresponds to an identical walk $\hat W_v$ from
$v$ to $\g$ in $G_S$, and both walks have exactly the same
probability. Thus we obtain
\begin{align}\label{eq:expectedHittingSet}
\hit[\pi,S] \; = \; \hit[\hat\pi,\g] 
\end{align}

Let $\hat{P}$ denote the transition matrix for a random walk on
$\G_S$, let $T_{\G}$ denote the mixing time for $\G_S$, and define
\begin{equation}\label{returnToSet}
R_{\g} \; = \; \sum_{t=0}^{T_{\G}-1} \hat{P}_\g^{(t)}(\g).
\end{equation}
Note that $R_{\g}$ is the number of returns to $\g$ in the
contracted graph $\G_S$, rather than a somewhat different notion of the number of returns to $S$ in $G$.
Note also that in the summation in~\eqref{returnToSet} we use
$T_{\G}$, the mixing time in $\G$, rather than the mixing time in
$G$.

For any vertices $u$ and $x$ in $\G$, we will assume
\begin{equation}\label{mixingInContracted}
  |\hat{P}_{u}^{(T_{\G})}(x)-\hat{\pi}_x| \leq n^{-2}.
\end{equation}
We will always choose mixing times $T$ satisfying $T \ge T_G$, and
$T \ge T_{\G}$. Because of this, the distributions $\r_G=\r(G,T,u)$
of $\cW_u$, and $\r_{\G}=\r(\G,T,u)$ of $\hat{\cW}_u$, satisfy
\[
\hit[\r_G,S] = \ooi \hit[\r_{\G},\g] = \ooi \hit[\pi, S] \; =
\; \ooi  \hit[\hat\pi, \g].
\]
We also note that provided $T_{\G}\hat{\pi}_{\g}=o(1)$ then Lemma
\ref{Eval} and Lemma \ref{crude} both apply to $\g$ in $\G$.

\section{Required Graph Properties}
Let
\begin{equation}\label{om}
\omega = \log\log\log n
\end{equation}
We use this value throughout our proof. 

\subsection{Mixing Time}
\label{MixingTimeSection}
We shall define the \emph{mixing time} to be the smallest $T$ such that \eqref{4} is satisfied.

From  inequality (\ref{distcondbound}) we have

\begin{equation}\label{mix}
|P_{u}^{(t)}(x)-\pi_x| \leq \left(\frac{\p_x}{\p_u}\right)^{1/2}\left(1-\frac{\Phi^2}{2}\right)^t.
\end{equation}

As explained in section \ref{conductancesectionprelims}, we require the walk to be lazy, meaning it only makes a
move to a neighbour with probability 1/2 at any step.
This halves the conductance but if $\Phi(G)>1/100$ (which will be the case \whp, by Theorem \ref{condLowerboundTheorem}) we can still set
\begin{equation}\label{mixx}
T = \omega^3\log n
\end{equation}
to satisfy \eqref{4}. The cover time is doubled. Asymptotically
the values $R_v$ are doubled too. Otherwise, it has a negligible effect on the
analysis and we will ignore this for the rest of the chapter and continue
as though there are no lazy steps.\\
Note that $\Delta = O(n^c)$ for some constant $c<1$ implies $T^2\pi_v = o(1)$ for any vertex $v$.

\subsection{Structural Properties}\label{Structural Properties}
\subsubsection{Some definitions}
\label{props}
We show some properties held by almost all graphs in $\mathcal{G}(\textbf{d})$. 
First we define some terms:
Let
\begin{equation}\label{ell}
\ell =B \log^2 n
\end{equation}
for some large constant $B$.

\begin{description}
\item[(i)] A vertex $v$ is \emph{light} if it has degree at most \(\ell\), otherwise it is \emph{heavy}.
\item[(ii)] A cycle $C$ or path $P$ is \emph{small} if it has at most $2\omega+1$ vertices, otherwise it is \emph{large}.
\item[(iii)] A small path is \emph{light} if all of it's vertices are light, otherwise it is \emph{heavy}.
\item[(iv)] A small cycle is \emph{light} if it has at most one vertex that is heavy, otherwise it is \emph{heavy}.
\end{description}

For a vertex $v$, let $G_v$ be the subgraph induced by the set of vertices within a distance $\omega$ of $v$. A vertex $v$ is 
\begin{description}
\item[(v)] \emph{locally tree-like} if the subgraph induced by all vertices within distance $\omega$ of $v$ is a tree,
\item[(vi)] \emph{$r$-regular} if it is locally tree-like and each vertex in the tree (possibly with the exception of $v$) has degree $r$,
\item[(vii)] \emph{$r$-tree-compliant} or \emph{$r$-compliant} if there exists a tree subgraph of $G_v$, $T_v$ rooted at $v$ in which both of the following are true
\begin{enumerate}
	\item In $T_v$, vertex $v$ has degree $d_v$ (that is, all edges incident on $v$ in $G_v$ are part of $T_v$). 
	\item Every other vertex in $T_v$ has degree $r$ in $T_v$ (where it may be the case that $r=d_v$).
\end{enumerate}
\end{description}

\subsubsection{The local graph $\Gamma_v$}\label{localgrahgamma}
For a light vertex $v$, define the graph $\Gamma_v$ as follows: If $u \in G_v$ is heavy, delete an edge $(u,w) \in G_v$ iff there is no path from $w$ to $v$ (inclusive of $w$) in $G_v$ that is light. $\Gamma_v$ is then taken to be the connected component of $v$ in $G_v$ after this process has completed. We denote by $\Gamma_v^{\circ}$ the union of the set of (heavy) vertices that have been pruned and the vertices distance $\omega$ from the root $v$
(when $v$ is heavy, we do not need to consider this process, as shall be seen in section \ref{ubsec}).

\subsubsection{Structural Properties}
\begin{lemma}\label{lightCPC}
With high probability,
\begin{description}
\item[(a)] No pair of light cycles are connected by a light path.
\item[(b)] No pair of light cycles intersect.
\end{description}
\end{lemma}

\proofstart
Before we prove the lemma, we show a useful inequality. For integer $x>0$, let $\mathcal{F}(2x) = \frac{(2x)!}{2^xx!}$, then

\begin{eqnarray*}
\frac{\mathcal{F}(\theta n - 2x)}{\mathcal{F}(\theta n)} &=& \frac{(\theta n - 2x)!}{\left(\frac{\theta n}{2}-x\right)!2^{\frac{\theta n}{2}-x}}\frac{\left(\frac{\theta n}{2}\right)!2^{\frac{\theta n}{2}}}{(\theta n)!}\\
&=& 2^x\left(\prod_{i=0}^{x-1}\frac{\theta n}{2}-i\right)\left(\prod_{i=0}^{2x-1}\theta n-i\right)^{-1}\\ 
&=& \left(\prod_{i=1}^{x}\theta n-2i+1\right)^{-1}\\ 
&\leq& \left(\frac{1}{\theta n-2x+1}\right)^x
\end{eqnarray*}
 
\textbf{(a)} The expected number of light cycle-path-cycle (c-p-c) subgraphs $\mu$ has the following upper bound:
\begin{equation} \mu \leq \sum_{a=3}^{2\omega+1}\sum_{b=3}^{2\omega+1}\sum_{c=1}^{2\omega+1}\binom{n}{a}\binom{n}{b}\binom{n}{c}\frac{(a-1)!}{2}\frac{(b-1)!}{2}c!ab\ell^{2(a+b+c-2)}\Delta^6\frac{\mathcal{F}(\theta n - 2(a+b+c+1))}{\mathcal{F}(\theta n)}\label{NoCPC}
\end{equation}
\\

\textbf{Explanation} This is an upper bound on the number of configurations with a light c-p-c, divided by the total number of configurations. We choose $a$ vertices for one cycle, $b$ vertices for the other and $c$ vertices for the path. We multiply by the number of ways of permuting vertices in the two cycles and the path. Every vertex but at most one in a cycle is light, hence has at most degree $\ell$. If there is a heavy vertex, it can have any degree up to $\Delta$. Hence, each light vertex has up to $\ell(\ell-1)$ ways to connect to a neighbour. We therefore multiply $\ell^2$ by itself $((a-1) + (b-1) + c)$ times for the number of ways to connect to a neighbour for the two cycles and the path. The remaining (unaccounted for), possibly heavy vertex in each cycle can connect in up to $\Delta(\Delta-1)$ ways to a neighbour in the cycle and $\Delta-2$ ways to a vertex on the path, hence $\Delta^6$ serves as an upper bound. Since a c-p-c subgraph connects $2(a+b+c+1)$ ports, the number of ways of connecting the remaining free ports in the configuration is counted by the function $\mathcal{F}(\theta n - 2(a+b+c+1))$.
\\

Thus, (\ref{NoCPC}) is bounded by
\begin{eqnarray*} &&\frac{\Delta^{6}}{\ell^6}\sum_{a=3}^{2\omega+1}\sum_{b=3}^{2\omega+1}\sum_{c=1}^{2\omega+1}\binom{n}{a}\binom{n}{b}\binom{n}{c}a!b!c!\ell^{2(a+b+c+1)}\frac{\mathcal{F}(\theta n - 2(a+b+c+1))}{\mathcal{F}(\theta n)}\\
&\leq& \frac{\Delta^{6}}{\ell^6}\sum_{a=3}^{2\omega+1}\sum_{b=3}^{2\omega+1}\sum_{c=1}^{2\omega+1}\binom{n}{a}\binom{n}{b}\binom{n}{c}a!b!c!\left(\frac{\ell^2}{\theta n-12\omega-7}\right)^{a+b+c+1}\\
&\leq& 
\frac{\Delta^{6}}{\ell^4}\frac{1}{\theta n-12\omega-7}\sum_{a=3}^{2\omega+1}\sum_{b=3}^{2\omega+1}\sum_{c=1}^{2\omega+1}\frac{n!}{(n-a)!}\frac{n!}{(n-b)!}\frac{n!}{(n-c)!}\left(\frac{\ell^2}{\theta n-12\omega-7}\right)^{a+b+c}\\
&\leq&
\frac{\Delta^{6}}{\ell^4}\frac{1}{\theta n-12\omega-7}\sum_{a=3}^{2\omega+1}\sum_{b=3}^{2\omega+1}\sum_{c=1}^{2\omega+1}\left(\frac{n}{\theta n-12\omega-7}\right)^{a+b+c}\ell^{2(a+b+c)}\\
&\leq&
\frac{\Delta^{6}\ell^{6(2\omega+1)}(2\omega+1)^3}{\theta n-12\omega-7}\\ 
&=& O\left(\frac{\Delta^{6}\ell^{12\omega+2}\omega^3}{\theta n}\right)
\end{eqnarray*}\\
Thus \(\Pr(\mu>0)=o(n^{-\e})\), for some constant \(\e>0\), since \(\D=O(n^{\kappa(d-1)/d})\) where \(\kappa<1/11\).

\textbf{(b)} The expected number of intersecting light cycles has the following upper bound:
\begin{equation}
\sum_{a=3}^{2\omega+1}\sum_{b=1}^{2\omega}\binom{n}{a}\binom{n}{b}\binom{a}{2}\frac{(a-1)!}{2}b!\ell^{2(a+b-2)}\Delta^{6}\frac{\mathcal{F}(\theta n - 2(a+b+1))}{\mathcal{F}(\theta n)} \label{NoInterCyc}
\end{equation}
\\
\textbf{Explanation} We have a cycle of size $a$ and another cycle which intersect with it. We treat the non-intersecting part of the second cycle as a chord of length $b$ on the $a$ cycle. Since at least one vertex is shared between the cycles, we only have $2\omega$ as the upper limit of $b$ in the inner sum. We need to pick one or two vertices of the $a$ cycle to connect the chord to, thus the factor of $\binom{a}{2}$. The factor of $b!$ counts (an upper bound of) the number of permutations of the vertices on the chord and the powers of $\ell$ and $\Delta$ play a similar role as per part (a) of the lemma. $\mathcal{F}(\theta n - 2(a+b+1))$ counts the number of configurations of the rest of the ports because $2(a+b+1)$ ports are occupied for a cycle and chord of size $a$ and $b$ respectively.

Thus, (\ref{NoInterCyc}) is bounded by 
\begin{eqnarray*}
&&\frac{\Delta^6}{\ell^6}\sum_{a=3}^{2\omega+1}\sum_{b=1}^{2\omega}\frac{n!}{(n-a)!}\frac{n!}{(n-b)!}(a-1)\ell^{2(a+b+1)}\frac{\mathcal{F}(\theta n - 2(a+b+1))}{\mathcal{F}(\theta n)}\\
&\leq& 
\frac{\Delta^6}{\ell^6}\frac{2\omega \ell^2}{\theta n-8\omega-3}\sum_{a=3}^{2\omega+1}\sum_{b=1}^{2\omega}\frac{n!}{(n-a)!}\frac{n!}{(n-b)!}\left(\frac{\ell^2}{\theta n-8\omega-3}\right)^{a+b}\\
&\leq&
\frac{\Delta^6}{\ell^4}\frac{2\omega}{\theta n-8\omega-3}\sum_{a=3}^{2\omega+1}\sum_{b=1}^{2\omega}\left(\frac{n}{\theta n-8\omega-3}\right)^{a+b}\ell^{2(a+b)}\\
&\leq&
\frac{\Delta^6\ell^{8\omega}(2\omega+1)^3}{\theta n-8\omega-3} = O\left(\frac{\Delta^6\ell^{8\omega}\omega^3}{\theta n}\right) = o(1)
\end{eqnarray*}
Thus as in part \textbf{(a)}, we have \(\Pr(\mu>0)=o(n^{-\e})\), for some constant \(\e>0\).
\proofend
\\

\begin{corollary}
\label{GamStruct}
With high probability, $\Gamma_v$ is either a tree or contains a unique cycle.
\end{corollary}
\proofstart
Clearly, if $v$ is locally tree-like, then $\Gamma_v$ will be a tree. If $G_v$ consists entirely of light vertices, then it can't have more than one cycle, since all cycles would be light and any two cycles would either intersect or be connected by a light path. By the lemma, \whp \, no such pair of light cycles exist, and so $\Gamma_v$ will contain a unique cycle. If $G_v$ contains a heavy vertex $h$ that has two light paths to the root $v$, then $h$ is part of a light cycle. No other heavy vertex $h'$ can have more than one light path to $v$ as $h'$ will also be part of a light cycle, and by the lemma, \whp \,  no such pair of light cycles exists. Therefore, any light cycle in $G_v$ will be unique (any other cycle will be heavy). Furthermore, there will be no more than two light paths from $h$ because of $h$ has at least three light paths to $v$, those three light paths would form an intersecting pair of cycles, which by the lemma, is forbidden \whp.\\ 
We see, therefore, that every heavy vertex in $G_v$ has one neighbour on a light path to $v$, with the exception of at most one heavy vertex, which can have two neighbours with such a property.  Since the pruning process deletes edges to neighbours of a heavy vertex that do not have a light path to the root, a pruned (heavy) vertex will have at most degree 1 in $\Gamma_v$, with at most one vertex ($h$), having degree 2 (and being part of a unique light cycle).\\
\proofend

\begin{lemma}\label{nondltl} With high probability,
\begin{description}
\item[(a)] The number of vertices \(v \in V\) that are not \(d\)-compliant is  at most \(n^{4\kappa(d-1)/d}\).
\item[(b)] There is no {\em small} vertex \(v\), \( \d \le d(v)< d\) which is not \(d\)-compliant.
\end{description}
\end{lemma}
\proofstart
\textbf{(a)} We lower bound the probability \(P\) that \(v\) is \(d\)-compliant
by the success, in the configuration model,
of the following  random process.
\\
{\bf Process \(\cP\):} 
\begin{enumerate}
	\item Each stub (half-edge) of $v$, is paired up with a stub of a vertex that has not been chosen thus far in the process. 
	\item For \(1 \leq i \leq \omega-1\),  and for each vertex \(w\) at level \(i\), the first \(d-1\) unpaired stubs of \(w\)
  			pair with stubs of vertices \(u\) of degree \(d(u) \ge d\) that have not been chosen thus far in the process.
\end{enumerate}

The tree created by process \(\cP\) involves \(N_1-1 = d_v\sum_{i=0}^{i=\omega-1}(d-1)^i \leq
\D (d-1)^{\omega}\) pairings. Let \(\sigma\) represent the sum
of degrees of vertices of degree less than \(d\). Thus

\begin{eqnarray*}
P &\geq& \prod_{i=1}^{N_1}\frac{\theta n-i\Delta-\sigma}{\theta n - 2i+1}\\ 
&\geq& \left(\frac{\theta n-N_1\Delta-\sigma}{\theta n - 2N_1+1}\right)^{N_1}\\ 
&\geq& \left(1-\frac{N_1\Delta+\sigma}{\theta n}\right)^{N_1}
\end{eqnarray*}
Let \(X\) count the number of vertices \(v\) 
that are
not \(d\)-compliant. Using the inequality \(1-(1-x)^y \leq xy\) for real \(x,y\), \(0 \leq
x \leq 1,y\geq 1\), we have
\begin{equation*}
\E[X] \leq n(1-P) = n\left(1-\left(1-
\frac{N_1\Delta+\sigma}{\theta
n}
\right)^{N_1}\right) \le \frac{N_1(N_1\Delta+\sigma)}{\theta}.
\end{equation*}
Now
\begin{equation}
\sum_{i=0}^{d-1}ic^i = \frac{(d-1)c^{d+1} - dc^d+c}{(c-1)^2} \leq \frac{dc^d}{c-1}+\frac{c}{(c-1)^2}\label{wqdf}
\end{equation}
for real $c \neq 1$. Hence under the assumption \(n_i=O(n^{\kappa i/d})\), we find that \eqref{wqdf} gives $\sigma = O(n^{\kappa (d-1)/d})$.

Given that \(\D=O(n^{\kappa(d-1)/d})\), we therefore have
\[
E[X] = \tilde O (\D^3+\D\s)=\tilde O(n^{3\kappa(d-1)/d})\leq c_1(\log^{c_2}n) n^{3\kappa(d-1)/d}
\]
for some \(c_1,c_2>0\). Then,
\begin{equation*}
\Pr(X \geq  n^{4\kappa(d-1)/d})=\tilde O(n^{-\kappa(d-1)/d}).
\end{equation*}

\textbf{(b)} In this case we have that the number of small vertices is \(O(n^{\kappa(d-1)/d})\)
and so \(\E[X]=\tilde O (n^{2\kappa(d-1)/d}/n)\). Hence, by Markov's inequality, the proposition holds with probability at least $1-n^{\Omega(1)}$.
\proofend
\\

\begin{lemma}\label{dltlLB} With high probability, there are \(n^{1-o(1)}\)
\(d\)-regular vertices \(v \in V\) with \(d(v)=d\).
\end{lemma}
\proofstart
We consider \(d\)-regular
vertices that have a root vertex \(v\) of degree \(d\).
Recall that \(n_d=|V_d|=\a n+o(n)\) for some constant \(\a>0\).
Let  \(N_2 =1+ d(d-1)^{\omega}\).
A \(d\)-regular tree of depth \(\om\) contains \(N_2\) vertices.
We proceed in a similar manner to
Lemma \ref{nondltl}, and bound the probability  \(P\) that a vertex \(v\)
 is \(d\)-regular by bounding the probability of success of the
construction of a \(d\)-regular tree in the configuration
model.
\begin{equation*}
P = \Pr(\text{a vertex }v\text{ is }d\text{-regular}) =
\prod_{i=1}^{N_2-1}\frac{d(n_d-i)}{\theta n - 2i+1} \geq
\left(d\frac{n_d-N_2}{\theta n}\right)^{N_2}
\end{equation*}
Let  \(M\) count the number of \(d\)-regular vertices,
then \(\E[M] = \mu = n_d P\), and
\begin{equation} \label{big}
\mu=\E[M] \ge   n^{1-o(1)}.
\end{equation}
To estimate \(\mathbf{Var}[M]\),
let \(I_v\) be the indicator that vertex \(v\) is \(d\)-regular. We have
\begin{eqnarray}
\E[M^2] = \mu + \sum_{v\in V_d}\sum_{w\in V_d, w \neq v}\E[I_vI_w],
 \label{jgnkjdn}
\end{eqnarray}
and
\[
\E[I_vI_w]=
\Pr(v,w\text{ are }d\text{-regular}, G_v \cap G_w = \emptyset)
+\Pr(v,w\text{ are }d\text{-regular}, G_v \cap G_w \neq
\emptyset).
\]
Now
\begin{equation}
\Pr(v,w\text{ are }d\text{-regular}, G_v \cap G_w =
\emptyset) = \prod_{i=1}^{2N_2-2}\frac{d(n_d-i-1)}{\theta n - 2i+1}
\leq {P}^2 \label{gfsfq}.
\end{equation}
For any vertex \(v\), the number of vertices \(u\) such that \(G_v
\cap G_w \neq \emptyset\) is bounded from above by \(N_2 +
dN_2^2\). Using this and (\ref{gfsfq}), we can bound
(\ref{jgnkjdn}) from above by \(\mu + \mu^2 + \mu(N_2+ dN_2^2)\).

By the Chebychev Inequality, for some constant \(0<\tilde{\epsilon}<1\),
\begin{equation*}
\Pr\left(|M-\mu| > \mu^{\frac{1}{2} + \tilde{\epsilon}}\right) \leq
\frac{\mathbf{Var}[M]}{\mu^{1 + 2\tilde{\epsilon}}}
 = \frac{\E[M^2]- \E[M]^2}{\mu^{1 + 2\tilde{\epsilon}}}\\ \leq \frac{\mu  + \mu N_2+\mu
dN_2^2}{\mu^{1 + 2\tilde{\epsilon}}} =O(n^{-\e})
\end{equation*}
for some constant $\e$.
 
The lemma now follows from \eqref{big}.
\proofend
\\

\begin{definition}\label{gooddef}
We call a nice graph $G$ \emph{good} if $\Phi(G)>1/100$, and if the statements of lemmas \ref{lightCPC}, \ref{nondltl} and \ref{dltlLB} hold for $G$. 
\end{definition}

The following is an immediate consequence of Definition \ref{gooddef} and the previously proved theorems and lemmas:
\\
\begin{proposition}
A graph $G\in \mathcal{G}(\mathbf{d})$ is good \whp.
\end{proposition}

As we shall see, $G$ being good allows us to use Lemma \ref{L3}, and, in particular, Corollary \ref{geom}. That is, when $G$ is good, the conditions of the lemma hold (to be demonstrated in section \ref{lemtrue}), and we have the tools to calculate the relevant quantities, in particular, (a bound on) the mixing time $T$ (see section \ref{MixingTimeSection}), and $R_v$, the expected number of returns in the mixing time to a vertex $v$. The latter quantity is calculated in the next section.

\section{Expected Number of Returns in the Mixing Time}
\label{Rvconds}
In this section we calculate bounds for the expected number of returns $R_v$ to a vertex $v$ within the mixing time $T$. First we require the following lemma.
\\

\begin{lemma}\label{R*Lemma}
Assume $G$ is good. Let $\cW_v^*$ denote the walk on $\Gamma_v$ starting at $v$ with $\Gamma_v^{\circ}$ made into absorbing states. Let $R_v^* = \sum_{t=0}^{\infty} r^*_t$ where $r^*_t$ is the probability that $\cW_v^*$ is at vertex $v$ at time $t$. There exists a constant $\zeta \in (0,1)$ such that
\begin{equation*}
R_v = R_v^* + O(\zeta^{\omega}).
\end{equation*}
\end{lemma}
\proofstart
We bound $|R_v-R_v^*|$ by using
\begin{equation}\label{A2}
R_v - R_v^* =\left(\sum_{t=0}^{\om} r_t-r^*_t\right)+ \left(\sum_{t=\om+1}^{T} r_t-r^*_t\right) -
\sum_{t=T+1}^{\infty}r^*_t .
\end{equation}
\underline{{\em Case \(t \le \om\).}}
When a particle starting from \(v\) lands at some vertex \(u \in \Gamma_v^{\circ}\),  this is either at
at distance \(\om\), or $u$ is
a heavy vertex at distance less than \(\om\) from \(v\).
In the case of \(u\) being heavy,
 by the light cycle condition,
 there are  at most two light paths back to \(v\) from \(u\) of length
 at most \(\om\). All other paths of length at most \(\om\) go via
 other heavy vertices.
 Hence if a particle is at \(u\), with probability at most \(2/\ell\) it will enter a light path
 to \(v\).
 Thus the probability of reaching \(v\) in time \(\om\) after having
 landed on a heavy vertex of \(\Gamma_v^{\circ}\) is at most \(O(\om/\ell)=o(\zeta^{\omega})\).

In the alternative  case that $u$ is at distance  \(\om\) from \(v\),
then  for \(t<\om\), \(r^*_t=r_t\). Thus we can write
\begin{equation}\label{wot}
 \left(\sum_{t=0}^{\om} r_t-r^*_t\right)=o(\zeta^{\omega}).
\end{equation}
\underline{{\em Case \(\om+1 \le  t \le T\).}}
Using (\ref{mix}) with \(x=u=v\) and \(\zeta = (1-\Phi^2/2) < 1\), we have  for
\(t \ge \om\), that \(r_t=\pi_v+O(\zeta^{t})\).
Since \(\Delta = O(n^{a}), a<1\), we have \(T \pi_v=o(\z^{\om})\) and so
\begin{equation}\label{wot2}
\sum_{t=\omega+1}^{T} |r_t-r^*_t| = \sum_{t=\omega+1}^{T} r_t
\leq \sum_{t=\omega+1}^{T} (\pi_v + \zeta^t )=O(\zeta^{\omega}).
\end{equation}

\underline{{\em Case \(t \ge T+1\).}}
It remains to estimate \(\sum_{t=T+1}^{\infty}r^*_t\).
Let $\s_t$ be the probability that the walk $\cW_v^*$ on $\G_v$ has not been absorbed by step $t$.
Then $r^*_t \leq \s_t$, and so
\[
\sum_{t=T+1}^{\infty}r^*_t \leq \sum_{t=T+1}^\infty \s_t,
\]

Assume first that $\Gv$ is a tree.
We estimate an upper bound for $\s_t$ as follows:
Consider an unbiased random walk $X_0^{(b)},X_1^{(b)},\ldots$
starting at $|b|< a\leq \omega$ on the finite line $(-a,-a+1,...,0,1,...,a)$,  with
absorbing states  $-a,a$. $X_m^{(0)}$ is the sum of $m$ independent $\pm 1$ random variables.
So the central limit theorem implies that there exists a constant $c>0$
such that
$$\Pr(X_{ca^2}^{(0)} \geq a \text{ or } X_{ca^2}^{(0)}\leq -a)
\geq 1-e^{-1/2}.$$

Now for any $t$ and $b$ with $|b|<a$, we have 
\begin{equation*}
\Pr(|X_{t}^{(b)}|< a) \leq \Pr(|X_{t}^{(0)}|< a)
\end{equation*}
which is justified with the following game: We have two walks, $A$ and $B$ coupled to each other, with $A$ starting at position $0$ and $B$ at position $b$, which, w.l.o.g, we shall assume is positive. The walk is a simple random walk which comes to a halt when either of the walks hits an absorbing state (that being, $-a$ or $a$). Since they are coupled, $B$ will win iff they drift $(a-b)$ to the right from $0$ and $A$ will win iff they drift $-a$ to the left from $0$. Therefore, given the symmetry of the walk, $B$ has a higher chance of winning.\\
Thus
\begin{equation*}
\Pr(|X_{2t}^{(0)}|< a) \leq \Pr(|X_{t}^{(0)}|< a)^2,
\end{equation*} 
since after $t$ steps, the worst case position for the walk to be at is the origin, $0$. Consequently, for any $b$ with $|b|<a$,
\begin{equation}\label{central}
\Pr(|X_{2ca^2}^{(b)}|\geq a)\geq 1-e^{-1}.
\end{equation}
Hence, for $t>0$,
\begin{equation}\label{at}
\s_t = \Pr(|X_{\t}^{(0)}|< a,\,\t=0,1,\ldots,t)\leq e^{-\rdown{t/(2ca^2)}}.
\end{equation}

Thus
\begin{equation}
\sum_{t=T+1}^\infty \sigma_t \leq \sum_{t=T+1}^\infty e^{-t/(3c\omega^2)} \leq \frac{e^{-T/(3c\omega^2)}}{1-e^{-1/(3c\omega^2)}} = O(\omega^2/n^{\Theta(\omega)}) = O(\zeta^{\omega})\label{65g46sg4s6dg}
\end{equation}

We now turn to the case where  $\G_v$ contains a unique cycle $C$.
The place where we have used the fact that $\G_v$ is a tree is in
\eqref{at} which relies on \eqref{central}.
Let $x$ be the furthest vertex of $C$ from $v$ in $\G_v$. This is the
only possible place where the  random walk is more likely
to get closer to $v$ at the next step. We can see this by
considering the breadth first construction of $\G_v$. Thus we can compare our walk with
 random walk on $[-a,a]$ where there is a unique value $x<a$ such that only
at $\pm x$ is the walk more likely to move towards the origin and
even then this probability is at most 2/3.

Since we have
$$\Pr(\exists \t\leq ca^2:\;|X_{\t}^{(b)}|= x)\geq 1-e^{-1/2}.$$
We now consider the probability the particle walks from e.g. $x$ to $a$ without returning to the cycle: There is $1/3$ chance of moving away from the origin (towards $a$). Now from standard results, an unbiased random walk in one dimension between $x$ and $a$ starting at position $x+1$ would have probability $1/(a-x)$ of hitting $a$ before hitting $x$. Since the actual walk will be biased in favour of movement toward $a$, a probability of $1/3(a-x)$ serves as a lower bound on the probability of going from $x$ to $a$ without returning to $x$.\\
Thus
$$\Pr(\exists \t\leq ca^2:\;|X_{\t+a-x}^{(b)}|= a)\geq
(1-e^{-1/2})/3 a,$$
and
\begin{equation}
\s_t= \Pr(|X_{\t}^{(0)}|< a,\,\t=0,1,\ldots,t)\leq
(1-(1-e^{-1/2})/3a)^{\rdown{t/(2ca^2)}}\leq e^{-t/(20c\omega)^3}. \label{att6465}
\end{equation}

Thus
\begin{equation*}
\sum_{t=T+1}^\infty \sigma_t \leq \sum_{t=T+1}^\infty e^{-t/(20c\omega)^3} \leq \frac{e^{-T/(20c\omega)^3}}{1-e^{-1/(20c\omega)^3}} = O(\omega^3/n^{\Theta(1)}) = O(\zeta^{\omega})
\end{equation*}
\proofend

\begin{lemma}\label{RvRltl}
Assume $G$ is good. For a vertex $v \in V$,
\begin{description}
\item[(a)] if $v$ is $d$-regular, then $R_v = \frac{d-1}{d-2} + O(\zeta^{\omega})$, 
\item[(b)] if $v$ is $d$-tree-compliant then $R_v \leq \frac{d-1}{d-2} + O(\zeta^{\omega})$,
\item[(c)] if $d_v \leq \ell$, $R_v \leq \frac{\delta^2 + \delta - 4}{\delta(\delta-2)} = \mathcal{R} \leq 8/3$
\end{description}
for some constant $\zeta \in (0,1)$.
\end{lemma}
\proofstart
We use Lemma \ref{R*Lemma}.

\textbf{(a)} We calculate $R_v^*$ for a walk $\cW_v^*$ on an $r$-regular tree $\Gamma_v$ with $\Gamma_v^{\circ}$ made into absorbing states. For a biased random walk on $(0,1,...,k)$, starting at vertex 1,  with  absorbing
states $0,k$, and with transition probabilities at vertices $(1,\ldots,k-1)$
of $q=\Pr($move left$)$,
$p=\Pr($move right$)$; then
\begin{equation}\label{abso}
\Pr(\text{absorption at } k)=\frac{(q/p)-1}{(q/p)^{k}-1}.
\end{equation}
This is the escape probability $\rho$ - the probability that after the particle moves from $v$ to an adjacent vertex, it reaches an absorbing state without having visited $v$ again. $R_v^* = 1+1/\rho-1 = 1/\rho$.\\
We project $\cW_v^*$ onto $(0,1,\ldots,\omega)$ with $p=\frac{d-1}{d}$ and $q = \frac{1}{d}$ giving 
\begin{equation}
R_v^* = \left(1-\frac{1}{(d-1)^{\omega}}\right)\frac{d-1}{d-2} = \frac{d-1}{d-2} - O\left((d-1)^{-\omega}\right)
\end{equation}
and part \textbf{(a)} of the lemma follows.\\

\textbf{(b)} If \(v\) is \(d\)-tree-compliant, we can prune \(G_v\) removing
edges from each vertex (other than \(v\)) until \(v\) is \(d\)-regular. Treating
the edges as having unit resistance, this pruning process cannot
decrease the effective resistance between \(v\) and a hypothetical
vertex \(\z\) that is connected by a zero-resistance edge to each of the
vertices in \(\Gamma_v^{\circ}\) (and no others). Then by part
\textbf{(a)} and Rayleigh's monotonicity law part \textbf{(b)} of
the lemma follows. (Here we are using the the fact that the probability
of reaching \(\z\) before returning to \(v\) is equal to \(\frac{1}{d(v)R}\) where \(R\) is the
effective resistance between \(v\) and \(\z\). Rayleigh's Law states that deleting edges increases \(R\)).

\textbf{(c)} Using \eqref{mix} and the fact that $T\pi_v = o(1)$, we have $\sum_{t=0}^T r_t \leq \sum_{t=0}^T  (\p_v+O(\z^t)) = O(1)$, quickly demonstrating that $R_v$ is bounded from above by a constant. To get the tighter bound, we observe that for any locally tree-like vertex $v$, $R_v \leq \frac{\delta-1}{\delta-2} + O(\zeta^{\omega})$ by part \textbf{(a)}. If there is a cycle (which will be unique) $\Gamma_v$ will be tree-like up to some distance $L-1<\omega$ and non-tree-like at level (distance) $L$, then there will be one vertex $v_c$ on the cycle at level $L$ if the cycle is an even size, and two vertices $v_c, w_c$ if it is an odd size. The subtree up to level $L-1$ will be $\delta$-tree-compliant, and so if the particle does reach $u_c \in \{v_c,w_c\}$, it will halt if $u_c$ is an absorbing vertex and move one level down (i.e., further away from the root) with probability at least $(\delta-2)/\delta$ if $u_c$ is not an absorbing state. If it is an even cycle, it will move up the tree (i.e., closer to the root) back to level $L-1$ with probability at most $2/\delta$, and if it is an odd cycle, it will move up with probability at most $1/\delta$ and stay at the same level (moving to the other vertex in $\{v_c, w_c\}$) with probability at most $1/\delta$. Levels below $L$ (should they exist) continue to be tree like and so maintain the probability distribution of going down with probability at least $(\delta-1)/\delta$ and going up with probability at most $1/\delta$.\\

Given that there exists a unique cycle in $\Gamma_v$, $R_v^*$ can be no larger than that for the case where $\Gamma_v$ is such that all non-leaf vertices have degree $\delta$ and the absorbing states $\Gamma_v^{\circ}$ are all at distance $\omega$ from the root. The only other variable is the location of the cycle.\\

We shall resort again to electrical network theory to bound $R_v^*$. We assign a unit flow to $\Gamma_v$ with the root vertex $v$ as the source and a hypothetical sink vertex that is connected by a zero-resistance edge to each of the vertices in $\Gamma_v^{\circ}$ (and no others). We treat each edge in $\Gamma_v$ as having unit resistance. When $\Gamma_v$ is a $r$-regular tree (apart from leaves), we split the flow coming into a vertex evenly across the $r-1$ neighbours down the next level. Denote this energy by $P_T(r)$. Note $P_T(r) > P_T(s)$ for $r<s$. To calculate bounds for the minimal energy when $\Gamma_v$ contains a cycle, we assign the same flow as for a tree, except at the cycle. Consider first an even-sized cycle. There will be two sources of flow into $v_c$ (defined above), and by symmetry, these will be equal - call them $f$. If $\delta = 3$, these flows combine  at $v_c$ to become a $2f$ flow going through one edge, then split off into two flows size $f$ at the next vertex. The flow assignment at subsequent vertices is as per usual, with an even split at each vertex. If $\delta > 3$, then from $v_c$, we assign a flow of $f$ to each of two edges into the next level down and then continue as normal. For $\delta=3$ the energy $P_{\Gamma_v}$ of the system is bound by $P_T(\delta) + (2f)^2 + 2f^2 = P_T(\delta) + 6f^2$ and for $\delta>3$ the energy is bounded by $P_T(\delta) + 2f^2$ and since $f \leq \frac{1}{(\delta-1)\delta}$, we have a bound of $P_T(\delta) + \frac{6}{(\delta-1)^2\delta^2}$.\\ 

For the case where the cycle is odd size, we assign the flows on the same basis as a tree (as in the case of the even cycle) except that the flow $f$ which comes into each of vertices $\{v_c, w_c\}$ will be directed from each vertex through a single edge going to the next level down. Thus for all $ \delta \geq 3$, $P_{\Gamma_v} \leq P_T(\delta) + 2f^2$, and since $f \leq 1/\delta$ we have $P_{\Gamma_v} \leq P_T(\delta) + 2/\delta^2 \ignore{\leq P_T(3) + 2/9}$. We bound $P_T(\delta)$ by calculating the effective resistance of an infinite $\delta$-regular tree: $P_T(\delta) \leq \sum_{i=0}^{\infty}\frac{1}{\delta(\delta-1)^i} = \frac{1}{\delta}\frac{\delta-1}{\delta-2}$. Thus $P_{\Gamma_v} \leq \frac{1}{\delta}\frac{\delta-1}{\delta-2} + \frac{2}{\delta^2}$. This is an upper bound on the effective resistance between the root vertex and a hypothetical vertex connected to each vertex in $\Gamma_v^{\circ}$. The escape probability $\rho$ of a walk starting at the root is related as $\rho = \text{Res}(v)/P_{\Gamma_v}$ where $\text{Res}(v)$ is the parallel resistance of the edges incident with $v$ (or alternatively, the reciprocal of the sum of conductances of the edges incident with $v$). This is $1/\delta$ when $d_v = \delta$. Thus $R^*_v \leq \frac{\delta^2 + \delta - 4}{\delta(\delta-2)}$ which is monotonically decreasing and for $\delta=3$ is $8/3$.

\proofend

\section{Lemma Conditions}
\label{lemtrue}
We address the conditions of Lemma \ref{L3}. For condition \textbf{(b)}, we note that $n < m < n^2$ implies the condition $T\pi_v=\Omega(n^{-2})$ holds since $T$ is larger than constant and $T^2\pi_v=o(1)$ holds since  $T^2=(\omega^3\log n)^2=o(n^c)$ for any positive constant $c$, and $\pi_v=O(n^{-c})$ for some positive constant $c$. It remains to show condition \textbf{(a)} holds:
\begin{lemma}\label{RTRs} 
Assume $G$ is good. For $|z|\leq  1+\l$, 
there exists a constant $\psi>0$ such that $|R_T(z)|\ge \psi$.
\end{lemma}
\proofstart
As before, let  $\Gvc$ be  the set of absorbing states of $\Gv$.
We consider a walk in $\G_v$ that starts at $v$, and each time it leaves $v$ it is terminated if it doesn't return to it within time $T$. It is also terminated if it reaches an absorbing state (states in $\Gvc$) at any time. Let
$\b(z)=\sum_{t=1}^T \b_t z^t$ where $\b_t$ is the
probability of a first return to $v$  at time $t\le T$.
Let $\a(z)=1/(1-\b(z))$, and write
$\a(z)=\sum_{t=0}^{\infty} \a_t z^t$,
so that  $\a_t$ is the probability that the walk is at $v$ at time $t$. Note that $\alpha_t$ is not quite the same as $r^*_t$ in Lemma \ref{RvRltl}, which would be generated as $1/(1-\sum_{t=1}^{\infty}b^*_t)$, where $b^*_t$ represents the probability of a first return of $\cW_v^*$ at time $t$. Observe that $\alpha_t \leq r^*_t \leq r_t$.\\ 
We shall prove below that the radius of convergence of $\a(z)$ is at least
$1+\OM(1/\omega^2)$.

We can write
\begin{eqnarray}
R_T(z)&=&\a(z)+Q(z)\nonumber\\
&=&\frac{1}{1-\b(z)}+Q(z),\label{Bs}
\end{eqnarray}
where $Q(z)=Q_1(z)+Q_2(z)$, and
\begin{eqnarray*}
Q_1(z)&=&\sum_{t=0}^T(r_t-\a_t)z^t\\
Q_2(z)&=&-\sum_{t=T+1}^\infty \a_t z^t.
\end{eqnarray*}
We note that $Q(0)=0$, $\a(0)=1$ and $\b(0)=0$.

We claim that the expression \eqref{Bs} is well defined for
$|z|\le 1+\l$.
We will show  below that
\begin{equation}\label{Q2}
|Q_2(z)|=o(1)
\end{equation}
for $|z|\leq 1+2\l$ and thus the radius of convergence of $Q_2(z)$
(and hence $\a(z)$) is
greater than $1+\l$.
This will imply that $|\b(z)|<1$ for $|z|\leq
1+\l$. For suppose there exists $z_0$ such that $|\b(z_0)|\geq 1$.
Then $\b(|z_0|)\geq |\b(z_0)|\geq 1$ and we can assume (by scaling) that $\b(|z_0|)=1$. We have
$\b(0)<1$ and so we can assume that $\b(|z|)<1$ for $0\leq
|z|<|z_0|$. But as $\r$ approaches 1 from below, \eqref{Bs} is valid
for $z=\r|z_0|$ and then $|R_T(\r|z_0|)|\to \infty$, contradiction.

Recall that $\l=1/K T$.
Clearly $\b(1)\le 1$ and so  for $|z|\leq 1+\l$
$$\b(|z|)\leq\b(1+\l)\leq\b(1)(1+\l)^T\leq e^{1/K}.$$
Using $|1/(1-\b(z))| \ge 1/(1+\b(|z|))$ we obtain
\begin{equation}\label{Bs<1}
|R_T(z)|\geq \frac{1}{1+\b(|z|)}-|Q(z)|
\geq \frac{1}{1+e^{1/K
}}-|Q(z)|.
\end{equation}
We now prove that $|Q(z)|=o(1)$ for $|z|\le 1+\l$ and the lemma will follow.

Turning our attention first to $Q_1(z)$, the following results hold both for $\Gv$ a tree, and $\Gv$ containing a cycle.
\begin{equation}\label{Q1}|Q_1(z)| \leq(1+\l)^T |Q_1(1)| \leq e^{2/K}\sum_{t=0}^{T} |r_t-\a_t|\end{equation}
Now using part of the proof of lemma \ref{RvRltl}, we see that $\sum_{t=0}^{T} |r_t-\a_t| = o(1)$, hence $|Q_1(z)| = o(1)$.\\

We now turn our attention to $Q_2(z)$. As in lemma \ref{RvRltl}, let $\sigma_t$ be the probability that a walk on $\G_v$ (with $\Gvc$ absorbing states) starting at $v$ has not been absorbed at step $t$. Then $\sigma_t \geq \alpha_t$, so
\begin{equation*}
|Q_2(z)| \leq \sum_{t=T+1}^\infty \s_t |z|^t,
\end{equation*}
If $G_v$ is a tree we can use \eqref{at} to see that the radius of convergence of $Q_2(z)$ is at least $e^{1/(3ca^2)}$. As $a\le \omega$, $e^{1/(3ca^2)} \gg
1+2\l$ and for $|z|\leq 1+2\l$,
$$|Q_2(z)|\leq \sum_{t=T+1}^\infty e^{2\l t-\rdown{t/(2ca^2)}}=o(1).$$

This lower bounds the radius of convergence of $\a(z)$,
proves \eqref{Q2} and then \eqref{Q1}, \eqref{Q2} and
\eqref{Bs<1} complete the proof of the case where $\G_v$
is a tree.\\

In the case that $G_v$ contains a unique cycle, we can use \eqref{att6465} to see that the radius of convergence of $Q_2(z)$ is  at least
$1+\frac{1}{25c\omega^3}>1+2\l$, assuming that $K$ (defined in
\eqref{lamby}) is sufficiently
large. Finally, if $z\in C_\l$ then
$$|Q_2(z)|\leq \sum_{t=T+1}^\infty e^{(\l -1/(20c\omega^3))t}\leq
\frac{e^{-T/(25c\omega^3)}}{1-e^{-1/(25c\omega^3)}}=O\left(w^3/n^{\Theta(1)}\right)$$. \proofend\\

\section{Cover Times}
We now prove the main theorem of the chapter, Theorem \ref{CovTh}, which we restate here,

\textbf{Theorem \ref{CovTh} (restated)}\emph{
Let  \(G\) be chosen \uar\ from \(\mathcal{G}(\textbf{d})\), where \textbf{d} is nice.
Then \whp,
\begin{equation}
\cov[G] \sim \frac{d-1}{d-2}\frac{\theta}{d}\;n\log n, 
\end{equation}
where $d$ is the effective minimum degree and $\theta$ is the average degree. 
}

In the following, we assume $G$ is good.
\subsection{Upper Bound on Cover Time } \label{ubsec}
As per the definition in section \ref{Random Walks on Graphs}, $c_u$ is the time taken by the random walk \(\cW_u\) starting from vertex $u$ to visit
every vertex of a connected graph \(G\). Let \(U_t\) be the number of
vertices of \(G\) which have not been visited by \(\cW_u\) by step \(t\).
We note the following:
\begin{eqnarray}
\cov_u[G]=\E[c_u]&=& \sum_{t > 0} \Pr(c_u \ge t), \label{ETG} \\
\label{TG} \Pr(c_u\geq t)=\Pr(c_u > t-1)&=&\Pr(U_{t-1}>0)\le
\min\{1,\E[U_{t-1}]\}.
\end{eqnarray}
Recall from \eqref{frat} that
\(\ul A_s(v)\) is the event that vertex \(v\) has not been visited by time \(s\).
It follows from (\ref{ETG}), (\ref{TG}) that
\begin{equation}\label{shed11}
\cov_u[G] \leq t+1+ \sum_{s \geq t} \E[U_{s}] = t+1+\sum_v \sum_{s \ge
t}\Pr(\ul A_s(v)).
\end{equation}
Let \(t_0 = \brac{\frac{d-1}{d-2}\frac{\theta}{d}} n\log n\) and \(t_1 =
\left(1+\e\right)t_0\), were \(\e=o(1)\) is sufficiently large that all
inequalities claimed below hold. We assume that the high probability claims of
Sections \ref{seca}, \ref{Rvconds} hold.

Recall from \eqref{pv} that \( p_v=(1+O(T\pi_v)) d_v/(\theta nR_v)\) and that Corollary \ref{geom} gives the probability that \(\cW_u\) has not visited
\(v\) during \([T,t]\) as
\begin{eqnarray*}
\Pr(\ul A_t(v))&=&\frac{(1+O(T\p_v))}{(1+p_v)^{t}} +O(T^2\p_v e^{-\lambda t/2})\\
&=&(1+o(1))(1+p_v)^{-t}
\end{eqnarray*}

Thus
\begin{eqnarray}
\sum_{t \geq t_1}(1+o(1))(1+p_v)^{-t} &=&(1+o(1))(1+p_v)^{-t_1}\sum_{t \geq t_1}(1+p_v)^{-(t-t_1)}\nonumber\\
&\leq& 2p_v^{-1}e^{-t_1p_v}\nonumber\\
&=& O(1)\frac{\theta
nR_v}{d_v}\exp\left(-(1+\Th(\e))\frac{d_v}{d}\frac{d-1}{d-2}
\frac{\log n}{R_v}\right)\label{A1}.
\end{eqnarray}

We consider the following subsets of \(V\):\\
(i) \(V_A=\bigcup_{\d \le i <d}V_i\).\\
(ii) \(V_B=\bigcup_{d \leq i}\{v\in V_i: v \text{ is } d\text{-compliant}\}\).\\
(iii)\(V_C= \bigcup_{d \leq i \leq \ell}\{v\in V_i: v \text{ is not } d\text{-compliant}\}\).\\
(iv) \(V_D=\bigcup_{\ell < i}V_i\).\\

{\em Case (i): \(\d \le d_v < d\).}\\
For these vertices, \(\G_v\) is \(d\)-compliant by Lemma  \ref{nondltl}.
Consider vertices in \(V_i\), \(i <d\). By Lemma
\ref{RvRltl} \textbf{(b)}, \(R_v \leq (1+o(1))\frac{d-1}{d-2}\) so for \(v \in V_i\)
\eqref{A1} is bounded by \(O(\theta n)
n^{-(1+o(1))\frac{i}{d}}\).
Recall that \(|V_i|=O(n^{\kappa i/d})\) where \(\kappa <1\). Thus
\begin{eqnarray*}
\sum_{v \in V_i}\sum_{t \geq t_1}(1+o(1))e^{-tp_v}
\leq O(\theta n)n^{\kappa i/d}n^{-(1+o(1))i/d}
&=& o(t_1).
\end{eqnarray*}

{\em Case (ii): \(d \le d_v \), \(v\) is \(d\)-compliant.}\\
For \(v \in V_B\)
\eqref{A1} is bounded by \(O(\theta )
n^{-\Theta(\e)}\). Therefore
\begin{eqnarray*}
\sum_{v \in V_B}\sum_{t \geq t_1}(1+o(1))e^{-tp_v} \leq
\sum_{v \in V_B}O(\theta )n^{-\Theta(\e)}
=O(\theta n) n^{-\Theta(\e)} = o(t_1).
\end{eqnarray*}

{\em Case (iii): \(d \le d_v \leq \ell\), \(v\) is  not \(d\)-compliant.}\\
For vertices \(v \in V_C\) \eqref{A1} is bounded by
\(O(\theta  n)
n^{-(1+\Theta(\e))\frac{1}{\mathcal{R}}\frac{d-1}{d-2}}\) where $\mathcal{R}\leq 8/3$.
By Lemma \ref{nondltl}, \(|V_C| \leq n^{4\kappa(d-1)/d}\), where $0<\kappa<1/11$
hence
\begin{eqnarray*}
\sum_{v \in V_C}\sum_{t \geq t_1}(1+o(1))
e^{-tp_v} &=& \sum_{v \in V_C}
O(\theta n)
n^{-(1+\Theta(\e))\frac{1}{\mathcal{R}}\frac{d-1}{d-2}}\\
&=&O(n^{4\kappa(d-1)/d}\theta n)
n^{-(1+\Theta(\e))\frac{1}{\mathcal{R}}\frac{d-1}{d-2}}\\
&=& o(t_1),\\
\end{eqnarray*}
since \[4\kappa \frac{(d-1)}{d} - (1+\Theta(\e))\frac{1}{\mathcal{R}}\frac{d-1}{d-2} < (d-1)\left(\frac{4\kappa}{d}-(1+\Theta(\e))\frac{3}{8(d-2)}\right)<-\frac{d-1}{88}.\]

\emph{Case (iv): $\ell < d_v$}\\
Clearly, $R_v \leq T = w^3\log n$ thus, for vertices \(v \in V_D\) \eqref{A1} is bounded by $n^{-\Theta(\log n/\omega^3)}$. Hence
\begin{eqnarray*}
\sum_{v \in V_D}\sum_{t \geq t_1}(1+o(1))e^{-tp_v} = \sum_{v \in V_D}n^{-\Theta(\log n/\omega^3)} \leq n^{-\Theta(\log n/\omega^3)}=o(t_1).
\end{eqnarray*}\\

In each of the cases above, the term \(\sum_v \sum_{s \ge
t}\Pr(\ul A_s(v))=o(t_1)\) and thus, from \eqref{shed11},
\(\cov_u[G] \le \ooi t_1\) as required. This completes the proof of the
upper bound on cover time of \(G\).
\proofend

\subsection{Lower Bound on Cover Time } \label{lbsec}
Let \(t_2=(1-\e)t_0\), were \(\e=o(1)\) is sufficiently large that all inequalities claimed
below hold.
For vertex \(u\) of degree \(d\), we exhibit a set of vertices \(S\) such that at
time \(t_2\) the probability the set \(S\) is
covered by the walk \(\cW_u\) tends to zero. Hence \(c_u > t_2\),
\whp\ which implies that \(\cov[G]\geq t_0-o(t_0)\).

We construct \(S\) as follows.
Let \(S_d\) be the set of \(d\)-regular vertices of degree \(d\).
Lemma \ref{dltlLB}
tells us that  \(|S_d|=n^{1-o(1)}\).
Let \(S\) be a maximal subset of \(S_d\)
such that the distance between any two
elements of \(S\) is least \(\om\).
Thus \(|S|=\Om(n^{1-o(1)}/d^{\omega})\).

Let \(S(t)\) denote the subset of \(S\) which has not been visited by
\(\cW_u\) after step \(t\).
Let \(v \in S\), then
\begin{eqnarray*}
\Pr(\ul A_{t_2}(v))&=&(1+o(1))(1+p_v)^{-t_2}+ o(n^{-2})\\
&=&(1+o(1))e^{-t_2 p_v(1-O(p_v))}+ o(n^{-2}).
\end{eqnarray*}
(The above is justified because $(1+x)^y=e^{y\log(1+x)}=e^{y(x-\frac{x^2}{2}+\frac{x^3}{3}-\ldots)}=e^{yx(1-O(x))}$ for real $x,y$ with $|x|<1$).

Hence,
\begin{eqnarray}
\E(|S(t_2)|) &\geq&(1+o(1))|S|e^{-(1-\e)t_0p_v}\nonumber\\
&=& \Omega \bfrac{n^{\e/2-o(1)}}{d^{\om}} \rightarrow \infty \label{beta}.
\end{eqnarray}
Let \(Y_{v,t}\) be the indicator for the event \(\ul A_t(v)\). Let
\(Z=\{v,w\} \subset S\). We will show (below) that that for \(v,w \in
S\)
\begin{equation}\label{poi}
\E(Y_{v,t_2}Y_{w,t_2})=\frac{1+O(T\p_v)}{(1+p_Z)^{t_2}}+o(n^{-2}),
\end{equation}
where
\(p_Z \sim p_v+p_w\). Thus
\begin{equation*}
\E(Y_{v,t_2}Y_{w,t_2})=(1+o(n^{-2}))\E(Y_{v,t_2})\E(Y_{w,t_2})
\end{equation*}
which implies
\begin{equation}\label{zzz}
\E(|S(t_2)|(|S(t_2)|-1))\sim \E(|S(t_2)|)(\E(|S(t_2)|)-1).
\end{equation}
It follows from (\ref{beta}) and (\ref{zzz}), that
\begin{equation*}
\Pr(S(t_2)\neq \emptyset)\geq \frac{\E(|S(t_2)|)^2}{\E(|S(t_2)|^2)}
=\frac{1}{\frac{{\bf E}(|S(t_2)|(|S(t_2)|-1))}{{\bf
E}(|S(t_2)|)^2}+\E(|S(t_2)|)^{-1}} =1-o(1).
\end{equation*}

{\bf Proof of (\ref{poi}).}
Let $\wh{G}$ be obtained from $G$ by merging $v,w$ into a single node $Z$. This node has
degree $2d$.

We apply Lemma \ref{L3} to $\wh{G}$. That $\p_Z=\frac{2d}{\theta n}$ is clear.
Furthermore, the vertex  $Z$ is tree-like
up to distance $\om$ in $\wh{G}$. The derivation of $R_Z$ as in Lemma \ref{RvRltl}(a) is valid.
The fact that the root vertex of the corresponding infinite
tree has degree $2d$ does not affect the calculation of $R_Z^*$.

Thus 
\begin{eqnarray*}
\wh{\Pr}(\ul A_{t_2}(Z))&=&\frac{1+O(T\p_Z)}{(1+p_Z)^{t_2}}+o(n^{-2})\\
&=& \frac{1+O(T\p_v)}{(1+p_v)^{t_2}}\frac{1+O(T\p_w)}{(1+p_w)^{t_2}} +o(n^{-2})\\
&=& \Pr(\ul A_{t_2}(v))\Pr(\ul A_{t_2}(w))+o(n^{-2})
\end{eqnarray*}

Now
\begin{eqnarray}
\wh{\Pr}(\ul A_{t_2}(Z)) &=& \sum_{x \neq Z}\wh{P}_u^{(T)}(x)\wh{\Pr}(\mathcal{W}_x(t-T) \neq Z, T \leq t \leq t_2)\nonumber\\
&=& \sum_{x \neq Z}\left(\frac{d_v}{\theta n}+O(n^{-3})\right)\wh{\Pr}(\mathcal{W}_x(t-T) \neq Z, T \leq t \leq t_2)\nonumber\\
&=& \sum_{x \neq v,w}\left(P_u^{(T)}+O(n^{-3})\right)\Pr(\mathcal{W}_x(t-T) \neq v,w, T \leq t \leq t_2)\label{nmpm}\\
&=&\Pr(\mathcal{W}_u(t) \neq v,w, T \leq t \leq t_2)+O(n^{-3})\nonumber\\
&=&\Pr(\ul A_{t_2}(v) \wedge \ul A_{t_2}(w)) +O(n^{-3})\nonumber
\end{eqnarray}
Equation \eqref{nmpm} follows because there is a natural measure-preserving
mapping from the set of walks in $G$ which start at $u$ and do not visit $v$ or $w$, to
the corresponding set of walks in $\wh{G}$ which do not visit $Z$.

We have shown that
\begin{equation*}
\Pr(\ul A_{t_2}(v) \wedge \ul A_{t_2}(w)) = (1+o(n^{-2}))\Pr(\ul A_{t_2}(v))\Pr(\ul A_{t_2}(w))+o(n^{-2})
\end{equation*}
\proofend

\chapter{Weighted Random Walks}\label{Weighted Random Walks}
Chapters \ref{ch:The Cover Time of Cartesian Product Graphs} and \ref{CovDS} studied random walks on different graph structures, but in both cases the walks were simple, meaning that each step an edge is chosen \uar\ from the edges incident on the current vertex, and the walk transitions that edge. However, we may also define random walks on graphs $G=(V,E,c)$ with differently-weighted edges, where an edge $e=(u,v)$ has some non-negative weight, or \emph{conductance}, $c(e)$, and the probability of transitioning $e$ from $u$ is given by $c(e)/c(u)$ where $c(u) = \sum_{e =(u,v)\in E}c(e)$. Simple random walks are a special case where $c$ is a constant function. Chapters \ref{intro}, \ref{Theory of Markov Chains and Random Walks} and \ref{networktheory} discuss the details, and we refer the reader there. We shall call random walks on weighted graphs \emph{weighted random walks}.

In this chapter, to avoid confusion between conductance $c(e)$ of an edge, and conductance $\Phi(\mathcal{M})$ of a Markov chain $\mathcal{M}$ as defined in Definition \ref{MarkovchainCondDef}, we instead refer to edge conductance as \emph{weight} and use $w$ rather than $c$. We will still use \emph{resistance} $r(e)$ and \emph{effective resistance} $R(u,v)$ as there is no ambiguity.

As with chapters \ref{ch:The Cover Time of Cartesian Product Graphs} and \ref{CovDS}, we deal only with simple, connected, undirected graphs.

\section{Weighted Random Walks: Hitting Time and Cover Time}
For simple random walks, asymptotically tight bounds for cover time were given by \cite{Feige1} and \cite{Feige2}:
\[
(1+o(1))n\log n \leq \cov[G] \leq (1+o(1))\frac{4}{27}n^3,
\]
and these lower and upper bounds are met by the complete graph and the lollipop graph respectively. See section \ref{FeigeSection}.

However, the hitting and cover times for weighted random walks have not been studied in as much depth, and it is not difficult to formulate many open questions on their behaviour. In particular, what bounds exist for hitting times and cover time? This was addressed in part by \cite{ikeda1} and \cite{ikeda2}. The investigation is framed as follows. For a graph $G$ let $\mathcal{P}(G)$ denote the set of all transition probability matrices for $G$, that is, stochastic matrices (i.e., non-negative entries with each row summing to $1$), that respect the graph structure. The latter means that $P_{u,v}=0$ if $(u,v)\notin E(G)$. For $P\in \mathcal{P}(G)$ let $H_G(P)$ denote the maximum hitting time in $G$ with transition matrix $P$, and $C_G(P)$ similarly for cover time. Let 
\[
H_G = \inf_{P\in \mathcal{P}(G)}H_G(P)\quad \text{ and }\quad C_G = \inf_{P\in \mathcal{P}(G)}C_G(P).
\]

Note that if for a graph $G$ one knows a spanning tree $T_G$, a transition matrix $P$ can be constructed that is a simple random walk on $T_G$, and ignores all other edges of $G$. As per section \ref{General Bounds and Methods}, this implies a $O(n^2)$ upper bound on $H_G$ and $C_G$.

It is asked
\begin{description}
	\item[(1)] What are general upper bounds on 
		\[
		H(n) = \max_{G : V(G)=n}H_G \quad \text{ and }\quad C(n) = \max_{G : V(G)=n}C_G.
		\]
	\item[(2)] What is the minimum local topological information on $G$ that is always sufficient to construct a transition matrix that achieves each of $H(n)$ and $C(n)$?
\end{description}

In \cite{ikeda2}, it is shown that for a path graph $P_n$, any transition matrix will have $\Omega(n^2)$ maximum hitting time (and therefore, cover time). This, in conjunction with the spanning tree argument, implies $\Theta(n^2)$ for both $H(n)$ and $C(n)$. 

The second question is addressed with a particular weighting scheme that gives $O(n^2)$ maximum hitting time for any graph. Note, Matthews' upper bound technique (see section \ref{Upperbound: Matthews' Technique}) applies to weighted walks as well, implying a $O(n^2\log n)$ bound on cover time. The weighting scheme, which we call the \emph{Ikeda} scheme is as follows:
\begin{equation}
 P_{u,v} = 
  \begin{cases}
   \frac{1/\sqrt{d(v)}}{\sum_{w\in N(u)}1/\sqrt{d(w)}} & \text{if } v \in N(u) \\
   0       & \text{otherwise} 
  \end{cases}\label{ikedascheme}
\end{equation}
where $d(v)$ is the degree of $v$ and $N(v)$ is the neighbour set of $v$. 

\textbf{A note on application}  In the algorithmic context of graph exploration, simple random walks are generally considered to have the benefit of not requiring information beyond what is needed to choose the next edge \uar. Generally, this implies that a token making the walk can be assumed to know the degree of the vertex it is currently on, but no more information about the structure of the graph. In the Ikeda scheme, information required in addition to the vertex degree, is the degrees of neighbouring vertices. Other local information - such as connectivity, is not required.

\section{Minimum Degree Weighting Scheme}
We present an alternative weighting to the Ikeda scheme: Simply, for $G=(V,E)$, assign each edge $(u,v)$ weight $w(u,v) = 1/\min\{d(u), d(v)\}$ (equivalently, each edge is assigned resistance $r(u,v)=\min\{d(u), d(v)\}$). This weighting scheme defines the transition matrix of a weighted random walk. We may, as a matter of convenience, say that $w(u,v)=0$ if $(u,v)\notin E$ in calculations of transition probabilities.  

We call this scheme the \emph{minimum degree} (or \emph{min-deg}) scheme. It places no greater information requirement in applications than the Ikeda scheme; it also only uses limited, local knowledge. In fact, it is easily checked that equation \eqref{ikedascheme} is equivalent to weighting the edges $w(u,v) = 1/\sqrt{d(u)d(v)}$.

We have the following (in)equalities.
\begin{equation}
w(u,v) \leq \frac{1}{d(u)} + \frac{1}{d(v)} \leq 2w(u,v).\label{mindegeq1}
\end{equation}
\begin{equation}
w(u) = \sum_{v\in N(u)}w(u,v). \label{w(u)definition}
\end{equation}
\begin{equation}
\sum_{v\in N(u)}w(u,v)\geq \sum_{v\in N(u)}\frac{1}{d(u)}= 1. \label{mindegeq2}
\end{equation}
\begin{equation}
\sum_{v\in N(u)}w(u,v) \leq \sum_{v\in N(u)}1 = d(u). \label{mindegeq5}
\end{equation}

\begin{equation}
 \sum_{(u,v) \in E}\frac{1}{d(u)} + \frac{1}{d(v)} =\sum_{u \in V}d(u)\left(\frac{1}{d(u)}\right) = n. \label{mindegeq3}
\end{equation}

Equation \eqref{w(u)definition} is directly from the definition in section \ref{weighted graphs}, equation \eqref{c(u)}. From equation \eqref{c(G)} (and noting that we are using $w$ to stand in for $c$), 
\[
w(G)=\sum_{u \in V}w(u)= 2\sum_{(u,v)\in E}w(u,v),
\]
so summing over all edges and using (\ref{mindegeq1}),
\begin{equation}
w(G)=2\sum_{(u,v)\in E}w(u,v) \leq 2\sum_{(u,v) \in E}\frac{1}{d(u)} + \frac{1}{d(v)} \leq 4\sum_{(u,v)\in E}w(u,v) = 2w(G) \label{mindegeq4}
\end{equation}
Then using (\ref{mindegeq3}) in the middle part of (\ref{mindegeq4}),
\begin{equation}
n\leq w(G) \leq 2n. \label{jdnakdjn}
\end{equation}

Thus, the min-deg weighting scheme has the following transition probabilities:

\begin{equation}
 P_{u,v} = 
  \begin{cases}
   \frac{1/\min\{d(u), d(v)\}}{\sum_{w\in N(u)}1/\min\{d(u), d(w)\}} & \text{if } v \in N(u) \\
   0       & \text{otherwise} 
  \end{cases}\label{mindegscheme}
\end{equation}

\subsection{Hitting Time}
The proof of $O(n^2)$ hitting time in \cite{ikeda2} applies the following, generally useful lemma.
\\

\begin{lemma}\label{sumdegpath}
For any connected graph $G$ and any pair of vertices $u,v \in V(G)$, let $\rho=(x_0, x_1, \ldots, x_{\ell})$ where $x_0=u$ and $x_{\ell}=v$ be a shortest path between $u$ and $v$. Then 
\[
\sum_{i=0}^{\ell}d(x_i) \leq 3n
\]
where $d(x)$ is the degree of vertex $x$ and $n=|V(G)|$.
\end{lemma}
\proofstart
For any $0\leq i<i+2<j\leq \ell$, it is the case that $(x_i, x_j)\notin E(G)$ and $N(x_i)\cap N(x_j) = \emptyset$. If not, then there would be a shorter path. Therefore, each vertex in $G$ can be connected to at most $3$ vertices on $\rho$, and the lemma follows.
\proofend

We have the following 
\\

\begin{theorem}\label{weightedhittingtimeTHM}
For a graph $G$ under the min-deg weighting scheme, $\hit[u,v] \leq 6n^2$ for any pair of vertices $u,v\in V(G)$.
\end{theorem}
\proofstart
Let $\rho=(x_0, x_1, \ldots, x_{\ell})$ where $x_0=u$ and $x_{\ell}=v$ be a shortest path between $u$ and $v$. 
\begin{eqnarray}
\hit[u,v]&\leq& \sum_{i=0}^{\ell-1}\hit[x_i, x_{i+1}]\\ 
&\leq& \sum_{i=0}^{\ell-1}\com[x_i, x_{i+1}]\\ \label{kjnkdnjs}
&=&w(G)\sum_{i=0}^{\ell-1}R(x_i, x_{i+1})\label{oicmlv}
\end{eqnarray}
where $R(x,y)$ is the effective resistance between vertices $x$ and $y$ and \eqref{oicmlv} follows \eqref{kjnkdnjs} by Theorem \ref{Th:commtetimetheorem}. 

Now
\[
R(x,y)\leq r(x,y) = \min\{d(x), d(y)\}
\]
and so
\[
\sum_{i=0}^{\ell-1}R(x_i, x_{i+1})\leq \sum_{i=0}^{\ell-1} \min\{d(x_i), d(x_{i+1})\}\leq \sum_{i=0}^{\ell-1}d(x_i) \leq 3n,
\]
where the last inequality follows by Lemma \ref{sumdegpath}.

By (\ref{jdnakdjn}) we have $w(G)\leq 2n$, and the theorem follows.
\proofend

By Matthews' technique (Theorem \ref{MatLemma1}), 
\\

\begin{corollary}\label{weightedhittingtimeTHMCor}
$\cov[G]=O(n^2\log n)$.
\end{corollary}

\cite{ikeda2} does not use electrical network theory; the proof that the Ikeda weighting scheme results in an $O(n^2)$ upper bound on hitting time is different. Furthermore, \cite{ikeda2} conjectures that this weighting scheme in fact gives an $O(n^2)$ upper bound on cover time. We do the same for our weighting scheme:
\\

\begin{conjecture}\label{O(n^2)conjecture}
The minimum degree weighting scheme has $O(n^2)$ cover time on all graphs $G$.
\end{conjecture}

To our knowledge, no weighting scheme has been shown to meet an $O(n^2)$ on all graphs $G$ (where, as stated in the introduction to the chapter, $G$ is assumed to be simple, connected and undirected).

\section{Random Graphs of a Given Degree Sequence}
In chapter \ref{CovDS} we determined the following theorem:

\textbf{Theorem \ref{CovTh} (restated)}\emph{ 
Let  \(G\) be chosen \uar\ from \(\mathcal{G}(\textbf{d})\), where \textbf{d} is nice.
Then \whp,
\begin{equation}
\cov[G] \sim \frac{d-1}{d-2}\frac{\theta}{d}\;n\log n, \label{unweightedCovertime}
\end{equation}
where $d$ is the effective minimum degree and $\theta$ is the average degree.
}

Recall $\omega=\log\log\log n$ as defined in equation \eqref{om}. In this section we prove the following:
\\

\begin{theorem}\label{weightedupperboundcovertime}
Let  \(G\) be chosen \uar\ from \(\mathcal{G}(\textbf{d})\), where \textbf{d} is nice with the following extra restriction: $\Delta\leq \omega^{\frac{1}{4}}$ \ignore{and $\kappa<1/25$ (rather than $\kappa<1/11$, see section \ref{seca})}. Weight the edges of $G$ with the min-deg weighting scheme, that is, for an edge $(u,v)$, assign it weight $w(u,v) = 1/\min\{d(u), d(v)\}$. Denote the resulting graph $G_w$
Then \whp,
\begin{equation}
\cov[G_w] \leq (1+o(1))\frac{\d-1}{\d-2}\;8n\log n. \label{mindegcovtime}
\end{equation}
where $\delta$ is the minimum degree.
\end{theorem}

Note that the degree sequence assumptions of section \ref{seca} allow for the ratio $\theta/d$ to be unbounded. Hence, the ratio of the min-deg cover time to the simple cover time, that is, the \emph{speed up}, can be unbounded.

Our approach in this section is similar to that in the chapter \ref{CovDS}, and we borrow heavily from there.  We sometimes use $d_v$ for $d(v)$.

\subsection{Conductance}
We restate the definition of conductance from section \ref{conductancesectionprelims}: $Q(u,v) = \pi(u)P_{u,v}$, and for a set $S \subset V$, $\pi(S) = \sum_{u \in S}\pi(u)$, $S' = V \setminus S$ and $Q(S,S') = \sum_{u \in S, v \in S'}Q(u,v)$. The conductance $\Phi(G_w)$ of  $G_w$ is defined as 
\[
\Phi(G_w) =   \min_{\pi(S)\leq 1/2}\frac{Q(S,S')}{\pi(S)}.
\]
We show
\\

\begin{lemma}
For $G \in \mathcal{G}(\mathbf{d})$ where $\mathbf{d}$ is nice, $\Phi(G_w)\geq 1/(100\Delta)$ \whp.
\end{lemma}
\proofstart

Since $w(u) = \sum_{v \in V} w(u,v)$ (where $w(u,v)=0$ if $(u,v) \notin E$) and $\pi(u) = \frac{w(u)}{w(G_w)}$ and $P_{u,v} = \frac{w(u,v)}{w(u)}$ we have $Q(u,v) =  \pi(u)P_{u,v} = \frac{w(u,v)}{w(G_w)}$ and 
\[
Q(S,S') = \frac{1}{w(G_w)}\sum_{u \in S, v \in S'}w(u,v).
\]
Since 
\[
\pi(S) = \sum_{u \in S}\pi(u) = \frac{1}{w(G_w)}\sum_{u \in S}w(u),
\]
we have
\begin{equation}
\Phi(G_w) =   \min_{\pi(S)\leq 1/2}\frac{\sum_{u \in S, v \in S'}w(u,v)}{\sum_{u \in S}w(u)}.\label{lampofsjg}
\end{equation}

In section \ref{ConductanceLowerbound} we showed that \whp\, the conductance is bounded below by $\varepsilon = 1/100$: For a graph $G$ picked \uar\ from $\mathcal{G}(\mathbf{d})$, subject to $\mathbf{d}$ being nice,, \whp\
\begin{equation}
\mathcal{E}(S) = \frac{|E(S:\overline S)|}{d(S)} \geq \varepsilon \label{Eratio}
\end{equation}
for any set $S$ such that $\pi(S) \leq 1/2$ (recall $E(S:\overline S)$ is the set of edges with one end in $S$ and the other in $\overline S$, and $d(S) = \sum_{v\in S}d(v)$). This implied that for an unweighted (or uniformly weighted) graph, $\Phi(G) \geq \varepsilon$, equation (\ref{lampofsjg}) becomes
\[
\Phi(G) =   \min_{\pi(S)\leq 1/2}\frac{|E(S:\overline S)|}{d(S)}.
\]
If $\Delta$ is the maximum degree in $\mathbf{d}$, then $w(e)\geq 1/\Delta$ for any edge $e$. Therefore, $\Phi(G_w)\geq \Phi(G)/\Delta \geq 1/(100\Delta)$.
\proofend

We set 
\begin{equation}
T=\omega^2\log n \label{newTforweighted}.
\end{equation}
Then,
\\
 
\begin{corollary}
Given the value of $T$ from \eqref{newTforweighted} and the restriction on $\Delta$, by section \ref{MixingTimeSection}, equation \eqref{mix}, $t\geq T$ satisfies \eqref{4}.
\end{corollary}
\proofstart
Section \ref{MixingTimeSection}, equation \eqref{mix} restated:
\begin{eqnarray*}
|P_{u}^{(t)}(x)-\pi_x| &\leq& \left(\frac{\p_x}{\p_u}\right)^{1/2}\left(1-\frac{\Phi^2}{2}\right)^t\\
&\leq&\Delta^{1/2}\left(1-\frac{1}{K\Delta^2}\right)^T\\
&\leq&\omega^{1/8}\exp\left(-\frac{\omega^2\log n}{K\omega^{1/2}}\right)\\
&\leq& n^{-3}
\end{eqnarray*}
where $K$ is a constant.
\proofend

\subsection{The Stationary Distribution}
\begin{lemma}
For a vertex $u$, 
\begin{equation}
\frac{1}{2n} \leq \pi(u) \leq \frac{d(u)}{n}.\label{minddegstatdist}
\end{equation}
\end{lemma}
\proofstart
$\pi(u) = \frac{w(u)}{w(G)}$. Now use \eqref{mindegeq2} and \eqref{mindegeq5} with \eqref{jdnakdjn}.
\proofend

\ignore{Observe, given the value of $T$ from \eqref{newTforweighted} and the restriction on $\Delta$, we have}

\begin{corollary}
$T\pi(u)=o(1)$.
\end{corollary}

\subsection{The Number of Returns in the Mixing Time}\label{The number of returns in the mixing time}
Recall the definitions in section \ref{props} and the results in section \ref{Rvconds}. We first require a version of Lemma \ref{R*Lemma} for weighted $\Gamma_v$. We only require the case that it has no cycles.
\\

\begin{lemma}\label{weightedR*Lemma}
Assume $G_w$ is good. Let $\cW_v^*$ denote the walk on $\Gamma_v$ (weighted) starting at $v$ with $\Gamma_v^{\circ}$ made into absorbing states. Assume further that there are no cycles in $\Gamma_v^{\circ}$. Let $R_v^* = \sum_{t=0}^{\infty} r^*_t$ where $r^*_t$ is the probability that $\cW_v^*$ is at vertex $v$ at time $t$. 
\begin{equation*}
R_v = R_v^* + O(\sqrt{\omega} e^{-\Omega(\sqrt{\omega})}).
\end{equation*}
\end{lemma}
\proofstart
We use the proof of Lemma \ref{R*Lemma} to prove equation \eqref{A2}, reproduced here
\begin{equation*}
R_v - R_v^* =\left(\sum_{t=0}^{\om} r_t-r^*_t\right)+ \left(\sum_{t=\om+1}^{T} r_t-r^*_t\right) -
\sum_{t=T+1}^{\infty}r^*_t .
\end{equation*}
Firstly, observe that since $\Gamma_v^{\circ}$ has no cycles and no heavy vertices (by the restriction on $\Delta$), it is a tree with leaves at distance $\omega$. 

\underline{{\em Case \(t \le \om\).}}
For \(t\leq \om\), \(r^*_t=r_t\). Thus we can write
\begin{equation}
 \left(\sum_{t=0}^{\om} r_t-r^*_t\right)=0.
\end{equation}

\underline{{\em Case \(\om+1 \le  t \le T\).}}
We use (\ref{mix}) with \(x=u=v\) and observe \(z=(1-\Phi(G_w)^2/2)<(1-1/K\sqrt{\omega})\) for some constant $K$. Observe $T\pi_v=O(1/n^c)$ for some constant $c>0$, so $T\pi_v=o(\sqrt{\omega}e^{-\sqrt{\omega}})$. Hence,
\begin{equation}
\sum_{t=\omega+1}^{T} |r_t-r^*_t| = \sum_{t=\omega+1}^{T} r_t
\leq \sum_{t=\omega+1}^{T} (\pi_v + z^t )\leq T\pi_v+ \frac{z^{\omega}}{1-z} = O(\sqrt{\omega}e^{-\Omega(\sqrt{\omega})}).
\end{equation}

\underline{{\em Case \(t \ge T+1\).}}
The proof of Lemma \ref{R*Lemma} for this case up to \eqref{65g46sg4s6dg} applies, with the final equality in \eqref{65g46sg4s6dg} replaced with 
\[
O(\omega^2/n^{\Omega(1)}) = O(\sqrt{\omega}e^{-\Omega(\sqrt{\omega})}).
\]
\proofend

\begin{lemma}\label{R_vweightedscheme}
Subject to the extra restrictions on $\mathbf{d}$, for a min-deg weighted good $G \in \mathcal{G}(\mathbf{d})$ and a vertex $v \in V(G)$,
\begin{description}
\item[(a)] if $u$ is locally tree-like then $R_u = \frac{\delta-1}{\delta-2}+O(\sqrt{\omega} e^{-\Omega(\sqrt{\omega})})$.
\item[(b)] $R_u \leq \frac{\delta}{\delta-2}\frac{\delta-1}{\delta-2}+O(\sqrt{\omega} e^{-\Omega(\sqrt{\omega})})\leq 6+O(\sqrt{\omega} e^{-\Omega(\sqrt{\omega})})$.
\end{description}
\end{lemma}
\proofstart
\textbf{(a)}
Suppose $u$ is locally tree-like, but not necessarily regular. Let $G(u)$ be the locally (i.e., out to distance $\omega = \log\log\log n$) induced subgraph of $G_w$ and let $v \in G(v)$. Let $p$ be the probability of movement away from $u$ and $q$ the probability toward it. Observe the weighting scheme gives
\[
\frac{p}{q}\geq \frac{(d_v-1)1/d_v}{1/\delta} \geq \delta\frac{\delta-1}{\delta}  = \delta-1.
\]
This ratio is the same as an unweighted walk on a $\delta$-regular vertex $u$. Therefore, using Lemma \ref{weightedR*Lemma} and by the same principles as Lemma \ref{RvRltl}
\[
R_u \leq \frac{\delta-1}{\delta-2} +O(\sqrt{\omega} e^{-\Omega(\sqrt{\omega})}).
\]

\textbf{(b)} Since any cycle in $\Gamma_v^{\circ}$ is unique, at most two of the edges, $(u,v_1), (u,v_2)$ out of $u$ lead to vertices on a cycle. Let $bad$ be the event event of moving from $u$ to some $v \in \{v_1,v_2\}$  Then
\[
\Pr(bad) \leq \frac{2/\delta}{2/\delta+(d_u-2)(1/d_u)} = \frac{2/\delta}{2/\delta+1-2/d_u}\leq 2/\delta.
\]

The probability of moving from $u$ to some $v \in N(u)\setminus\{v_1,v_2\}$ is then at least $\frac{\delta-2}{\delta}$, implying $\frac{\delta}{\delta-2}$ returns to $u$ in expectation for every transition from $u$ to $N(u)\setminus\{v_1,v_2\}$. Assuming that a move from $u$ to $v_1$ or $v_2$ always results in an immediate return, we can bound $R_u \leq \frac{\delta}{\delta-2}(\frac{\delta-1}{\delta-2}+O(\sqrt{\omega} e^{-\Omega(\sqrt{\omega})}))\leq 6+O(\sqrt{\omega} e^{-\Omega(\sqrt{\omega})})$.
\proofend

\subsection{The Number of Vertices not Locally Tree-like}
We wish to bound the number of vertices $v$ that are not locally tree-like, i.e., for which $\Gamma_v$ has a cycle.
We use the following from the proof of Lemma \ref{lightCPC} of section \ref{Structural Properties}:
$\mathcal{F}(2x) = \frac{(2x)!}{2^xx!}$, then

\begin{equation*}
\frac{\mathcal{F}(\theta n - 2x)}{\mathcal{F}(\theta n)} \leq \left(\frac{1}{\theta n-2x}\right)^x
\end{equation*}

\begin{lemma}\label{weightednonlocallytreelike}
With probability at least $1-n^{-\Omega(1)}$, the number of vertices not locally tree-like is at most $n^{1/10}$ .
\end{lemma}
\proofstart
The expected number of small cycles has upper bound
\begin{eqnarray*}
&&\sum_{k=3}^{2\omega+1}\binom{n}{k}\frac{(k-1)!}{2}(\Delta(\Delta-1))^k\frac{\mathcal{F}(\theta n - 2k)}{\mathcal{F}(\theta n)}\\
&&\leq\sum_{k=3}^{2\omega+1}n^k\Delta^{2k}\left(\frac{1}{\theta n-2k}\right)^k\\
&&\leq\sum_{k=3}^{2\omega+1}\Delta^{2k}\left(\frac{n}{\theta n-4\omega-2}\right)^k\\
&&\leq \Delta^{4(\omega+1)}
\end{eqnarray*}
Therefore, the expected number of vertices within distance $\omega$ of a cycle is at most $\Delta^{4(\omega+1)}\Delta^{\omega}=\Delta^{5\omega+4}$. Since $\Delta\leq \omega^{1/4}$, the lemma follows by Markov's inequality.
\proofend

\subsection{The Probability a Vertex is Unvisited}
We use the same principles as in section \ref{ubsec}. Restating \eqref{shed11},
\begin{equation}
\cov_u[G] \leq  t+1+\sum_v \sum_{s \ge t}\Pr(\ul A_s(v)).\label{shed11againgfsdgs}
\end{equation}

We use Lemmas \ref{Eval} and \ref{crude}, which hold for weighted random walks (see chapter 2, General Markov Chains, in \cite{AlFi} for justification of \eqref{Z_wwEquation} and the inequality $D(s+t) \le 2D(s)D(t)$. All other expressions in the proofs hold for weighted random walks). Thus,
\[
\Pr({\ul A_t(v)}) \le  \exp\brac{\frac{-(1-o(1))\lfloor t/\tau_v\rfloor}{2}},
\]
where $\tau_v=T+2\hit[\pi, v]$ and $\hit[\pi,v] = \ooi R_v/\pi_v$.

Hence, for a given $v$,
\begin{eqnarray*}
\sum_{s \ge t}\Pr(\ul A_s(v)) &\leq& \sum_{s \ge t}\exp\brac{\frac{-(1-o(1))\lfloor s/\tau_v\rfloor}{2}}\\
&\leq&\tau_v\sum_{s \ge \lfloor t/\tau_v\rfloor}\exp\brac{\frac{-(1-o(1))s}{2}}\\
&\leq&3\tau_v\exp\brac{\frac{-(1-o(1))\lfloor t/\tau_v\rfloor}{2}}\\
&=&3\tau_v\exp\brac{\frac{-(1-o(1))}{2}\left\lfloor \frac{t\pi_v}{T\pi_v+\ooi 2R_v}\right\rfloor}.
\end{eqnarray*}
Since $T\pi_v=o(1)$ and $\pi_v \geq 1/2n$ from \eqref{minddegstatdist}, we get
\begin{equation*}
\sum_{s \ge t}\Pr(\ul A_s(v)) \leq 3\tau_v\exp\brac{\frac{-(1-o(1))}{2}\left\lfloor \frac{t}{\ooi 4nR_v}\right\rfloor}
\end{equation*}

Let $t=t^*=(1+\epsilon)8\frac{\delta-1}{\delta-2}n\log n$ where $\epsilon \rightarrow 0$ sufficiently slowly. Then
\begin{equation}
\sum_{s \ge t}\Pr(\ul A_s(v)) \leq 3\tau_v\exp\brac{-(1+\Theta(\epsilon))\frac{\delta-1}{\delta-2}\frac{\log n}{R_v}}\label{g6s8g466846847}
\end{equation}

We partition the double sum $\sum_v \sum_{s \ge t}\Pr(\ul A_s(v))$ from \eqref{shed11againgfsdgs} into 
\begin{equation*}
\sum_{v \in V_A}\sum_{s \ge t}\Pr(\ul A_s(v))+\sum_{v \in V_B} \sum_{s \ge t}\Pr(\ul A_s(v))
\end{equation*}
where $V_A$ are locally tree-like and $V_B$ are not.

If $v$ is locally tree-like, then using Theorem \ref{R_vweightedscheme} \textbf{(a)}, the RHS of \eqref{g6s8g466846847} is bounded by 
\begin{eqnarray*}
3\tau_v n^{-(1+\Theta(\epsilon))} &=& 3(T+2\ooi R_v/\pi_v) n^{-(1+\Theta(\epsilon))}\\
&\leq& \ooi 12nR_vn^{-(1+\Theta(\epsilon))}\\
&=&O(1)n^{-\Theta(\epsilon)}
\end{eqnarray*}
Thus, 
\begin{equation}
\sum_{v \in V_A}\sum_{s \ge t}\Pr(\ul A_s(v)) \leq O(1)n^{1-\Theta(\epsilon)}=o(t). \label{V_Ao(t)}
\end{equation}

For any $v$ (i.e., including those not locally tree-like),  \eqref{g6s8g466846847} is bounded by
\begin{equation}
3\tau_v n^{-(1+\Theta(\epsilon))\frac{\delta-1}{6(\delta-2)}} \leq O(1) n^{1-(1+\Theta(\epsilon))\frac{\delta-1}{6(\delta-2)}}\label{fjaknqfafaf456}
\end{equation}

Using Lemma \ref{weightednonlocallytreelike} to sum the bound \eqref{fjaknqfafaf456} over all non locally tree-like vertices, we get 
\begin{equation}
\sum_{v \in V_B}\sum_{s \ge t}\Pr(\ul A_s(v)) \leq O(1) n^{\frac{1}{10}+1-(1+\Theta(\epsilon))\frac{\delta-1}{6(\delta-2)}}=O(n^{\frac{1}{2}})=o(t).\label{V_Bo(t)}
\end{equation}

Hence, combining \eqref{shed11againgfsdgs}, \eqref{V_Ao(t)} and \eqref{V_Bo(t)} for $t=t^*$, Theorem \ref{weightedupperboundcovertime} follows.
\proofend

Compare this with \eqref{unweightedCovertime}, we see that the speed up,
\[
\mathcal{S}=\frac{\cov[G]}{\cov[G_w]} = \Omega(\theta),
\]

Therefore $\mathcal{S}\rightarrow \infty$ as $n \rightarrow \infty$ if $\theta\rightarrow \infty$ as $n \rightarrow \infty$. That is, we can have an unbounded speed up.

We conjecture that the following tighter bound holds:
\begin{conjecture}\label{tightcovDSweightedconjecture}
Equation \eqref{mindegcovtime} can be replaced by 
\[
\cov[G_w] \leq (1+o(1))\frac{d-1}{d-2}\;n\log n.
\]
\end{conjecture}

\chapter{Conclusion}
\section{Main Results}
In this thesis, we have addressed the following question, first posed in the Introduction: 

\emph{For a random walk $\mathcal{W}_u$ on a simple, connected, undirected graph $G=(V,E)$, what is the expected number of steps required
to visit all the vertices in $G$, maximised over starting vertices $u$?}

We have given answers to this question for three cases: 
{\begin{enumerate}
	\item Simple random walks on graphs $G$ of the form $G=G_1\Box G_2$ where $\Box$ is the Cartesian product operator,
	\item Simple random walks on graphs $G$ picked \uar\ from the space of graphs $\mathcal{G}(\mathbf{d})$ where $\mathbf{d}$ is a given degree sequence satisfying various constraints,
	\item Weighted random walks for general graphs, and graphs $G$ picked \uar\ from $\mathcal{G}(\mathbf{d})$.
\end{enumerate}}
In particular, we have proved the following theorems \footnote{Theorems, lemmas and conjectures presented in this chapter are restatements of ones presented previously.}.
 
For the first case, we showed in chapter \ref{ch:The Cover Time of Cartesian Product Graphs}:

\textbf{Theorem \ref{main}}\emph{
Let $F= (V_F, E_F) = G \Box H$ where $G=(V_G, E_G)$ and $H = (V_H, E_H)$ are simple, connected, unweighted, undirected graphs. 
We have
\begin{equation*}
\cov[F] \geq \max\left\{\left(1+\frac{\delta_G}{\Delta_H}\right) \cov[H], \left(1+\frac{\delta_H}{\Delta_G}\right) \cov[G]\right\}. 
\end{equation*}
Suppose further that $n_H \geq D_G+1$, then
\begin{equation*}
\cov[F] \leq K\left(\left(1+\frac{\Delta_G}{\delta_H}\right)\bcov[H] + \frac{Mm_Gm_Hn_H\ell^2}{\cov[H]D_G}\right)
\end{equation*}
where $M =|E_F|= n_Gm_H + n_Hm_G$, $\ell=\log (D_G+1) \log(n_GD_G)$ and $K$ is some universal constant.}

(We remind the reader that in the above, for a graph $\Gamma$, $n_{\Gamma}=|V({\Gamma})|$, $m_{\Gamma}=|E({\Gamma})|$, $D_{\Gamma}$ is the diameter of ${\Gamma}$, and $\delta_{\Gamma}$ and  $\Delta_{\Gamma}$ are the minimum and maximum degrees respectively). 

For the second case, we showed in chapter \ref{CovDS}:

\textbf{Theorem \ref{CovTh}}\emph{ 
Let  \(G\) be chosen \uar\ from \(\mathcal{G}(\textbf{d})\), where \textbf{d} is nice.
Then \whp,
\begin{equation*}
\cov[G] \sim \frac{d-1}{d-2}\frac{\theta}{d}\;n\log n, 
\end{equation*}
where $d$ is the effective minimum degree and $\theta$ is the average degree.
}

Theorem \ref{main} answers the question within bounds, whilst Theorem \ref{CovTh} gives a precise asymptotic value. We reiterate some definitions: a statement $\mathcal{P}(n)$ parameterised on an integer $n$ holds \emph{with high probability} (\whp) if $\Pr(\mathcal{P}(n) \text{ is true})\rightarrow 1$ as $n \rightarrow \infty$; the notation $f(n) \sim g(n)$ means $f(n)/g(n) \rightarrow 1$ as $n \rightarrow \infty$.

The third case was addressed in chapter \ref{Weighted Random Walks}. Each edge $e$ is given a non-negative weight $w(e)$, and the probability of transitioning $e$ from vertex $u$ is $w(e)/w(u)$ where $w(u)=\sum_{e=(u,v)\in E}w(e)$. We presented a simple weighting scheme, the \emph{minimum degree} (min-deg) weighting scheme, where each edge $e=(u,v)$ has $w(e)=1/\min\{d(u), d(v)\}$. That is, the transition probability matrix was

\begin{equation*}
 P_{u,v} = 
  \begin{cases}
   \frac{1/\min\{d(u), d(v)\}}{\sum_{w\in N(u)}1/\min\{d(u), d(w)\}} & \text{if } v \in N(u) \\
   0       & \text{otherwise} 
  \end{cases}
\end{equation*} 

We proved the following theorem for general graphs ($G$ is assumed to be simple, connected and undirected):

\textbf{Theorem \ref{weightedhittingtimeTHM}}
\emph{For a graph $G$ under the min-deg weighting scheme, $\hit[u,v] \leq 6n^2$ for any pair of vertices $u,v\in V(G)$.}

Consequently, by Matthews' technique (Theorem \ref{MatLemma1}),

\textbf{Corollary \ref{weightedhittingtimeTHMCor}}
\emph{$\cov[G]=O(n^2\log n)$.}

We then studied this weighted walk on the specific class of graphs that is the subject of Theorem \ref{CovTh}, applying the min-deg weighting scheme to $G \in \mathcal{G}(\mathbf{d})$. We presented the following theorem:

\textbf{Theorem \ref{weightedupperboundcovertime}}\emph{ 
Let  \(G\) be chosen \uar\ from \(\mathcal{G}(\textbf{d})\), where \textbf{d} is nice with the following extra restriction: $\Delta\leq \omega^{\frac{1}{4}}$ \ignore{and $\kappa<1/25$ (rather than $\kappa<1/11$, see section \ref{seca})}. Weight the edges of $G$ with the min-deg weighting scheme, that is, for an edge $(u,v)$, assign it weight $w(u,v) = 1/\min\{d(u), d(v)\}$. Denote the resulting graph $G_w$
Then \whp,
\begin{equation}
\cov[G_w] \leq (1+o(1))\frac{\d-1}{\d-2}\;8n\log n. 
\end{equation}
where $\delta$ is the minimum degree.
}

We have, by the above, extended the body of knowledge on random walks on graphs, in particular, the cover time of random walks on graphs. 

\section{Secondary Results}
Additionally, in the course of proving the above theorems, we have presented a number of results which may be of independent interest. Specifically, we have given results related to effective resistance of the Cartesian product of two graphs, a frame work for analysing walks on a graph by analysing ``local observations'' of the walk, and structural results for \(G \in \mathcal{G}(\textbf{d})\).

The following lemmas were proved in chapter \ref{ch:The Cover Time of Cartesian Product Graphs}:

\textbf{Lemma \ref{lg1x}}\emph{  
For a graph $G$ and tree $T$, $R_{max}(G\Box T) < 4R_{max}(G\Box P_r)$ 
where $|V(T)| \leq r \leq 2|V(T)|$ and $P_r$ is the path on $r$ vertices.
}

\textbf{Lemma \ref{lg2x}}\emph{
For graphs $G,H$ suppose $D_G+1 \leq n_H \leq \alpha (D_G+1)$, for some $\alpha$. Then $R_{max}(G \Box H) < \zeta \alpha \log(D_G+1)$, where $\zeta$ is some universal constant.
}

(We remind the reader that $R_{max}(G)$ is the maximum effective resistance between a pair of vertices in a graph $G$).

In chapter \ref{ch:The Cover Time of Cartesian Product Graphs}, to prove the Theorem \ref{main}, we presented a framework to bound the cover time of a random walk $\mathcal{W}(G)$ on a graph $G$ by dividing the graph up into regions - sets of vertices $S \subseteq V(G)$ - and analysing the behaviour of the walk whilst on those regions. It does so by relating $\mathcal{W}(G)$ whilst walking on $S$ to a random walk $\mathcal{W}(H)$ on $H$, where $H$ is weighted graph derived from $(G,S)$. The analysis of the $\mathcal{W}(H)$ walks can then be composed over the whole graph. Thus the analysis of the whole graph is reduced to the analysis of outcomes on local regions and subsequent compositions of those outcomes. This framework can be applied more generally than Cartesian products. The details are presented in section \ref{DefnSubgraphWalk}.

In chapter \ref{CovDS}, we presented the following result on the conductance of a graph $G \in \mathcal{G}(\mathbf{d})$:

\textbf{Theorem \ref{condLowerboundTheorem}}\emph{
Subject to assumptions \textbf{(i)}-\textbf{(iv)} listed in section \ref{ConductanceLowerbound}, for a graph $G \in \mathcal{G}(\mathbf{d})$, $\Phi(G)>1/100$ \whp.
}

\section{Future Work}
We discuss some possible avenues for further study. 

\subsection{Cover Time of other Random Graph Models}
Similar techniques to those in chapter  \ref{CovDS} have been used to study the cover times of other models of random graphs, for example, the preferential attachment graph \cite{CGweb}, the random geometric graph \cite{ColinGeometric} and the giant component of the Erd\H{o}s--R\'{e}nyi random graph \cite{CooperFriezeGiant} (see section \ref{Random Graphs: Models and Cover Time} for a discussion of these models). We envisage that the techniques can also be extended to other models of random graphs, for example, the random intersection graph \cite{blackburn}, \cite{stark}, the Bollob\'{a}s--Chung (B--C) model \cite{bcmodel}, and the Watts-Strogatz (W--S) model \cite{wattsstrogatz}\footnote{The  Watts-Strogatz model is motivated by the desire to find graph models whose properties more accurately reflect ``real-world'' networks such as social networks. The W-S model has some of these desired properties, such as small distance between nodes (hence the term \emph{small world} model). See \cite{wattsstrogatz}. The Bollob\'{a}s--Chung model was motivated by the desire to create bounded-degree graphs with small diameter. \cite{bcmodel} demonstrates approximately $\log_2 n$ diameter for the B--C model.}. 

We quote \cite{blackburn} for the definition of the (uniform) random intersection graph\footnote{\cite{blackburn} motivates random intersection graphs by wireless networks, where a colour is a unique cryptographic key, and a pair of sensors can communicated if and only if they have a key in common.}:
\begin{quote}
The uniform random intersection graph $G(n, m, k)$ is a random graph defined as follows. Let $V$ be a set of n nodes, and let $M$ be a set of $m$ colours.
To each node $v \in V$ we assign a subset $F_v \subseteq M$ of $k$ distinct colours, chosen
uniformly and independently at random from the $k$-subsets of $M$. We join
distinct nodes $u, v \in V$ by an edge if and only if $F_u \cap F_v \neq \emptyset$.
\end{quote}

The B--C model is a $n$-cycle with a random matching. In the W--S model, the graph begins as an $n$ cycle with each vertex $v$ connected to $k$ others; the $k/2$ nearest neigbours on each side. Labelling the vertices round the cycle $(0,1,\ldots,n-1)$, for each vertex $i$ in turn, ``rewire'' each edge $(i,j)$, $i<j$ to $(i,k)$ where $k$ is chosen \uar\ from amongst the vertices that would avoid loops or parallel edges. 

As we envisage it, the crux of a proof, as in chapter \ref{CovDS}, would be the application of Lemma \ref{L3}. This would require verification that the conditions of the lemma hold, determining a suitable mixing time $T$, degree distribution of the vertices, and calculations for the expected number of returns $R_v$ to a vertex $v$ within the mixing time. 

We envisage that the aforementioned models will be rapidly mixing (i.e., $T=O(\log^an)$ for some constant $a$), and that, for sufficiently sparse graphs, it will be possible to determine degree distributions and demonstrate tree-like local graph structures (which will facilitate  the calculation of $R_v$). Indeed, for the B--C model, each vertex has degree $3$ by construction, and \cite{Durrett} shows $T=O(\log n)$, thus satisfying condition \textbf{(b)} of Lemma \ref{L3}.

\subsection{Weighted Random Walks}
The study of weighted random walks, in particular, the cover time, is not nearly as developed as for the unweighted case. We reiterate the following conjecture, made in chapter \ref{Weighted Random Walks}:

\textbf{Conjecture \ref{O(n^2)conjecture}}\emph{
The minimum degree weighting scheme has $O(n^2)$ cover time on all graphs $G$.
}

We do not believe that ``standard'' techniques (for examples, those presented in chapter \ref{Techniques and Results for Hitting and Cover Times}, and subsequent chapters) are sufficient to prove this conjecture. In fact, we believe that it may require new techniques or further development of the underlying theory of weighted random walks.

It may not be as demanding, however, to improve the bound of Theorem \ref{weightedupperboundcovertime}, or even prove the following conjecture, made in chapter \ref{Weighted Random Walks}:

\textbf{Conjecture \ref{tightcovDSweightedconjecture}}\emph{
Equation \eqref{mindegcovtime} can be replaced by 
\[
\cov[G_w] \leq (1+o(1))\frac{d-1}{d-2}\;n\log n.
\]
}

\cleardoublepage


\addcontentsline{toc}{chapter}{Bibliography}

\begin{thebibliography}{99}

\bibitem{CovDS} M. Abdullah, C.Cooper, A. Frieze, The cover time of a random graph with a given degree sequence.
\emph{21st International Meeting on Probabilistic, Combinatorial and Asymptotic Methods for the Analysis of Algorithms (AofA)} (2010).

\bibitem{CartProd} M. Abdullah, C.Cooper, T. Radzik, The cover time of Cartesian product graphs.
\emph{21st International Workshop on Combinatorial Algorithms (IWOCA)} (2010).

\bibitem{Albert} R. Albert, H. Jeong, A.-L. Barab\'{a}si, Internet: Diameter of the world-wide web.
\emph{Nature}, 401:130--131 (1999).
          
\bibitem{AlFi} D. Aldous and J. Fill, Reversible Markov Chains and
Random Walks on Graphs, http://stat-www.berkeley.edu/pub/users/aldous/RWG/book.html.

\bibitem{AKLLR} R. Aleliunas, R.M. Karp, R.J. Lipton, L. Lov\'asz and
C. Rackoff, Random walks, universal traversal sequences, and the
complexity of maze problems. {\em Proceedings of the 20th Annual
IEEE Symposium on Foundations of Computer Science}, 218-223 (1979). 

\bibitem{chenAdhoc} C. Avin, G. Ercal, Bounds on the mixing time and partial cover of ad-hoc and sensor networks. \emph{EWSN-05. European Workshop on Wireless Sensor Networks} (2005). 

\bibitem{chenGeometric} C. Avin, G. Ercal, On the cover time and mixing time of random geometric graphs. 
\emph{Theor. Comput. Sci.}, 380(1-2):2–22 (2007).

\bibitem{Barabasi} A.-L. Barab\'{a}si, R. Albert, Emergence of scaling in random networks.
\emph{Science}, 286, 509–512 (1999).

\bibitem{Barabasi2} A.-L. Barab\'{a}si, R. Albert, H. Jeong, Mean-field theory for scale-free random networks.
\emph{Physica A}, 272, 173–187 (1999).

\bibitem{blackburn}S.R. Blackburn, S. Gerke, Connectivity of the uniform random intersection graph. 
\emph{Discrete Mathematics}, 309 5130-5140 (2009).

\bibitem{Boll1}
B. Bollob\'as,  A probabilistic proof of an asymptotic formula for the number of labelled regular graphs. 
\emph{European Journal on Combinatorics}, 1 311-316 (1980). 

\bibitem{Bellabook} B. Bollob\'{a}s, \emph{Random Graphs}. Cambridge University Press, 2nd edition (2001). 

\bibitem{Bollobas} B. Bollob\'{a}s, G. Brightwell, Random walks and electrical resistances in product graphs.
\emph{Discrete Applied Mathematics}, 73 69-79 (1997).

\bibitem{bcmodel} B. Bollob\'{a}s, F. R. K. Chung, The diameter of a cycle plus a random matching. 
\emph{SIAM J. Disc. MATH.}, Vol. 1, No. 3 (1988).

\bibitem{BollPref} B. Bollob\'{a}s, O. Riordan, J. Spencer, G. Tusn\'{a}dy, The degree sequence of a
scale-free random graph process.
\emph{Random Structures and Algorithms}, 18 279-290 (2001).

\bibitem{Broutin} N. Broutin, L. Devroye, N. Fraiman, G. Lugosi, Connectivity threshold for Bluetooth graphs. Preprint (2011).

\bibitem{BC} J. Brown and R. Churchill, {\em Complex Variables and Applications,}
(Sixth Edition)  McGraw-Hill (1996).

\bibitem{Chandra} A.K. Chandra, P. Raghavan, W.L. Ruzzo, R. Smolensky, P. Tiwari, The electrical resistance of a graph captures its commute and cover times. \emph{Computational Complexity}, 6 312-340 (1997) 

\bibitem{ConwayGuy} J.H. Conway, R.K. Guy, 
\emph{The Book of Numbers}, Springer-Verlag (1996).

\bibitem{CFreg}  C. Cooper, A.M. Frieze, The cover time of random
regular graphs, {\em SIAM Journal on Discrete Mathematics}, 18
728-740 (2005).

\bibitem{CooperER}  C. Cooper, A.M. Frieze, The cover time of sparse random graphs.
\emph{Random Structures and Algorithms}, 30 1-16 (2007).

\bibitem{CGweb} C. Cooper and A.M. Frieze, The cover time of the
preferential attachment graph. {\em Journal of Combinatorial Theory Series B},
97 269-290 (2007).

\bibitem{CooperFriezeGiant}
C. Cooper and A. M. Frieze. The cover time of the giant component of a random graph.
{\em Random Structures and Algorithms}, 32 401-439 (2008).

\bibitem{Dolev} S. Dolev, E. Schiller, J. Welch, Random walk for self-stabilizing group communication in ad-hoc networks.
\emph{21st IEEE Symposium on Reliable Distributed Systems}, IEEE Computer Society (2002).

\bibitem{cooperp2p} C. Cooper, R. Klasing, T. Radzik, A randomized algorithm for the joining protocol in dynamic distributed networks. \emph{Theor. Comput. Sci.}, 406(3): 248-262 (2008).

\bibitem{coopermult} C. Cooper, A.M. Frieze, T. Radzik, Multiple random walks in random regular graphs. 
\emph{SIAM J. Discrete Math.}, 23(4): 1738-1761 (2009).

\bibitem{ColinGeometric} C.Cooper, A.M. Frieze, The cover time of random geometric graphs.
\emph{Proceedings of SODA}, 48-57 (2009). 

\bibitem{Coppersmith}
D. Coppersmith, U. Feige, J. Shearer, Random walks on regular and irregular graphs.
{\em SIAM Journal on Discrete Mathematics}, 9(2) 301--308 (1996).

\bibitem{Dembo}A. Dembo, Y. Peres, J. Rosen, O. Zeitouni, Cover times for Brownian motion and random walks in two dimensions.
\emph{Ann. Math.}, 160 433-464 (2004).

\bibitem{deistel} R. Diestel,
{\em Graph Theory}, Springer-Verlag (2010).

\bibitem{Ding} J. Ding, J. Lee, Y. Peres, Cover times, blanket times, and majorizing measures.
\emph{Annals of Mathematics}, (to appear).

\bibitem{Durrett} R. Durrett, 
{\em Random Graph Dynamics}, Cambridge University Press (2007).

\bibitem{DoyleSnell} P.G. Doyle and J.L. Snell, 
\emph{Random Walks and Electrical Networks}, (2006).

\bibitem{ER} P. Erd\H{o}s, A. R\'{e}nyi,  On Random Graphs I.  
\emph{Publ. Math. Debrecen}, 6 (1959).

\bibitem{Faloutsos} M. Faloutsos, P. Faloutsos, C. Faloutsos, On power-law relationships of the Internet topology. 
\emph{Proceedings of the conference on Applications, technologies, architectures, and protocols for computer communication}, (1999).


\bibitem{Feige1} U. Feige, A tight upper bound for the cover time of
random walks on graphs. {\em Random Structures and Algorithms}, 6
 51-54 (1995).

\bibitem{Feige2} U. Feige, A tight lower bound for the cover time of
random walks on graphs. {\em Random Structures and Algorithms}, 6
 433-438 (1995).

\bibitem{Feige3} U. Feige, Collecting coupons on trees, and the cover time of random walks.
{\em Computational Complexity}, 6 341--356. (1996/1997)

\bibitem{Fe} W. Feller, {\em An Introduction to Probability Theory, Volume I}
(Second edition), Wiley (1960).

\bibitem{FR} N. Fountoulakis, B. Reed, The evolution of the mixing rate. (To appear).


\bibitem{Gilbert} E.N. Gilbert, Random Graphs. 
\emph{Annals of Mathematical Statistics}, 30 (1959).

\bibitem{Gkants2003} C. Gkantsidis, M. Mihail, A. Saberi, Conductance and congestion in power law graphs. \emph{ACM SIGMETRICS International Conference on Measurement and Modeling of Computer Systems}, 148-159 (2003).


\bibitem{Gkan} C. Gkantsidis, M. Mihail, A. Saberi, Random walks in peer-to-peer networks. \emph{INFOCOM 2004. Twenty-third Annual
Joint Conference of the IEEE Computer and Communications Societies}, vol. 1 (2004).

\bibitem{Gupta} P. Gupta, P.R. Kumar, Critical power for asymptotic connectivity in wireless networks. \emph{Stochastic Analysis, Control, Optimization and Applications}, Birkha\"user, Boston (1998).

\bibitem{hornmatrix} R.A. Horn, C.R. Johnson, 
\emph{Matrix Analysis}. Cambridge University Press (1989).

\bibitem{numbertheorybook} J.F. Humphreys, M.Y. Prest, 
\emph{Numbers, Groups and Codes}. Cambridge University Press (1990).

\bibitem{ikeda1}  S. Ikeda, I. Kubo, N. Okumoto, M. Yamashita, Impact of local topological information on random walks on finite graphs. \emph{30th International Colloquium on Automata, Languages and Programming (ICALP)}, 1054-1067 (2003).

\bibitem{ikeda2} S. Ikeda, I. Kubo, M. Yamashita, The hitting and the cover times of random walks on finite graphs using local degree information. \emph{Theoretical Computer Science}, 410 94-100 (2009).

\bibitem{Janson} S. Janson, T. Luczak, A. Rucinski, 
\emph{Random Graphs}.  Wiley, New York (2000).

\bibitem{JerrumSinclair} M. Jerrum, A. Sinclair, Approximating the permanent.  
\emph{SIAM Journal on Computing}, 18 1149 - 1178 (1989).

\bibitem{JonassonCover} J. Jonasson, On the cover time for random walks on random graphs. 
\emph{Combinatorics, Probability and Computing}, 7 265-279 (1998).

\bibitem{Jonasson}J. Jonasson, An upper bound on the cover time for powers of graphs. 
\emph{Discrete Mathematics}, 222 181-190 (2000).

\bibitem{LovaszBlanket} J. Kahn, J. H. Kim, L. Lov\'{a}sz, V. H. Vu. The cover time, the blanket time, and the
Matthews bound. \emph{41st Annual Symposium on Foundations of Computer Science (Redondo
Beach, CA, 2000)}, 467–475. IEEE Comput. Soc. Press, Los Alamitos, CA (2000).

\bibitem{Kahn}
J.D. Kahn, N. Linial, N.Nisan, M.E. Saks, On the cover time of random walks on graphs.
{\em Journal of Theoretical Probability}, 2(1):121-128 (1989).

\bibitem{Kleinberg} J. Kleinberg, E. Tardos, \emph{Algorithm Design}, Addison Wesley (2005).

\bibitem{LawlerSokal} G.F. Lawler, A.D. Sokal, Bounds on the $L^2$ spectrum for Markov chains and Markov Processes: a generalization of Cheeger's inequality. \emph{Tran. Amer.Math. Soc.}, 309, 557 - 580 (1988).

\bibitem{LPW} D. A. Levin, Y. Peres, E. L. Wilmer.
{\em Markov Chains and Mixing Times}, AMS Press (2008).

\bibitem{Lovaszcombbook} L. Lov\'{a}sz, Combinatorial Problems and Exercises.
\emph{North Holland; 2 edition} (1993).

\bibitem{LovaszSurvey} L. Lov\'{a}sz, Random walks on graphs: a survey. \emph{Combinatorics, Paul Erdos is Eighty},
 Vol. 2 (ed. D. Mikl\'{o}s, V. T. S\'{o}s, T. Szonyi), János Bolyai Mathematical Society, Budapest, 353-398 (1996).

\bibitem{McW}
B. McKay and N. Wormald, Asymptotic enumeration by degree sequence
of graphs with degrees \(o(n^{1/2})\). {\em Combinatorica}, 11 369-382 (1991).

\bibitem{mitz} M. Mitzenmacher, E. Upfal
{\em Probability and Computing: Randomized Algorithms and Probabilistic Analysis}, Cambridge University Press (2005).

\bibitem{Matthews} P. Matthews, Covering problems for Brownian motion on spheres.
{\em Ann. Prob.}, 16:189-199 (1988). 

\bibitem{Nash-Williams}C. St.J A. Nash-Williams, Random walk and electric currents in
networks. \emph{Proc. Camb. Phil. Soc.}, 55:181-194 (1959).

\bibitem{Norris} J. R. Norris,
{\em Markov Chains}, Cambridge University Press (1998).

\bibitem{Omer} O. Reingold, Undirected ST-connectivity in Log-Space. 
\emph{Electronic Colloquium on Computational Complexity}, 94 (2004)


\bibitem{Penrose} M. D. Penrose, \emph{Random Geometric Graphs}. Oxford University Press (2003)

\bibitem{Servetto} S.D. Servetto, G. Barrenechea, Constrained random walks
on random graphs: Routing Algorithms for Large Scale Wireless
Sensor Networks. \emph{ACM Int. workshop on Wireless sensor networks and applications}, ACMPress (2002).


\bibitem{Sin} A. Sinclair, Improved bounds for mixing rates of Markov chains
and multicommodity flow. {\em Combinatorics, Probability and Computing}, 1 351-370 (1992).

\bibitem{stark} D. Stark, The vertex degree distribution of random intersection
graphs. 
\emph{Random Structures and Algorithms}, 24 249–258 (2004).

\bibitem{wattsstrogatz} D. J. Watts, S. H. Strogatz, Collective dynamics of `small-world' networks. 
\emph{Nature}, 393 (6684): 409–10 (1998).

\bibitem{Wi} H. Wilf, {\em Generatingfunctionology,}  Academic Press (1990).

\bibitem{Winkler} P. Winkler, D. Zuckerman, Multiple cover time.
\emph{Random Structures and Algorithms}, 9 403-411 (1996).

\bibitem{Zu} D. Zuckerman, On the time to traverse all edges of a graph.
{\em Information Processing Letters}, 38 335-337 (1991).

\bibitem{ZuLower} D. Zuckerman, A Technique for Lower Bounding the Cover Time. 
\emph{SIAM J. Discrete Math.}, 5(1): 81-87 (1992).

\end{thebibliography}
\end{document}